\documentclass[a4paper,12pt]{article}

\usepackage{indentfirst,latexsym,bm}

\usepackage{amsmath}
\usepackage{amssymb}
\usepackage[pdftex]{graphicx}
\usepackage{amscd}

\numberwithin{equation}{section}

\title{\Large{\textbf{Ricci Flow with Surgery on Four-manifolds
with Positive Isotropic Curvature}}}
\author{Bing-Long Chen and Xi-Ping Zhu \\[8pt]
Department of Mathematics \\
Zhongshan University \\
Guangzhou, P.R.China}
\date {(Revised version)}
\begin{document}
\maketitle

\centerline{\Large{\textbf{Abstract}}}
 \vskip 0.4cm

  In this paper we
study the Ricci flow on compact four-manifolds with positive
isotropic curvature and with no essential incompressible space
form. Our purpose is two-fold. One is to give a complete proof of
the main theorem of Hamilton in \cite{Ha7}; the other is to extend
some results of Perelman \cite{P1}, \cite{P2} to four-manifolds.
During the proof we have actually provided, parallel to the paper
of the second author with H.-D. Cao \cite{CaZ}, all necessary
details for the part from Section 1 to Section 5 of Perelman's
second paper \cite{P2} on the Ricci flow.

\pagebreak[4] \centerline{\Large{\textbf{1. Introduction}}} \vskip
0.5cm Let $M^{n}$ be a compact $n$-dimensional Riemannian manifold
with metric $g_{ij}(x)$. The Ricci flow is the following evolution
equation
$$\frac{\partial}{\partial t}g_{ij}(x,t)=-2R_{ij}(x,t), \mbox{
for }x\in M \mbox{ and }t>0, \eqno (1.1)$$ with
$g_{ij}(x,0)=g_{ij}(x)$, where $R_{ij}(x,t)$ is the Ricci
curvature tensor of the evolving metric $g_{ij}(x,t)$. This
evolution system was initially introduced by Hamilton in
\cite{Ha1}. Now it has been found to be a powerful tool to
understand the geometry, topology and complex structure of
manifolds (see for example \cite{Ha1}, \cite{Ha2}, \cite{Ha3},
\cite{Ha7}, \cite{Ha8}, \cite{C}, \cite{Hu}, \cite{Cao85}
\cite{CZ}, \cite{CTZ}, \cite{P1}, \cite{P2}, \cite{CaZ} etc.)

 One of the main topics in
modern mathematics is to understand the topology of compact three
dimensional and four dimensional manifolds. The idea to approach
this problem via the Ricci flow is to evolve the initial metric by
the evolution equation (1.1), and try to study the geometries
under the evolution. The key point of this method is to get the
long-time behavior of the solutions of the Ricci flow. For a
compact three (or four) dimensional Riemannian manifold with
positive Ricci curvature (or positive curvature operator,
respectively) as initial data, Hamilton \cite{Ha1} (or \cite{Ha2}
respectively) proved that the solution to the Ricci flow keeps
shrinking and tends to a compact manifold with positive constant
curvature before the solution vanishes. Consequently, a compact
three-manifold with positive Ricci curvature or a compact
four-manifold with positive curvature operator is diffeomorphic to
the round sphere or a quotient of it by a finite group of fixed
point free isometrics in the standard metric. In these classical
cases, the singularities are formed everywhere simultaneously and
with the same rates.

Note that even though the Ricci flow may develop singularities
everywhere at the same time, the singularities can still be formed
with
 different rates.  The general case is that the Ricci flow may
 develop singularities in some
parts while keeps smooth in other parts for general initial
metrics. This suggests that we have to consider the structures of
all the singularities (fast or slow forming). For the general
case, naturally one would like to cut off the singularities and to
continue the Ricci flow. If the Ricci flow still develops
singularity after a while, one can do the surgeries and run the
Ricci flow again. By repeating this procedure, one will get a kind
of ``weak" solution to the Ricci flow. Furthermore, if the ``weak"
solution has only a finite number of surgeries at any finite time
interval and one can remember what had been cut during the
surgeries, as well as the ``weak" solution has a well-understood
long-time behavior, then one will also get the topology structure
of the initial manifold. This surgerically modified Ricci flow was
initially developed by Hamilton \cite{Ha7} for compact
four-manifolds. More recently, the idea of the Ricci flow with
surgery was further developed by Perelman \cite{P2} for compact
three-manifolds (see \cite{CaZ} for complete detail).

Let us give a brief description for the arguments of Hamilton in
\cite{Ha7}. Recall that a Riemannian four-manifold is said to have
\textbf{positive isotropic curvature} if for every orthonormal
four-frame the curvature tensor satisfies
$$R_{1313}+R_{1414}+R_{2323}+R_{2424}>2R_{1234}.$$
An \textbf{incompressible space form} $N^3$ in a four-manifold $M^4$
is a three-dimensional submanifold diffeomorphic to
$\mathbb{S}^3/\Gamma$ (the quotient of the three-sphere by a group
of isometries without fixed point) such that the fundamental group
$\pi_1(N^3)$ injects into $\pi_1(M^4)$. The space form is said to be
\textbf{essential} unless $\Gamma=\{1\}$, or $\Gamma=\mathbb{Z}_2$
and the normal bundle is non-orientable. In \cite{Ha7}, Hamilton
considered a compact four-manifold $M^4$ with no essential
incompressible space-form and with a metric of positive isotropic
curvature. He used this metric as initial data, and evolved it by
the Ricci flow. From the evolution equations of curvatures, one can
easily see that the curvature will become unbounded in finite time.
Under the positive isotropic curvature assumption, he proved that as
the time tends to the first singular time, either the solution has
positive curvature operator everywhere, or it contains a neck, a
region where the metric is very close to the product metric on
$\mathbb{S}^3\times \mathbb{I}$, where $\mathbb{I}$ is an interval
and $\mathbb{S}^3$ is a round three-sphere, or a quotient of this by
a finite group acting freely. When the solution has positive
curvature operator everywhere, it is diffeomorphic to $\mathbb{S}^4$
or $\mathbb{RP}^4$ by \cite{Ha2}, so the topology of the manifold is
understood and one can throw it away. When there is a neck in the
solution, he used the no essential incompressible space form
assumption to conclude that the neck must be $\mathbb{S}^3\times
\mathbb{I}$ or $\mathbb{S}^3\times\mathbb{I} /\mathbb{Z}_2$ where
$\mathbb{Z}_2$  acts antipodally on $\mathbb{S}^3$ and by reflection
on $\mathbb{I}$. For the first case, one can replace
$\mathbb{S}^3\times \mathbb{I}$ with two caps (i.e. two copies of
the differential four-ball $\mathbb{B}^4$) by cutting the neck and
rounding off the neck. While for the second case, one can do the
quotient surgery to eliminate an $\mathbb{RP}^4$ summand. In
\cite{Ha7}, Hamilton performed these cutting and gluing surgery
arguments so carefully that the positive isotropic curvature
assumption and the improved pinching estimates are preserved under
the surgeries. It is not hard to show that, after surgery, the new
manifold still has no essential incompressible space form. Then by
using this new manifold as initial data, one can run the Ricci flow
and do the surgeries again. These arguments were given in Section
A-D of \cite{Ha7}. In the last section (Section E) of \cite{Ha7},
Hamilton showed that after a finite number of surgeries in finite
time, and discarding a finite number of pieces which are
diffeomorphic to $\mathbb{S}^4,$ $\mathbb{RP}^4$, the solution
becomes extinct. This concludes that the four-manifold is
diffeomorphic to $\mathbb{S}^4$, $\mathbb{RP}^4$, $\mathbb{S}^3
\times \mathbb{S}^1$, the twisted product
$\mathbb{S}^3\widetilde{\times} \mathbb{S}^1$ ( i.e.,
$\mathbb{S}^{3}\tilde{\times}\mathbb{S}^{1}=\mathbb{S}^{3}\times
\mathbb{S}^{1}/\mathbb{Z}_{2},$ where $\mathbb{Z}_{2}$ flips
$\mathbb{S}^{3}$ antipodally and rotates $\mathbb{S}^{1}$ by
$180^{0}$), or a connected sum of them.

The celebrated paper \cite{P2} tells us how to recognize the
formation of singularities and how to perform the surgeries. One
can see from Section A to D of \cite{Ha7} that every statement is
accurate and every proof is complete, precise and detailed.
Unfortunately, the last section (Section E) contains some
unjustified statements, which have been known for the experts in
this field for several years. For example one can see the comment
of Perelman in \cite{P2} (Page 1, the second paragraph) and one
can also check that the proof of Theorem E 3.3 of \cite{Ha7} is
incomplete (in Proposition 3.4 of the present paper, we will prove
a stronger version of Theorem E 3.3 of \cite{Ha7}). The key point
is how to prevent the surgery times from accumulated (furthermore,
it requires to perform only a finite number of surgeries in each
finite time interval). By inspecting the last section of
\cite{Ha7}, it seems that surgeries were taken on the parts where
the singularities are formed from the global maximum points of
curvature. Intuitively, the other parts, where the curvatures go
to infinity also but not be comparable to the global maximums,
will still develop singularities shortly after surgery if one only
performs the surgeries for the global maximum points of curvature.
To prevent the surgery times from accumulated, one needs to cut
off those singularities (not just the curvature maximum points)
also. This says that one needs to perform surgeries for all
singularities. Another problem is that, when one performs the
surgeries with a given accuracy at each surgery time, it is
possible that the errors may add up to a certain amount which
causes the surgery times to accumulate. To prevent this from
happening, as time goes on, successive surgeries must be performed
with increasing accuracy.

Recently, Perelman \cite{P1}, \cite{P2} presented the striking
ideas how to understand the structures of all singularities of the
three-dimensional Ricci flow, how to find ``fine" necks, how to
glue ``fine" caps, and how to use rescaling to prove that the
times of surgery are discrete. When using rescaling arguments for
surgically modified solutions of the Ricci flow, one encounters
the difficulty of how to apply Hamilton's compactness theorem,
which works only for smooth solutions.  To overcome the
difficulty, Perelman argued in \cite{P2} by choosing the cutoff
radius in necklike regions small enough to push the surgical
regions far away in space. But it still does not suffice to take a
smooth limit since Shi's interior derivative estimate is not
available, and so one cannot be certain that Hamilton's
compactness result holds when only having the bound on curvatures.
This is discussed in \cite{CaZ} and this paper.

In this paper, inspired by Perelman's works, we will study the
Ricci flow on compact four-manifolds with positive isotropic
curvature and with no essential incompressible space-form. We will
give a complete proof for the main theorem of Hamilton in
\cite{Ha7}. One of our major contribution in this paper is to
establish several time-extension results for the surgical
solutions in the proof of the discreteness of surgery times so
that the surgical solutions are smooth on some uniform (small)
time intervals (on compact subsets) and Hamilton's compactness
theorem is still applicable. In Perelman's works \cite{P1, P2},
the universal noncollapsing property of singularity models is a
crucial fact to prove the surviving of noncollapsing property
under surgery. Another feature of this paper is our proof on this
crucial fact. In dimension three, one obtains this by using
Perelman's classification of three-dimensional shrinking Ricci
solitons with nonnegative curvature (see \cite{CaZ} for the
details). But in the present four dimension case, we are not able
to obtain a complete classification for shrinking solitons. In the
previous version, we presented an argument to obtain the universal
noncollapsing for shrinking solitons. But, as pointed out to us by
Joerg Enders, that argument contains a gap. Fortunately in the
present version, we find a new argument, without appealing a
classification of shrinking Ricci solitons, to get the universal
noncollapsing for all possible singularity models.

During the proof we have actually provided, up to slight
modifications, all necessary details for the part from Section 1 to
Section 5 of Perelman's second paper \cite{P2} on Ricci flow to
approach the Poincar\'e conjecture. The complete details of the
arguments in three-dimension can be found in the recent paper of
H.-D. Cao and the second author in \cite{CaZ}. Furthermore, a
complete proof to the Poincar\'e conjecture and Thurston's
geometrization conjecture has been given in \cite{CaZ}.

The main result of this paper is the following \vskip 0.3cm

$\underline{\mbox{\textbf{Theorem 1.1}}}$ \emph{Let $M^4$ be a
compact four-manifold with no essential incompressible space-form
and with a metric $g_{ij}$ of positive isotropic curvature. Then
we have a finite collection of smooth solutions $g^{(k)}_{ij}(t)$,
$k=0,1,\cdots,m$, to the Ricci flow, defined on
$M^4_k\times[t_k,t_{k+1})$, $(0=t_0<\cdots<t_{m+1})$ with
$M^4_0=M^4$ and $g^{(0)}_{ij}(t_0)=g_{ij}$, which go singular as
$t\rightarrow t_{k+1}$, such that the following properties hold:}

(i) \emph{for each $k=0,1,\cdots,m-1$, the compact (possible
disconnected) four-manifold $M^4_k$ contains an open set
$\Omega_k$ such that the solution $g^{(k)}_{ij}(t)$ can be
smoothly extended to $t=t_{k+1}$ over $\Omega_k$; }

(ii) \emph{for each $k=0,1,\cdots,m-1$,
$(\Omega_k,g^{(k)}_{ij}(t_{k+1}))$ and
$(M^4_{k+1},g^{(k+1)}_{ij}(t_{k+1}))$ contain compact (possible
disconnected) four-dimensional submanifolds with smooth boundary,
which are isometric and then can be denoted by $N^4_k$;}

(iii) \emph{for each $k=0,1,\cdots,m-1$, $M_{k}^{4}\setminus
N_{k}^{4}$ consists of a finite number of disjoint pieces
diffeomorphic to $\mathbb{S}^{3}\times \mathbb{I}$,
$\mathbb{B}^{4}$ or $\mathbb{RP}^{4}\setminus \mathbb{B}^{4}$,
while $M^4_{k+1}\setminus N^4_k$ consists of a finite number of
disjoint pieces diffeomophic to $\mathbb{B}^4$;}

(iv) \emph{for $k=m$, $M^4_m$ is diffeomorphic to the disjoint union
of a finite number of $\mathbb{S}^4$, or $\mathbb{RP}^4$, or
$\mathbb{S}^3\times \mathbb{S}^1$, or $\mathbb{S}^3
\widetilde{\times}\mathbb{S}^1$, or $\mathbb{RP}^4 \#
\mathbb{RP}^4$.} \vskip 0.3cm

 As a direct consequence we have the following
classification result of Hamilton \cite{Ha7}.\vskip 0.3cm

$\underline{\mbox{\textbf{Corollary 1.2}}}$( Hamilton \cite{Ha7})
\emph{A compact four-manifold with no essential incompressible
space-form and with a metric of positive isotropic curvature is
diffeomorphic to $\mathbb{S}^4$, or $\mathbb{RP}^4$, or
$\mathbb{S}^3\times \mathbb{S}^1$, or
$\mathbb{S}^3\widetilde{\times}\mathbb{S}^1$,or a connected sum of
them.} \vskip 0.3cm

This paper contains five sections and an appendix. In Section 2 we
recall the pinching estimates of Hamilton obtained in \cite{Ha7} and
present two useful geometric properties for complete noncompact
Riemannian manifolds with positive sectional curvature. The usual
way to understand the singularities of the Ricci flow is to take a
rescaling limit and to find the structure of the limiting models. In
Section 3 we study the limiting models, so called ancient
$\kappa$-solutions. We will establish the uniform
$\kappa$-noncollapsing, compactness and canonical neighborhood
structures for ancient $\kappa$-solutions. These generalize the
analogs results of Perelman \cite{P1} from three-dimension to
four-dimension. In Section 4 we will extend the canonical
neighborhood characterization to any solution of the Ricci flow with
positive isotropic curvature, and  describe  the structure of the
solution at the singular time. In Section 5, we will define the
Ricci flow with surgery as Perelman in \cite{P2}. By a long
inductive argument, we will obtain a long-time existence result for
the surgerically modified Ricci flow so that the solution becomes
extinct in finite time and takes only a finite number of surgeries.
This will give the proof of the main theorem. In the appendix we
will prove the curvature estimates for the standard solutions and
give the canonical neighborhood description of the standard solution
in dimension four, which are used in Section 5 for the surgery
arguments. \vskip 0.8cm
 \centerline{{\textbf{Table of Contents}}}
\  \\
1. Introduction\\
2. Preliminaries \\
3. Ancient solutions

   3.1 Splitting lemmas

   3.2 Elliptic type estimate, canonical neighborhood decomposition

   \ \ \ \ \  for noncompact $\kappa$-solutions

   3.3 Universal noncollapsing of ancient $\kappa$-solutions

   3.4 Canonical neighborhood structures\\
4. The structure of the solutions at the singular time\\
5. Ricci flow with surgery for four manifolds\\
Appendix. Standard solutions

 \vskip 0.8cm

 We are grateful to Professor H.-D. Cao
for many helpful discussions and Professor S.-T. Yau for his
interest and encouragement. We also thank Joerg Enders for telling
us an error in the previous version. The first author is partially
supported by FANEDD 200216 and NSFC 10401042 and the second author
is partially supported by NSFC 10428102 and the IMS of The Chinese
University of Hong Kong.

\vskip 1cm \centerline{\large{\textbf{2. Preliminaries}}} \vskip
0.5cm
 Consider a four-dimensional compact Riemannian
manifold $M^4$. The curvature tensor of $M^4$ may be regarded as a
symmetric bilinear form $M_{\alpha\beta}$ on the space of real
forms $\Lambda^2$. It is well known that one can decompose
$\Lambda^2$ into $\Lambda^2_+\oplus \Lambda^2_-$ as eigen-spaces
of the Hodge star operator with eigenvalues $\pm 1$. This gives a
block decomposition of the curvature operate ($M_{\alpha\beta}$)
as

$$(M_{\alpha\beta})=
\left(
\begin{array}{cc}
A & B\\
^tB & C\\
 \end{array}\right).$$
It was shown in Lemma A2.1 of \cite{Ha7} that a four-manifold has
positive isotropic curvature if and only if $$a_1+a_2>0 \mbox{ and
} c_1+c_2>0$$ where $a_3\geq a_2\geq a_1$, $c_3\geq c_2 \geq c_1$
are eigenvalues of the matrices $A$ and $C$ respectively.

Let $\{X_1, X_2, X_3, X_4\}$ be a positive oriented orthonormal
basis of one-forms. Then $\varphi_1=X_1\wedge X_2+X_3\wedge X_4$,
$\varphi_2=X_1\wedge X_3+X_4\wedge X_2$, $\varphi_3=X_1\wedge
X_4+X_2\wedge X_3$ is a basis of $\Lambda^2_+$ and
$\psi_1=X_1\wedge X_2-X_3\wedge X_4$, $\psi_2=X_1\wedge
X_3-X_4\wedge X_2$, $\psi_3=X_1\wedge X_4-X_2\wedge X_3$ is a
basis of $\Lambda^2_-$. It is easy to check $trA=trC=\frac{1}{2}R$
by using this orthonormal basis and the Bianchi identity.

Since $B$ may not be symmetric, its eigenvalues need to be
explained as follows. For an appropriate choice of orthonormal
bases $y^+_1, y^+_2,y^+_3$ of $\Lambda^2_+$ and $y^-_1, y^-_2,
y^-_3$ of $\Lambda^2_-$ the matrix
$$B= \left(
\begin{array}{ccc}
b_1 & 0 & 0\\
0 & b_2 & 0\\
0 & 0 & b_3 \\
 \end{array}\right).$$
with $0\leq b_1 \leq b_2 \leq b_3$. They are actually the
eigenvalues of the symmetric matrices $\sqrt{B^tB}$ or $\sqrt{^t
BB}$.

In \cite{Ha7} Hamilton proved that the Ricci flow on a compact
four-manifold preserves positive isotropic curvature and obtained
the following improving pinching estimate.\vskip 0.3cm

$\underline{\mbox{\textbf{Lemma 2.1}}}$ (Theorem B1.1 and Theorem
B2.3 of \cite{Ha7})

\emph{Given an initial metric on a compact four-manifold with
positive isotropic curvature, there exist positive constants
$\rho$, $\Lambda,P<+\infty$ depending only on the initial metric,
such that the solution to the Ricci flow satisfies
$$a_1+\rho>0 \mbox{ and }c_1+\rho>0, \eqno (2.1)$$
$$\max\{a_3,b_3,c_3\}\leq \Lambda (a_1+\rho) \mbox{ and
}\max\{a_3,b_3,c_3\}\leq \Lambda (c_1+\rho), \eqno (2.2)$$ and
$$\frac{b_3}{\sqrt{(a_1+\rho)(c_1+\rho)}}\leq
1+\frac{\Lambda e^{Pt}}{\max\{\log\sqrt{(a_1+\rho)(c_1+\rho)},2\}}
\eqno (2.3)$$ at all points and all times.}\vskip 0.2cm

This lemma tells us that as we consider the Ricci flow for a
compact four-manifold with positive isotropic curvature, any
rescaling limit along a sequence of points where the curvatures
become unbounded must still have positive isotropic curvature and
satisfies the following restricted isotropic curvature pinching
condition $$a_3\leq \Lambda a_1,\ \  c_3\leq \Lambda c_1,\ \
b_3^2\leq a_1c_1. \eqno (2.4)$$

In the rest of this section, we will give two useful geometric
properties for Riemannian manifolds with nonnegative sectional
curvature.

Let $(M^n,g_{ij})$ be an $n$-dimensional complete Riemannian
manifold and let $\varepsilon$ be a positive constant. We call an
open subset $N\subset M^n$ to be an \textbf{$\varepsilon$-neck of
radius $r$} if $(N, r^{-2}g_{ij})$ is $\varepsilon$-close, in
$C^{[\varepsilon^{-1}]}$ topology, to a standard neck
$\mathbb{S}^{n-1}\times \mathbb{I}$ with $\mathbb{I}$ of the length
$2\varepsilon^{-1}$ and $\mathbb{S}^{n-1}$ of the scalar curvature
1.\vskip 0.3cm

$\underline{\mbox{\textbf{Proposition 2.2}}}$ \emph{There exists a
constant $\varepsilon_0=\varepsilon_0(n)>0$ such that every complete
noncompact Riemannian manifold $(M^{n}, g_{ij})$ of nonnegative
sectional curvature has a positive constant $r_0$ such that any
$\varepsilon$-neck of radius $r$ on $(M^{n}, g_{ij})$ with
$\varepsilon \leq \varepsilon_0$ must have $r\geq r_0$. }\vskip
0.2cm

$\underline{\mbox{\textbf{Proof}}}$.
 We argue by contradiction.
Suppose there exists a sequence of positive constants
$\varepsilon^\alpha \rightarrow 0$ and a sequence of $n$-dimensional
complete noncompact Riemannian manifolds $(M^\alpha,
g^{\alpha}_{ij})$ with nonnegative sectional curvature such that for
each fixed $\alpha$, there exists a sequence of
$\varepsilon^\alpha$-necks $N_{k}$ of radius at most $1/k$ on
$M^\alpha$ with centers $P_{k}$ divergent to infinity.

Fix a point $P$ on the manifold $M^\alpha$ and connect each $P_k$
to $P$ by a minimizing geodesic $\gamma_k$. By passing to
subsequence we may assume the angle $\theta_{kl}$ between geodesic
$\gamma_k$ and $\gamma_l$ at $P$ is very small and tends to zero
as $k, l\rightarrow+\infty$, and the length of $\gamma_{k+1}$ is
much bigger than the length of $\gamma_k$. Let us connect $P_k$ to
$P_l$ by a minimizing geodesic $\eta_{kl}$. For each fixed $l>k$,
let $\tilde{P}_k$ be a point on the geodesic $\gamma_l$ such that
the geodesic segment from $P$ to $\tilde{P}_k$ has the same length
as $\gamma_k$ and consider the triangle $\Delta PP_k\tilde{P}_k$
in $M^\alpha$ with vertices $P$, $P_k$ and $\tilde{P}_k$. By
comparing with the corresponding triangle in the Euclidean plane
$\mathbb{R}^2$ whose sides have the same corresponding lengths,
Toponogov comparison theorem implies
$$d(P_k,\tilde{P}_k)\le 2\sin(\frac{1}{2}\theta_{kl})\cdot d(P_k,P).$$
Since $\theta_{kl}$ is very small, the distance from $P_k$ to the
geodesic $\gamma_l$ can be realized by a geodesic $\zeta_{kl}$
which connects $P_k$ to a point $P_k'$ on the interior of the
geodesic $\gamma_l$ and has length at most
$2\sin(\frac{1}{2}\theta_{kl})\cdot d(P_k,P).$ Clearly the angle
between $\zeta_{kl}$ and $\gamma_l$ at the intersection point
$P_k'$ is $\frac{\pi}{2}$. Consider $\alpha$ to be fixed and
sufficiently large. We claim that as $k$ large enough, each
minimizing geodesic $\gamma_l$ with $l>k$, connecting $P$ to
$P_l$, goes through the neck $N_k$.

Suppose not, then the angle between $\gamma_k$ and $\zeta_{kl}$ at
$P_k$ is close to either zero or $\pi$ since $P_k$ is in the
center of an $\varepsilon^\alpha$-neck and $\alpha$ is
sufficiently large. If the angle between $\gamma_k$ and
$\zeta_{kl}$ at $P_k$ is close to zero, we consider the triangle
$\Delta PP_kP_k'$ in $M^\alpha$ with vertices $P$, $P_k$, and
$P_k'$. By applying Toponogov comparison theorem to compare the
angles of this triangle with those of the corresponding triangle
in the Euclidean plane with the same corresponding lengths, we
find that it is impossible. Thus the angle between $\gamma_k$ and
$\zeta_{kl}$ at $P_k$ is close to $\pi$. We now consider the
triangle $\Delta P_kP_k'P_l$ in $M^\alpha$ with the three sides
$\zeta_{kl}$, $\eta_{kl}$ and the geodesic segment from $P_k'$ to
$P_l$ on $\gamma_l$. We have seen that the angle of $\Delta P_k
P_k'P_l$ at $P_k$ is close to zero and the angle at $P_k'$ is
$\frac{\pi}{2}$. By comparing with corresponding triangle
$\bar{\Delta}\bar{P_k}\bar{P_k'}\bar{P_l}$ in the Euclidean plane
$\mathbb{R}^2$ whose sides have the same corresponding lengths,
Toponogov comparison theorem implies
$$\angle \bar{P_l}\bar{P_k}\bar{P_k'}+\angle \bar{P_l}\bar{P_k'}\bar{P_k}\le
\angle P_lP_kP_k'+\angle P_lP_k'P_k<\frac{3}{4}\pi.$$ This is
impossible since the length between $\bar{P_k}$ and $\bar{P_k'}$
is much smaller than the length from $\bar{P_l}$ to either
$\bar{P_k}$ or $\bar{P_k'}$. So we have proved each $\gamma_l$
with $l>k$ passes through the neck $N_k$.

Hence by taking a limit, we get a geodesic ray $\gamma$ emanating
from $P$ which passes through all the necks $N_k$, $k = 1, 2,
\cdots,$  except a finite number of them. Throwing these finite
number of necks, we may assume $\gamma$ passes through all necks
$N_k$, $k=1,2,\cdots.$ Denote the center sphere of $N_k$ by $S_k$,
and their intersection points with $\gamma$ by  $p_{k}\in S_k\cap
\gamma$, for $k=1,2,\cdots.$

Take a sequence points  $\gamma(m)$ with $m=1,2,\cdots.$ For each
fixed neck $N_k$, arbitrarily choose a point $q_{k}\in N_k$ near the
center sphere $S_k$, draw a geodesic segment $\gamma^{km}$ from
$q_{k}$ to $\gamma(m)$. Now we claim that for any fixed neck $N_l$
with $l>k$, $\gamma^{km}$ will pass through $N_l$ for all
sufficiently large $m$.

 We argue by contradiction. Let us place the
all necks $N_i$ horizontally so that the geodesic $\gamma$ passes
through each $N_i$ from the left to the right. We observe that the
geodesic segment $\gamma^{km}$ must pass through the right half of
$N_k$; otherwise $\gamma^{km}$ can not be minimal. Then as $m$
large enough, the distance from $p_{l}$ to the geodesic segment
$\gamma^{km}$ must be achieved by the distance from $p_l$ to some
interior point ${p_k}'$ of $\gamma^{km}$. Let us draw a minimal
geodesic $\eta$ from $p_{l}$ to the interior point ${p_k}'$ with
the angle at the intersection point ${p_k}'\in \eta\cap
\gamma^{km}$ to be $\frac{\pi}{2}.$ Suppose the claim is false.
Then the angle between $\eta$ and $\gamma$ at $p_{l}$ is close to
$0$ or $\pi$ since $\varepsilon^\alpha$ is small.

If the angle between $\eta$ and $\gamma$ at $p_{l}$ is close to
$0$, we consider the triangle ${\Delta} {p}_{l}{p_k}'{\gamma}(m)$
and construct a comparison triangle $ \bar{\Delta}
\bar{p}_{l}\bar{p_k}'\bar{\gamma}(m)$ in the plane with the same
corresponding length. Then by Toponogov comparison, we see the sum
of the inner angles of the comparison triangle $ \bar{\Delta}
\bar{p}_{l}\bar{p_k}'\bar{\gamma}(m)$ is less than $3\pi/4$, which
is impossible.

If the angle between $\eta$ and $\gamma$ at $p_{l}$ is close to
$\pi$, by drawing a minimal geodesic from $\xi$ from $q_k$ to
$p_{l}$, we see that $\xi$ must pass through the right half of
$N_k$ and the left half of $N_l$; otherwise $\xi$ can not be
minimal. Thus the three inner angles of the triangle $\Delta
p_{l}{p_k}'q_k$ are almost $0,\pi/2,0$ respectively. This is also
impossible by Toponogov comparison theorem.

Hence we have proved that the geodesic segment $\gamma^{km}$
passes through $N_l$ as $m$ large enough.

Consider the triangle $\Delta p_{k}q_k\gamma(m)$ with two long
sides $\overline{p_{k}\gamma(m)}(\subset\gamma)$ and
$\overline{q_{k}\gamma(m)}(= \gamma^{km})$. For any $s>0$, choose
two points ${\tilde{p}_{k}}$ on $\overline{p_{k}\gamma(m)}$  and
${\tilde{q}_{k}}$ on $\overline{q_{k}\gamma(m)}$ with
$d(p_{k},{\tilde{p}_{k}})=d(q_{k},{\tilde{q}_{k}})=s$. By
Toponogov comparison theorem, we have
 \begin{eqnarray*} &
&(\frac{d({\tilde{p}_{k}},{\tilde{q}_{k}})}{d(p_{k},q_{k})})^{2}\\[4mm]
&=&
\frac{d({\tilde{p}_{k}},\gamma(m))^{2}+d({\tilde{q}_{k}},\gamma(m))^{2}-2d({\tilde{p}_{k}},\gamma(m))d({\tilde{q}_{k}},
\gamma(m))\cos\bar{\measuredangle}
({\tilde{p}_{k}}\gamma(m){\tilde{q}_{k}})}
{d({p_{k}},\gamma(m))^{2}+d({q_{k}},\gamma(m))^{2}-2d({p_{k}},\gamma(m))d({q_{k}},
\gamma(m))\cos\bar{\measuredangle}
({p_{k}}\gamma(m){q_{k}})}\\[4mm]
&\geq&
\frac{d({\tilde{p}_{k}},\gamma(m))^{2}+d({\tilde{q}_{k}},\gamma(m))^{2}-2d({\tilde{p}_{k}},\gamma(m))d({\tilde{q}_{k}},
\gamma(m))\cos\bar{\measuredangle}
({\tilde{p}_{k}}\gamma(m){\tilde{q}_{k}})}
{d({p_{k}},\gamma(m))^{2}+d({q_{k}},\gamma(m))^{2}-2d({p_{k}},\gamma(m))d({q_{k}},
\gamma(m))\cos\bar{\measuredangle}
({\tilde{p}_{k}}\gamma(m){\tilde{q}_{k}})}\\[4mm]
&=& \frac{(d({\tilde{p}_{k}},\gamma(m)) -
d({\tilde{q}_{k}},\gamma(m)))^{2}+2d({\tilde{p}_{k}},\gamma(m))d({\tilde{q}_{k}},
\gamma(m))(1 - \cos\bar{\measuredangle}
({\tilde{p}_{k}}\gamma(m){\tilde{q}_{k}}))}
{(d({\tilde{p}_{k}},\gamma(m)) -
d({\tilde{q}_{k}},\gamma(m)))^{2}+2d({{p}_{k}},\gamma(m))d({{q}_{k}},
\gamma(m))(1 - \cos\bar{\measuredangle}
({\tilde{p}_{k}}\gamma(m){\tilde{q}_{k}}))}\\[4mm]
&\geq& \frac{d({\tilde{p}_{k}},\gamma(m))d({\tilde{q}_{k}},
\gamma(m))} {d({p_{k}},\gamma(m))d({q_{k}}, \gamma(m))}\\[4mm]
& \rightarrow & 1
\end{eqnarray*}
 as $m\rightarrow\infty$, where $\bar{\measuredangle}
 ({{p}_{k}}\gamma(m){{q}_{k}})$ and $\bar{\measuredangle}
 ({\tilde{p}_{k}}\gamma(m){\tilde{q}_{k}})$ are the
 the corresponding angles in the corresponding comparison triangles.

Letting $m\rightarrow \infty$, we see that $\gamma^{km}$ has a
convergent subsequence whose limit $\gamma^{k}$ is a geodesic ray
passing through all $N_l$ with $l > k$. Denote by
$p_{j}=\gamma(t_j), j=1, 2, \cdots$. From the above computation,
we deduce that
$$d(p_{k},q_{k})\leq d(\gamma(t_k+s),\gamma^{k}(s)).
$$
for all $s>0$.

Let $\varphi(x)=\lim_{t\rightarrow+\infty}(t-d(x,\gamma(t)))$ be
the Busemann function constructed from the ray $\gamma$. Note that
the level set $\varphi^{-1}(\varphi(p_{j}))\cap N_j$ is close to
the center sphere $S_j$ for any $j = 1, 2, \cdots$. Now let $q_k$
be any fixed point in $\varphi^{-1}(\varphi(p_{k}))\cap N_k$. By
the definition of Busemann function $\varphi$ associated to the
ray $\gamma$, we see that
$\varphi(\gamma^{k}(s_1))-\varphi(\gamma^{k}(s_2))=s_1-s_2$ for
any $s_1$, $s_2\geq0$. Consequently, for each $l>k$, by choosing
$s=t_l-t_k$, we see $\gamma^{k}(t_l-t_k)\in
\varphi^{-1}(\varphi(p_{l}))\cap N_l.$ Since
$\gamma(t_k+t_l-t_k)=p_{l}$, it follows that
$$d(p_{k},q_{k})\leq d(p_{l},\gamma^{k}(s)).
$$
with $s =t_l-t_k >0$. This implies that the diameter of
$\varphi^{-1}(\varphi(p_{k}))\cap N_k$ is not greater the diameter
of $\varphi^{-1}(\varphi(p_{l}))\cap N_l$ for any $l>k$, which is
a contradiction as $l$ much larger than $k$.

Therefore we have proved the proposition.

$$\eqno \#$$
\vskip 0.3cm

In \cite{Ha6}, Hamilton discovered an interesting result, called
finite bump theorem, about the influence of a bump of strictly
positive curvature in a complete noncompact Riemannian manifold with
nonnegative sectional curvature. Namely, minimal geodesic paths that
go past the bump have to avoid it. The following result is in the
same spirit as Hamilton's finite bump theorem.\vskip 0.3cm

$\underline{\mbox{\textbf{Proposition 2.3}}}$ \emph{Suppose
$(M^{n},g)$ is a complete $n$-dimensional Riemanian manifold with
nonnegative sectional curvature. Let $P\in M^{n}$ be fixed, and
$P_j\in M^{n}$ a sequence of points and $R_j$ a sequence of positive
numbers with $d(P,P_j)\rightarrow +\infty$ and
$R_jd(P,P_j)^{2}\rightarrow +\infty$. If the sequence of marked
manifolds $(M^{n}, R_jg,P_{j})$ converges in $C^{\infty}_{loc}$
topology (in Cheeger sense) to a smooth manifold
$(\tilde{M}^{n},\tilde{g})$, then the limit
$(\tilde{M}^{n},\tilde{g})$ splits as the metric product of the form
$\mathbb{R}\times N$, where $N$ is a nonnegatively curved manifold
of dimension $n-1$.}\vskip 0.2cm

$\underline{\mbox{\textbf{Proof}}}$: Let us denote by
$|OQ|=d(O,Q)$ for the distance of two points $O,Q\in M^{n}$.
Without loss of generality, we may assume that for each $j$
$$1+2|PP_j|\leq|PP_{j+1}|. \eqno (2.5)$$ Draw a minimal geodesic
$\gamma_j$ from $P$ to $P_j$ and a minimal geodesic $\sigma_j$
from $P_j$ to $P_{j+1}$, both parameterized by the arclength. By
the compactness of unit sphere of the tangent space at $P$,
$\{\gamma'_j(0)\}$ has a convergent subsequence. We may further
assume
$$\theta_j=|\measuredangle(\gamma'_j(0),\gamma'_{j+1}(0))|<\frac{1}{j}.
\eqno (2.6)$$

Since $(M^{n}, R_jg,P_{j})$ converges in $C^{\infty}_{loc}$
topology (in Cheeger sense) to a smooth marked manifold
$(\tilde{M}^{n},\tilde{g},\tilde{P})$, by further choices of
subsequences, we may also assume $\gamma_j$ and $\sigma_j$
converge to two geodesic rays $\widetilde{\gamma}$ and
$\widetilde{\sigma}$ starting at $\widetilde{P}$. We claim that
 that $\tilde{\gamma}\cup\tilde{\sigma}$ forms a line in
$\tilde{M}^{n}$. Since the sectional curvature of $\tilde{M}^{n}$
is nonnegative, then by Toponogov splitting theorem \cite{CE} the
limit $\tilde{M}^{n}$ must  split as $\mathbb{R}\times N$
isometrically.

To prove the claim, we argue by contradiction. Suppose
$\widetilde{\gamma}\cup\widetilde{\sigma}$ is not a line, then for
each $j$, there exist two points $A_j\in \gamma_j$ and $B_j\in
\sigma_j$ such that as $j\rightarrow +\infty$,
$$\left\{
\begin{array}{llll}
R_jd(P_j,A_j)\rightarrow A>0,\\
R_jd(P_j,B_j)\rightarrow B>0,\\
R_jd(A_j,B_j)\rightarrow C>0,  \ \ \ \ \ \ \ \ \ \ \ \ \ \ \ \ \ \ \ \ \ \ \ \ \ (2.7)\\
\mbox{but  }A+B>C.\\
 \end{array}\right.$$

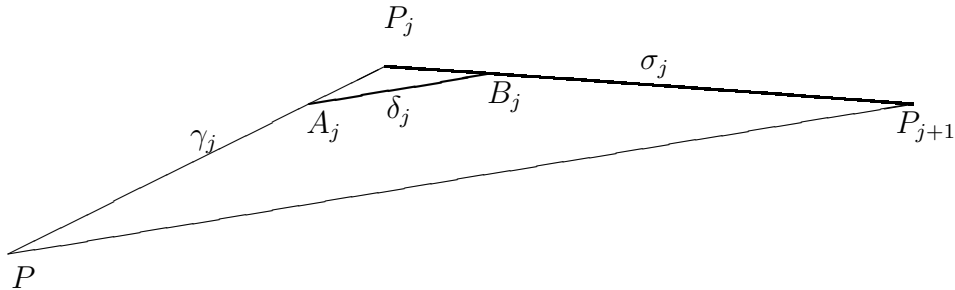
\begin{center}
\setlength{\unitlength}{2mm}
\begin{picture}(60,20)
\linethickness{1pt} \put(0,0){\line(6,1){60}}

\put(0,0){\line(2,1){25}} \thicklines

\qbezier(25,12.5)(40,11.5)(60,10)

\put(20,10){\line(6,1){12}}

\put(0,-2){\makebox(2,1)[c]{$P$}}

\put(25,15){\makebox(2,1)[c]{$P_j$}}

\put(60,8){\makebox(2,1)[c]{$P_{j+1}$}}

\put(20,8){\makebox(2,1)[c]{$A_j$}}

\put(25,9){\makebox(2,1)[c]{$\delta_j$}}

\put(32,10){\makebox(2,1)[c]{$B_j$}}

\put(42,12){\makebox(2,1)[c]{$\sigma_j$}}

\put(12,7){\makebox(2,1)[c]{$\gamma_j$}}

\end{picture}

\end{center}

Now draw a minimal geodesic $\delta_j$ from $A_j$ to $B_j$.
Consider comparison triangle
$\bar{\triangle}{\bar{P}_j\bar{P}\bar{P}_{j+1}}$ and
$\bar{\triangle}{\bar{P}_j\bar{A}_j\bar{B}_j}$ in $\mathbb{R}^2$
with
$$|\bar{P}_j\bar{P}|=|P_jP|, |\bar{P}_j\bar{P}_{j+1}|=|P_jP_{j+1}|,
|\bar{P}\bar{P}_{j+1}|=|PP_{j+1}|,$$ $$\mbox{ and
}|\bar{P}_j\bar{A}_j|=|P_jA_j|, |\bar{P}_j\bar{B}_j|=|P_jB_j|,
|\bar{A}_j\bar{B}_j|=|A_jB_j|.$$ By Toponogov comparison theorem
\cite{CE}, we have
$$\measuredangle\bar{A}_j\bar{P}_j\bar{B}_j\geq
\measuredangle \bar{P} \bar{{P}_j}\bar{P}_{j+1}. \eqno (2.8)$$ On
the other hand, by (2.6) and using the Toponogov comparison
theorem again, we have
$$\measuredangle\bar{P}_j\bar{P}\bar{P}_{j+1}\leq \measuredangle P_jPP_{j+1}<\frac{1}{j}, \eqno
(2.9)$$ and since $|\bar{P}_j\bar{P}_{j+1}|>|\bar{P}\bar{P}_j|$ by
(2.5), we further have
$$\measuredangle\bar{P}_j\bar{P}_{j+1}\bar{P}\leq \measuredangle\bar{P}_j\bar{P}\bar{P}_{j+1}<\frac{1}{j}. \eqno
(2.10)$$ Thus the above inequalities (2.8)-(2.10) imply that
$$\measuredangle\bar{A}_j\bar{P}_j\bar{B}_j>\pi-\frac{2}{j}.$$
Hence $$|\bar{A}_j\bar{B}_j|^2\geq
|\bar{A}_j\bar{P}_j|^2+|\bar{P}_j\bar{B}_j|^2-2|\bar{A}_j\bar{P}_j|\cdot|\bar{P}_j\bar{B}_j|\cos(\pi-\frac{2}{j}).
\eqno (2.11)$$

Multiplying the above inequality by $R_j$ and letting
$j\rightarrow+\infty$, we get $$C\geq A+B$$ which contradicts with
(2.7). Therefore we have proved the proposition.$$\eqno \#$$\vskip
0.3cm

$\underline{\mbox{\textbf{Corollary 2.4}}}$ \emph{Suppose $(X,d)$ is
a complete $n$-dimensional Alexandrov space with nonnegative
curvature. Let $P\in X$ be fixed, and $P_j\in X$ a sequence of
points and $R_j$ a sequence of positive numbers with
$d(P,P_j)\rightarrow +\infty$ and $R_jd^{2}(P,P_j)\rightarrow
+\infty$. Then the marked spaces $(X, R_j^{\frac{1}{2}}d,P_j)$ have
a (Gromov-Hausdorff) convergent subsequence such that the limit
splits as the metric product of the form $\mathbb{R}\times N$, where
$N$ is a nonnegatively curved Alexandrov space.}\vskip 0.2cm

$\underline{\mbox{\textbf{Proof}}}$: By the compactness theorem of
Alexandrov spaces (see \cite{BGP}), there is a subsequence of
$(X,R_j^{\frac{1}{2}}d,P_j)$, which converges (in the sense of
Gromov-Hausdorff) to a nonnegatively curved Alexandrov space
$(\widetilde{X},\widetilde{d},\widetilde{P})$ of dimension $\leq
n$. By Toponogov splitting theorem \cite{Mi} for Alexandrov
spaces, we only need to show that the limit $\widetilde{X}$
contains a line. Note that the same inequality (2.6) now follows
from the compactness of the space of directions at a fixed point
\cite{BGP}. Since the Toponogov triangle comparison theorem still
holds on Alexsandrov Spaces (in fact, the notion of the curvature
of general metric spaces is defined by Toponogov triangle
comparison), the same argument of the Proposition 2.3 proves the
corollary. $$\eqno \#$$

 \vskip 1cm \centerline{\large{\textbf{3. Ancient
Solutions}}} \vskip 0.5cm
 A solution to the Ricci flow on a compact
four-manifold with positive isotropic curvature develops
singularities in finite time. The usual way to understand the
formations of the singularities is to rescale the solution along
the singularities and to try to take a limit for the rescaled
sequences. According to Lemma 2.1, a rescaled limit will be a
complete non-flat solution to the Ricci flow
$$\frac{\partial}{\partial t}g_{ij}=-2R_{ij},$$ on an ancient time
interval $-\infty<t\leq 0$, called an \textbf{ancient solution},
which has positive isotropic curvature and satisfies the
restricted isotropic curvature pinching condition (2.4). We remark
that as we consider the general singularities (not be necessarily
those points coming from the maximum of the curvature ), we don't
know whether at a priori, the rescaled limit exists, and even
assuming the existence, whether the limit has bounded curvature
for each $t$. Nevertheless, in this section we will take the
attention to those rescaled limits with bounded curvature.

According to Perelman \cite{P1}, a solution to the Ricci flow is
$\textbf{$\kappa$-noncollapsed }$ $\textbf{for scale $r_0$}>0$ if we
have the following statement: whenever
 we have
$$|Rm|(x,t)\leq r_0^{-2},$$
$\mbox{for all } t\in [t_0-r_0^2,t_0], x\in B_t(x_0,r_0),$  for some
$(x_0,t_0)$,
  then there holds $$\ Vol_{t_0}(B_{t_0}(x_0,r_0))\geq \kappa
r_0^4.$$
 Here we denote by $B_{t}(x_0,r_0)$  and
 $Vol_{t_0}$ the geodesic ball centered
at $x_0$ of radius $r_0$ with respect to the metric $g_{ij}(t)$ and
the volume with respect to the metric $g_{ij}(t_0)$ respectively. It
was shown by Perelman \cite{P1} that any rescaled limit obtained by
blowing up a smooth solution to the Ricci flow on a compact manifold
in finite time is $\kappa$-noncollapsed on all scales for some
$\kappa>0$.

We say a solution to the Ricci flow on a four-manifold is an
\textbf{ancient $\kappa$} $\textbf{-solution with restricted
isotropic curvature pinching}$ (for some $\kappa>0$) if it is a
smooth solution to the Ricci flow on the ancient time interval
$t\in(-\infty,0]$ which is complete,
 has positive isotropic curvature and bounded curvature,
and satisfies the restricted isotropic curvature pinching condition
(2.4), as well as is $\kappa$-noncollapsed on all scales.
 \subsection*{3.1 Splitting lemmas}
 To understand the structures of the
solutions to the Ricci flow on a compact four-manifold with positive
isotropic curvature, we are naturally led to investigate the ancient
solutions which have positive isotropic curvature, satisfy the
restricted isotropic curvature pinching condition (2.4) and are
$\kappa$-noncollapsed for all scales. Note that the restricted
isotropic curvature pinching condition (2.4) implies the curvature
operator is nonnegative. In this subsection we will derive two
useful splitting results \textbf{\emph{without}} assuming bounded
curvature condition.\vskip 0.3cm

$\underline{\mbox{\textbf{Lemma 3.1}}}$ \emph{Let $(M^4, g_{ij})$ be
a complete noncompact Riemannian manifold which satisfies the
restricted isotropic curvature pinching condition (2.4) and has
positive curvature operator. And let $P$ be a fixed point in $M^4$,
$\{P_l\}_{1\leq l<+\infty}$ a sequence of points in $M^4$ and
$\{R_l\}_{1 \leq l<+\infty}$ a sequence of positive numbers with
$d(P,P_l)\rightarrow+\infty$ and $R_ld^2(P,P_l)\rightarrow+\infty$
as $l\rightarrow+\infty$ , where $d(P,P_l)$ is the distance between
$P$ and $P_l$. Suppose $(M^4,R_lg_{ij},P_l)$ converges in
$C^{\infty}_{loc}$ topology to a smooth nonflat limit $Y$. Then $Y$
must be isometric to $\mathbb{R}\times \mathbb{S}^3$ with the
standard metric (up to a constant factor).}\vskip 0.2cm

$\underline{\mbox{\textbf{Proof}}}$. By Proposition 2.3, we see
that $Y=\mathbb{R}\times X$ for some smooth three-dimensional
manifold $X$. Thus the block decomposition of the curvature
operator has the form
$$(M_{\alpha\beta})=
\left(
\begin{array}{cc}
A & A\\
A & A\\
 \end{array}\right).$$ The assumption that $b^2_3\leq a_1c_1$ in
 (2.4) implies that the matrix $A(=B=C)$ is a multiple of the
 identity. Since this is true at every point, it follows from the
 contracted second Bianchi identity that $X$ has (positive)
 constant curvature, i.e. $X=\mathbb{S}^3/\Gamma$ for some group
 $\Gamma$ of isometries without fixed points. We remain to show
 $X=\mathbb{S}^3$ (i.e. $\Gamma=\{1\}$).

 Note that the original manifold $M^4$ is diffeomorphic to
 $\mathbb{R}^4$ by the positive curvature operator assumption. To
 show $X=\mathbb{S}^3$ (i.e., $\Gamma=\{1\}$) we only need to prove that
 $X=\mathbb{S}^3/\Gamma$ is incompressible in $M^4$. In the following we
 adapt Hamilton's argument in Theorem C 4.1 of \cite{Ha7} to noncompact
 manifolds.

 Suppose $X=\mathbb{S}^3/\Gamma$ is not incompressible in $M^4$, then for
 $j$ large, there exists a simply closed curve $\gamma\subset
 \mathbb{S}^3/\Gamma$ (the center space form passing through $P_j$) which
 is homotopically nontrivial in $\mathbb{S}^3/\Gamma$ but bounds a disk $D^2$
 in $M^4$. By Lemma C4.2 in \cite{Ha7}, we may assume the disk $D^2$ meets
 the center space form $\mathbb{S}^3/\Gamma$ only at $\gamma$, where it is
 transversal. Now construct a new manifold $\widehat{M}^4$ in the
 following way. As in \cite{Ha7} we can deform the metric in the neck
 around $P_j$ a little so it is standard in a smaller neck but
 still has positive isotropic curvature everywhere. Since $M^{4}$ is
 simply connected, the connected
 and closed submanifold $\mathbb{S}^3/\Gamma$ of codimension 1 separates
 $M^{4}$ (see for example \cite{Hir}). Cut $M^4$ open along
 the center space form $\mathbb{S}^3/\Gamma$ to get a (maybe disconnected)manifold with two
 boundary components $\mathbb{S}^3/\Gamma$, and double across the boundary
 to get $\widehat{M}^4$. The new manifold $\widehat{M}^4$ also has
 a metric of positive isotropic curvature since the boundary is flat
 extrinsically and we can double the metric. If $\widehat{M}^4$
 contains a compact connected component and the above disk $D^2$
 also lies in the compact connected component, then the same
 argument as in the proof of Theorem C4.1 of \cite{Ha7} derives a
 contradiction. Thus we may assume that the unique noncompact
 connected component of $\widehat{M}^4$, denoted by $\widehat{M}^4_1$,
  contains the disk $D^2$. Since $M^4$ has positive
 curvature operator, we know from \cite{CG} that there is a strictly
 convex exhausting function $\varphi$ on $M^4$. We can define a
 function $\widehat{\varphi}$ on $\widehat{M}^4$ so that
 $\widehat{\varphi}=\varphi$ on each copy of $M^4$. Then as $c>0$
 is sufficiently large, the level set $\widehat{\varphi}^{-1}(c)$
is contained in the unique noncompact connected component
$\widehat{M}^4_1$ and is strictly convex (in the sense that its
second fundamental is strictly positive). Take our disk $D^2$
bounding the curve $\gamma$ and perpendicular to the boundary in a
neighborhood of the boundary, and double it across the boundary to
get a sphere $\mathbb{S}^2$ which is $\mathbb{Z}_2$ invariant and
intersects the boundary component transversally in $\gamma$. The
homotopy class $[\gamma]$ is nontrivial in $\mathbb{S}^3/\Gamma$.
Clearly the above two-sphere $\mathbb{S}^2$ is contained in the
set $\{x\in \widehat{M}|\widehat{\varphi}(x)\leq c\}$ as $c>0$
large enough, since $\varphi$ is an exhausting on $M^4$. Now fix
such a large positive constant $c$. Among all spheres which are
$\mathbb{Z}_2$ invariant, contained in the manifold $\{x\in
\widehat{M}^4 |\widehat{\varphi}(x)\leq c\}$, and intersect the
$\mathbb{S}^3/\Gamma$ in the homotopy class $[\gamma]\neq 0$,
there will be one of least area since the boundary of the manifold
$\{x\in \widehat{M}^4 |\widehat{\varphi}(x)\leq c\}$ is strictly
convex. This sphere must even have least area among all nearby
spheres. For if a nearby sphere of less area divides in two parts
bounding $[\gamma]$ in $\mathbb{S}^3/\Gamma$, one side or the
other has less than half the area of the original sphere. We could
then double this half to get a sphere of less area which is
$\mathbb{Z}_2$ invariant, contradicting the assumption that ours
was of least area among this class. But the hypothesis of positive
isotropic curvature implies there are no stable minimal
two-spheres as was shown in \cite{MiMo}. Hence we get a
contradiction unless $X=\mathbb{S}^3/\Gamma$ is incompressible in
$M^4$.

Therefore we have proved the lemma.$$\eqno \#$$

$\underline{\mbox{\textbf{Lemma 3.2}}}$ \emph{Let $(M^4, g_{ij}(t))$
be an ancient solution which has positive isotropic curvature and
satisfies the restricted isotropic curvature pinching condition
(2.4). If its curvature operator has a nontrivial null eigenvector
somewhere at some time, then the solution is, up to a scaling, the
evolving round cylinder $\mathbb{R}\times \mathbb{S}^3$ or a metric
quotient of the round cylinder $\mathbb{R}\times
\mathbb{S}^3$}.\vskip 0.3cm

$\underline{\mbox{\textbf{Proof.}}}$ \ \  Recall that the solution
$g_{ij}(t)$ has nonnegative curvature operator everywhere and
every time. Because the curvature operator of the ancient solution
$g_{ij}(t)$ has a nontivial null eigenvector somewhere at some
time, it follows from \cite{Ha2} (by using Hamilton's strong
maximum principle) that at any earlier time the solution has null
eigenvector everywhere and the Lie algebra of the holonomy group
is restricted a proper subalgebra of $so(4)$. Since the ancient
solution is nonflat and has positive isotropic curvature, we rule
out the subalgebras $\{1\}$, $u(2),so(2)\times so(2),so(2)\times
\{1\}$ as on $\mathbb{R}^4,\mathbb{CP}^2, \mathbb{S}^2\times
\mathbb{S}^2, \mathbb{S}^2\times \mathbb{R}^2$ or a metric
quotient of them. The only remaining possibility for the Lie
subalgebra of the holonomy is $so(3)$.

Now the only way we get holonomy $so(3)$ is when in some basis we
have $A=B=C$ in the curvature operator matrix, so that
$$(M_{\alpha\beta})= \left(
\begin{array}{cc}
A & A\\
A & A\\
 \end{array}\right),$$ which corresponds to the fact that the
 metric $g_{ij}(t)$ is locally a product of $\mathbb{R}\times X$
 for some smooth three-dimensional manifold $X$ with curvature
 operator $A$. Then the inequality $b^2_3\leq a_1c_1$ in the
 restricted isotropic curvature pinching condition (2.4) implies
 that $A$ is a multiple of the identity. Moreover this is true at
 every point, it follows from the contracted second Bianchi
 identity that the factor $X$ has (positive) constant curvature.
 Consequently, $X$ is compact and then for each $t$, the metric
 $g_{ij}(t)$ is isometric to (up to a scaling) the evolving metric of the round
 cylinder
 $\mathbb{R}\times \mathbb{S}^3$ or a metric quotient of it.
 $$\eqno \#$$

\subsection*{3.2 Elliptic type estimate, canonical neighborhood decomposition  for noncompact $\kappa$-solutions}
\vskip 0.5cm
 The following elliptic type Harnack property for
four-dimensional ancient $\kappa$-solutions with restricted
isotropic curvature pinching will be crucial for the analysis of the
structure of singularities of the Ricci flow on four-manifold with
positive isotropic curvature. The analogous result for
three-dimensional ancient $\kappa$-solutions was implicitely given
by Perelman in Section 11.7 of \cite{P1} and Section 1.5 of
\cite{P2}. \vskip 0.3cm

$\underline{\mbox{\textbf{Proposition 3.3}}}$ \emph{ For any
$\kappa>0$, there exist a positive constant $\eta$ and a positive
function $\omega:[0,+\infty)\rightarrow(0,+\infty)$ with the
following properties. Suppose we have a four-dimensional ancient
$\kappa$-solution $(M^4,g_{ij}(t)),-\infty<t\leq 0$, with
restricted isotropic curvature pinching. Then}

(i) \emph{for every $x,y\in M^4$ and $t\in (-\infty,0]$, there
holds
$$R(x,t)\leq R(y,t)\cdot \omega(R(y,t)d^2_t(x,y));$$}
 \ \ (ii) \emph{for all $x\in M^4$ and $t\in (-\infty,0]$, there hold
$$|\nabla R|(x,t)\leq \eta R^{\frac{3}{2}}(x,t) \mbox{ and }|\frac{\partial R}{\partial t}|(x,t)\leq \eta
R^2(x,t).$$ }
\vskip 0.2cm$\underline{\mbox{\textbf{Proof}}}$.
Obviously we may assume the ancient $\kappa$-solution is not a
metric quotient of the round neck $\mathbb{R}\times \mathbb{S}^3$.

(i) We only need to establish the estimate at $t=0$. Let $y$ be
fixed in $M^4$. By rescaling, we can assume $R(y,0)=1$.

Let us first consider the case that $\sup\{R(x,0)d^2_0(x,y)|x\in
M^4\}>1$. Define $z$ to be the closest point to $y$ (at time
$t=0$) satisfying $R(z,0)d^2_0(z,y)=1$. We want to bound
$R(x,0)/R(z,0)$ from above for $x\in
B_0(z,2R(z,0)^{-\frac{1}{2}})$.

Connect $y$ to $z$ by a shortest geodesic and choose a point
$\widetilde{z}$ lying on the geodesic satisfying
$d_0(\widetilde{z},z)=\frac{1}{4}R(z,0)^{-\frac{1}{2}}$. Denote by
$B$ the ball centered at $\widetilde{z}$ and with radius
$\frac{1}{4}R(z,0)^{-\frac{1}{2}}$ (with respect to the metric at
$t=0$). Clearly the ball $B$ lies in
$B_0(y,R(z,0)^{-\frac{1}{2}})$ and lies outside
$B_0(y,\frac{1}{2}R(z,0)^{-\frac{1}{2}})$. Thus as $x\in B$, we
have $$R(x,0)d^2_0(x,y)\leq 1 \mbox{ and }d_0(x,y)\geq
\frac{1}{2}R(z,0)^{-\frac{1}{2}}$$ which imply
$$R(x,0)\leq \frac{1}{(\frac{1}{2}R(z,0)^{-\frac{1}{2}})^2},\mbox{ on
}B.$$ Then by Li-Yau-Hamilton inequality \cite{Ha4} and the
$\kappa$-noncollapsing, we have $$Vol_0(B)\geq
\kappa(\frac{1}{4}R(z,0)^{-\frac{1}{2}})^4$$ and then
$$Vol_0(B_0(z,8R(z,0)^{-\frac{1}{2}})\geq
\frac{\kappa}{2^{20}}(8R(z,0)^{-\frac{1}{2}})^4.$$ So by Corollary
11.6 of \cite{P1}, there exists a positive constant $A_1$
depending only on $\kappa$ such that $$R(x,0)\leq A_1R(z,0),
\mbox{ for }x\in B_0(z,2R(z,0)^{-\frac{1}{2}}). \eqno (3.1)$$

We now consider the remaining case: $R(x,0)d^2_0(x,y)\leq 1$ for
all $x\in M^4$. We choose a point $z\in M^4$ satisfying
$R(z,0)\geq \frac{1}{2}\sup\{R(x,0)|x\in M^4\}$. Obviously we also
have the estimate (3.1) in the remaining case.

After having the estimate (3.1), we next want to bound $R(z,0)$
for the chosen $z\in M^4$. By combining with Li-Yau-Hamilton
inequality  \cite{Ha4}, we have $$R(x,t)\leq A_1R(z,0),$$ for all
$x\in B_0(z,2R(z,0)^{-\frac{1}{2}})$ and all $t\leq 0$. It then
follows from Shi's local derivative estimate  \cite{Sh1} that
$$\frac{\partial}{\partial t}R(z,t)\leq A_2R(z,0)^2, \ \ \mbox{ for all } -R^{-1}(z,0) \leq t\leq
0,$$ where $A_2$ is some constant depending only on $\kappa$. This
implies $$R(z,-cR^{-1}(z,0))\geq cR(z,0)$$ for some small positive
constant $c$ depending only on $\kappa$. On the other hand, by
using the Harnack estimate  \cite{Ha4} (as a consequence of
Li-Yau-Hamilton inequality), we have $$1=R(y,0)\geq
\widetilde{c}R(z,-cR^{-1}(z,0))$$ for some small positive constant
$\widetilde{c}$ depending only on $\kappa$. Thus we obtain
$$R(z,0)\leq A_3 \eqno (3.2)$$ for some positive constant $A_3$
depending only on $\kappa$.

The combination of (3.1) and (3.2) gives $$R(x,0)\leq A_1A_3,
\mbox{ on }B_0(y,A^{-\frac{1}{2}}_3).$$ Thus by the
$\kappa$-noncollapsing there exists a positive constant $r_0$
depending only on $\kappa$ such that $$Vol_0(B_0(y,r_0))\geq
\kappa r^4_0.$$ For any fixed $R_0\geq r_0$, we have
$$Vol_0(B_0(y,R_0))\geq \kappa(\frac{r_0}{R_0})^4\cdot R^4_0.$$ By
applying Corollary 11.6 of  \cite{P1} again, there exists a
positive constant $\omega(R_0^2)$ depending only on $R_0$ and
$\kappa$ such that $$R(x,0)\leq \omega (R^2_0), \mbox{ on
}B_0(y,\frac{1}{4}R_0).$$ This gives the desired estimate.

(ii) It immediately follows from the above assertion (i), the
Li-Yau-Hamilton inequality \cite{Ha4} and Shi's local derivative
estimates \cite{Sh1}.
$$\eqno \#$$\vskip 0.3cm

$\mathbf{Remark.}$ The argument in the last paragraph of the above
proof for (i) implies the following assertion:

\emph{For any $\zeta>0$, there is a positive function $\omega$
depending only on $\zeta$ such that if there holds
$$\frac{Vol_{t_0}(B_{t_0}(y,{R(y,t_0)}^{-\frac{1}{2}}))}
{{R(y,t_0)}^{-\frac{4}{2}}}\geq \zeta,$$ for some \textbf{fixed}
point $y$ and some $t_0\in(-\infty,0]$, then we have the following
the elliptic type estimate
$$R(x,t_0)\leq R(y,t_0)\cdot \omega(R(y,t_0)d^2_{t_0}(x,y))$$
for all $x\in M.$}\vskip 0.3cm

 This estimate will play a key role in deriving the
universal noncollapsing property in the next subsection.

Let $g_{ij}(t)$, $-\infty<t\leq 0$, be a nonflat solution to the
Ricci flow on a four-manifold $M^4$. Fix a small $\varepsilon>0$. We
say that a point $x_0\in M^4$ is the \textbf{center of an evolving
$\varepsilon$-neck}, if the solution $g_{ij}(t)$ in the set
$\{(x,t)|-\varepsilon^{-2} Q^{-1}<t\leq 0, d^2_0(x,x_0)<\varepsilon
^{-2}Q^{-1}\}$(with $Q=R(x_0,0)$) is, after scaling with factor $Q$,
$\varepsilon$-close (in $C^{[\varepsilon^{-1}]}$ topology) to the
corresponding subset of the evolving round cylinder
$\mathbb{R}\times \mathbb{S}^3$, having scalar curvature one at
$t=0$.

The following result generalizes Corollary 11.8 of Perelman
\cite{P1} to four-dimension and verifies Theorem E 3.3 of Hamilton
\cite{Ha7}. The crucial information in the following Proposition is
that the constant $C=C(\varepsilon)>0$ depends \textbf{\emph{only}}
on $\varepsilon$. \vskip 0.3cm

$\underline{\mbox{\textbf{Proposition 3.4}}}$ \emph{For any
$\varepsilon>0$, there exists $C=C(\varepsilon)>0$ such that if
$g_{ij}(t)$ is a nonflat ancient $\kappa$-solution with restricted
isotropic curvature pinching on a \textbf{noncompact} four-manifold
$M^4$ for some $\kappa>0$, and $M^4_{\varepsilon}$ denotes the set
of points of $M^4$, which are not centers of evolving
$\varepsilon$-necks, then either the whole $M^4$ is a metric
quotient of the round cylinder $\mathbb{R}\times \mathbb{S}^3$ or
$M^4_{\varepsilon}$ satisfies the following properties }

(i) \emph{$M^4_{\varepsilon}$ is compact, and}

(ii) \emph{$diam(M^4_{\varepsilon})\leq CQ^{-\frac{1}{2}}$ and
$C^{-1}Q\leq R(x,0)\leq CQ$, whenever $x\in M^4_{\varepsilon},$
where $Q=R(x_0,0)$ for some $x_0\in \partial M^4_{\varepsilon}$ and
$diam(M^4_{\varepsilon})$ is the diameter of the set
$M^4_{\varepsilon}$ with respect to the metric $g_{ij}(0)$.}\vskip
0.2cm

$\underline{\mbox{\textbf{Proof.}}}$   Note that the curvature
operator of the ancient $\kappa$-solution is nonnegative. We first
consider the easy case that the curvature operator has a nontrivial
null vector somewhere at some time. By Lemma 3.2, we know that the
ancient $\kappa$-solution is a metric quotient of the round cylinder
$\mathbb{R}\times \mathbb{S}^3$.

We then assume the curvature operator of the ancient
$\kappa$-solution is positive everywhere. Firstly we want to show
$M^4_{\varepsilon}$ is compact. Argue by contradiction. Suppose
there exists a sequence of points $z_k,\ k=1,2,\cdots$, going to
infinity (with respect to the metric $g_{ij}(0)$) such that each
$z_k$ is not the center of any evolving $\varepsilon$-neck. For
arbitrarily fixed point $z_0\in M^4$, it follows from Proposition
3.3 (i) that $$0<R(z_0,0)\leq R(z_k,0)\cdot
\omega(R(z_k,0)d^2_0(z_k,z_0))$$ which implies that
$$\lim_{k\rightarrow
\infty}R(z_k,0)d^2_0(z_k,z_0)=+\infty.$$  By Lemma 3.1 and
Proposition 3.3 and Hamilton's compactness theorem, we conclude that
$z_k$ is the center of an evolving $\varepsilon$-neck as $k$
sufficiently large. This is a contradiction, so we have proved that
$M^4_{\varepsilon}$ is compact.

Note that $M^4$ is diffeomorphic to $\mathbb{R}^4$ since the
curvature operator is positive. We may assume $\varepsilon>0$ so
small that Hamilton's replacement for Schoenflies conjecture and
its proof (Theorem G1.1 and Lemma G1.3 of  \cite{Ha7}) are
available. Since every point outside the compact set
$M^4_{\varepsilon}$ is the center of an evolving
$\varepsilon$-neck, it follows that the approximate round
three-sphere cross-section through the center divides $M^4$ into
two connected components such that one of them is diffeomorphic to
the four-ball $\mathbb{B}^4$.  Let $\varphi$ be a Busemann
function on $M^{4}$(constructed from all geodesic rays emanating
from a given point), it is a standard fact that $\varphi$ is
convex and proper. Since $M^{4}_{\varepsilon}$ is compact,
$M^{4}_\varepsilon$ is contained in a compact set
$K=\varphi^{-1}((-\infty,A])$ for some large $A$. We note that
each point $x\in M^{4}\setminus M_\varepsilon$ is the center of an
$\varepsilon$-neck. It is clear that there is an
$\varepsilon$-neck $N$ lying entirely outside $K$.  Consider a
point $x$ on one of its boundary components of the
$\varepsilon$-neck $N$. Since $x \in M^{4}\setminus
M^{4}_{\varepsilon}$, there is an $\varepsilon$-neck adjacent to
the initial $\varepsilon$-neck, producing a longer neck. We then
take a point on the boundary of the second $\varepsilon$-neck and
continue. This procedure can either terminate when we get into
$M_{\varepsilon}$ or go on infinitely to produce a semi-infinite
(topological) cylinder. The same procedure can be repeated for the
other boundary component of the initial $\varepsilon$-neck. This
procedure will give a maximal extended neck $\tilde{N}$. If
$\tilde{N}$ never touch $M^{4}_\varepsilon$, the manifold will be
diffeomorphic to the standard infinite cylinder, which is a
contradiction. If both of the two ends of $\tilde{N}$ touch
$M^{4}_\varepsilon$, then there is a geodesic connecting two
points of $M^{4}_{\varepsilon}$ and passing through $N$. This is
impossible since the function $\varphi$ is convex. So we conclude
that one end of $\tilde{N}$ will touch $M^{4}_\varepsilon$ and the
other end will tend to infinity to produce a semi-infinite
(topological) cylinder. Then one can find an approximate round
three-sphere cross-section which encloses the whole set
$M^{4}_{\varepsilon}$ and touches some point $x_{0}\in
\partial M^{4}_{\varepsilon}$. We now want to show that
$R(x_0,0)^{\frac{1}{2}}\cdot diam(M^{4}_{\varepsilon})$ is bounded
from above by some positive constant $C=C(\varepsilon)$ depending
\textbf{\emph{only}} on $\varepsilon$.

Suppose not, there exist a sequence of nonflat noncompact ancient
$\kappa_j$-solutions with restricted isotropic curvature pinching
and with positive curvature operator, for some sequence of positive
constants $\kappa_j$, such that for above chosen points $x_0\in
M^4_{\varepsilon}$ there would hold
$$R(x_0,0)^{\frac{1}{2}}\cdot diam(M^4_{\varepsilon})\rightarrow +\infty. \eqno
(3.3)$$

Since the point $x_0$ lies in some $2\varepsilon$-neck, clearly,
there is a universal positive lower bound for
${Vol_0(B_0(x_0,\frac{1}{\sqrt{R(x_0,0)}}))}/
({\frac{1}{\sqrt{R(x_0,0)}}})^{4}. $  By the remark after the proof
of the previous Proposition 3.3, we see that there is a universal
positive function $\omega:[0,+\infty)\rightarrow [0,+\infty)$ such
that the elliptic type estimate
$$R(x,0)\leq R(x_0,0)\cdot \omega(R(x_0,0)d^2_0(x,x_0))  \eqno(3.4)$$
holds for all $x\in M.$

 Let us scale the ancient solutions around the points
$x_0$ with the factors $R(x_0,0)$. By (3.4), Hamilton's compactness
theorem (Theorem 16.1 of \cite{Ha6})  and the universal noncollasing
property at $x_0$, we can extract a convergent subsequence. From the
choice of the points $x_0$ and (3.3), the limit contains a line.
Actually we may draw a geodesic ray from some point $x_1\in
M_{\varepsilon}^{4}$ which is far from $x_0$ (in the normalized
distance). This geodesic ray must across some  vertical three-sphere
containing $x_0$. The limit of these rays gives us a line.  Then by
Toponogov splitting theorem the limit is isometric to
$\mathbb{R}\times X^3$ for some smooth three-manifold $X^3$. As
before, by using the restricted isotropic curvature pinching
condition (2.4) and the contracted second Bianchi identity, we see
that $X^3=\mathbb{S}^3/\Gamma$ for some group $\Gamma$ of isometrics
without fixed points. Then we apply the same argument as in the
proof of Lemma 3.1 to conclude that $\Gamma=\{1\}$. This says that
the limit must be the evolving round cylinder $\mathbb{R}\times
\mathbb{S}^3$. This contradicts with the fact that each chosen
points $x_0$ is not the center of any evolving $\varepsilon$-neck.
Therefore we have proved
$$diam(M^4_{\varepsilon})\leq CQ^{-\frac{1}{2}}$$ for some
positive constant $C=C(\varepsilon)$ depending only on
$\varepsilon$, where $Q=R(x_0,0)$.

Finally by combining this diameter estimate and the remark after
proposition 3.3, we directly deduce $$\widetilde{C}^{-1}Q\leq
R(x,0)\leq \widetilde{C}Q,\mbox{ whenever }x\in M^4_{\varepsilon},$$
for some positive constant $\widetilde{C}$ depending only on
$\varepsilon$.
$$\eqno \#$$\vskip 0.3cm

 Consequently, by applying the standard volume comparison to Proposition 3.4,
 we conclude that all complete \textbf{\emph{noncompact}} four-dimensional ancient $\kappa$-solutions
with restricted isotropic curvature pinching and positive curvature
operator are $\kappa_0$-noncollapsing on all scales for some
universal constant $\kappa_0>0$. In the next subsection, we will
prove this universal noncollapsing property for both compact and
noncompact cases.\vskip 0.5cm

\subsection*{3.3 Universal noncollapsing of ancient $\kappa$-solutions}
First we note that the universal noncollapsing is not true for all
metric quotients of round $\mathbb{R}\times \mathbb{S}^{3}.$  The
main result of this section is to establish the universal
noncollapsing property for all ancient $\kappa$-solutions with
restricted isotropic curvature pinching which are not metric
quotients of round $\mathbb{R}\times \mathbb{S}^{3}.$ The analogous
result for three-dimensional ancient $\kappa$-solutions was claimed
by Perelman in Remark 11.9 of \cite{P1} and Section 1.5 of
\cite{P2}.\vskip 0.3cm

$\underline{\mbox{\textbf{Theorem 3.5}}}$ \emph{ There exists a
positive constant $\kappa_0$ with the following property. Suppose we
have a four-dimensional (compact or noncompact) ancient
$\kappa$-solution with restricted isotropic curvature pinching for
some $\kappa>0$. Then either the solution is $\kappa_0$-noncollapsed
for all scales, or it is a metric quotient of the round cylinder
$\mathbb{R}\times \mathbb{S}^3$.}\vskip 0.2cm

$\underline{\mbox{\textbf{Proof.}}}$ Let $g_{ij}(x,t)$, $x\in M^4$
and $t\in (-\infty,0]$, be an ancient $\kappa$-solution with
restricted isotropic curvature pinching for some $\kappa>0$. We had
known that the curvature operator of the solution $g_{ij}(x,t)$ is
nonnegative everywhere and every time. If the curvature operator of
the solution $g_{ij}(x,t)$ has a nontrivial null eigenvector
somewhere at some time, then we know from Lemma 3.2 that the
solution is a metric quotient of the round neck $\mathbb{R}\times
\mathbb{S}^3$.

We now assume the solution $g_{ij}(x,t)$ has positive curvature
operator everywhere and every time. We want to apply backward limit
argument of Perelman to take a sequence points $q_k$
and a sequence of times $t_k \rightarrow -\infty$ such that the
scalings of $g_{ij}(\cdot,t)$ around $q_k$ with factors $|t_k|^{-1}$
(and shifting the times $t_k$ to zero) converge in
$C^{\infty}_{loc}$ topology to a non-flat gradient shrinking
soliton.

Clearly, we may assume the nonflat ancient $\kappa$-solution is not
a gradient shrinking Ricci soliton.  For arbitrary point $(p,t_0)\in
M^4\times (-\infty,0]$, we define as in \cite{P1} that
$$\tau=t_0-t, \mbox{ for }t<t_0,$$
\begin{eqnarray*}\ \ \ \ \ \ l(q,\tau)=\frac{1}{2\sqrt{\tau}}\inf\{&\int^{\tau}_0\sqrt{s}
(R(\gamma(s),t_0-s)+|\dot{\gamma}(s)|^2_{g_{ij}(t_0-s)})ds|\\
  &\gamma:[0,\tau]\rightarrow M^4\mbox{ with
}\gamma(0)=p,\gamma(\tau)=q\},
\end{eqnarray*}
$$\mbox{and}\ \ \  \ \ \widetilde{V}(\tau)=\int_{M^4}(4\pi\tau)^{-2}\exp(-l(q,\tau))dV_{t_0-\tau}(q),\ \ \ \ \ \ \ \
$$
where $|\cdot|_{g_{ij}(t_0-s)}$ is the norm with respect to the
metric $g_{ij}(t_0-s)$ and $dV_{t_0-\tau}$ is the volume element
with respect to the metric $g_{ij}(t_0-\tau)$. According to
\cite{P1}, $l$ is called the \textbf{reduced distance} and
$\widetilde{V}(\tau)$ is called the \textbf{reduced volume}. Since
the manifold $M^4$ may be noncompact, one would ask whether the
reduced volume is finite. Since the scalar curvature is
nonnegative and the curvature is bounded, it is no hard to see
that the reduced distance is quadratically growth and then the
reduced volume is always finite. (Actually, by using Perelman's
Jacobian comparison theorem \cite{P1} one can show that the
reduced volume is always finite for any complete solution of the
Ricci flow (see \cite{CaZ} for the detail)). In \cite{P1},
Perelman proved that the reduced volume $\widetilde{V}(\tau)$ is
nonincreasing in $\tau$, and the monotonicity is strict unless the
solution is a gradient shrinking Ricci soliton.

From \cite{P1} (Section 7 of \cite{P1}), the function
$\overline{L}(q,\tau)=4\tau l(q,\tau)$ satisfies
$$\frac{\partial}{\partial\tau}\overline{L}+\triangle\overline{L}\leq 8.$$
It is clear that $\overline{L}(\cdot,\tau)$ achieves its minimum on
$M^4$ for each $\tau>0$ since the scalar curvature is nonnegative.
Then the minimum of $\overline{L}(\cdot,\tau)-8\tau$ is
nonincreasing, so in particular, the minimum of $l(\cdot,
\overline{\tau})$ does not exceed $2$ for each $\tau>0$. Thus for
each $\tau>0$ we can find $q=q(\tau)$ such that $l(q(\tau),\tau)\leq
2$. We can apply Perelman's Proposition 11.2 in \cite{P1} to
conclude that the scalings of $g_{ij}(\cdot,t_0-\tau)$ around
$q(\tau)$ with factors $\tau^{-1}$ converge in $C^{\infty}_{loc}$
topology along a subsequence $\tau\rightarrow+\infty$ to a non-flat
gradient shrinking soliton. Because the proof of this proposition in
\cite{P1} is just a sketch, we would like to give its detail in the
following for the completeness.

We first claim that for any $A \geq 1$, one can find
$B=B(A)<+\infty$ such that for every $\overline{\tau}>1$ there
holds
$$l(q,\tau)<B\mbox{ and }\tau R(q,t_0-\tau)\leq B \eqno (3.5)$$
whenever $\frac{1}{2}\overline{\tau}\leq \tau\leq
A\overline{\tau}$ and
$d^2_{t_0-\frac{\overline{\tau}}{2}}(q,q(\frac{\overline{\tau}}{2}))\leq
A\overline{\tau}$.

Indeed, by Section 7 of \cite{P1}, the reduced distance $l$
satisfies the following $$\left\{
\begin{array}{lll}
\frac{\partial}{\partial
\tau}l=-\frac{l}{\tau}+R+\frac{K}{2\tau^{3/2}},
\ \ \ \ \ \ \ \ \ \ \ \ \ \ \ \ \ \ \ \ \ \ \ \ \ \ \ \ \ \ \ \ \ \ \ \ \ (3.6)\\
|\nabla l|^2=-R+\frac{l}{\tau}-\frac{K}{\tau^{3/2}}, \ \ \ \ \ \ \ \ \ \ \ \ \ \ \ \ \ \ \ \ \ \ \ \ \ \ \ \ \ \ \ \ \ \ \ (3.7)\\
\triangle l\leq -R+\frac{2}{\tau}-\frac{K}{2\tau^{3/2}}, \ \ \ \ \
\ \ \ \ \ \ \ \ \ \ \ \ \ \ \ \ \ \ \ \ \ \ \ \ \ \ \ \ \ \ \ \
(3.8)\end{array}\right.$$ in the sense of distributions, and the
equality holds everywhere if and only if we are on a gradient
shrinking soliton, where $K=\int^{\tau}_0s^{3/2}Q(X)ds$ and $Q(X)$
is the trace Li-Yau-Hamilton quadratic given by
$$Q(X)=-\frac{\partial}{\partial \tau}R-\frac{R}{\tau}-2<\nabla
R,X>+2Ric(X,X)$$ and $X$ is the tangential (velocity) vector of a
$\mathcal{L}$-shortest curve $\gamma:[0,\tau]\rightarrow M^4$
connecting $p$ to $q$.

By applying the trace Li-Yau-Hamilton inequality \cite{Ha4} to the
ancient $\kappa$-solution, we have $$Q(X)\geq -\frac{R}{\tau}$$
and then
$$K\geq -\int^{\tau}_0\sqrt{s}Rds\geq -2\sqrt{\tau}l.$$ Thus by
(3.7) we get $$|\nabla l|^2+R\leq \frac{3l}{\tau}. \eqno (3.9)$$
At $\tau=\frac{\bar{\tau}}{2}$, we have
\begin{equation}\tag{3.10}
\begin{split}
 \sqrt{l (q,\frac{\overline{\tau}}{2})}&\leq
\sqrt{2}+\sup\{|\nabla \sqrt{l}|\}\cdot
d_{t_0-\frac{\overline{\tau}}{2}}(q,q(\frac{\overline{\tau}}{2})) \\
&\leq \sqrt{2}+\sqrt{\frac{3A}{2}},
\end{split}
 \end{equation}
 and
$$R(q,t_0-\frac{\overline{\tau}}{2})\leq \frac{6}{\overline{\tau}}(\sqrt{2}+\sqrt{\frac{3A}{2}})^2, \eqno
(3.11)$$ for $q\in
B_{t_0-\frac{\overline{\tau}}{2}}(q(\frac{\overline{\tau}}{2}),\sqrt{A\overline{\tau}})$.
Since the scalar curvature of an ancient solution with nonnegative
curvature operator is pointwisely nondecreasing in time (by the
trace Li-Yau-Hamilton inequality \cite{Ha4}), we further have
$$\tau R(q,t_0-\tau)\leq 6A(\sqrt{2}+\sqrt{\frac{3A}{2}})^2 \eqno
(3.12)$$ whenever $\frac{1}{2}\overline{\tau}\leq \tau\leq
A\overline{\tau}$ and
$d^2_{t_0-\frac{\overline{\tau}}{2}}(q,q(\frac{\overline{\tau}}{2}))\leq
A\overline{\tau}$.

By (3.6), (3.7) and (3.12), we have $$\frac{\partial}{\partial
\tau}l\leq
-\frac{l}{2\tau}+\frac{3A}{\tau}(\sqrt{2}+\sqrt{\frac{3A}{2}})^2$$
and by integrating this inequality and using the estimate (3.10),
we obtain $$l (q,\tau)\leq 7A(\sqrt{2}+\sqrt{\frac{3A}{2}})^2
\eqno (3.13)$$ whenever $\frac{1}{2}\overline{\tau}\leq \tau\leq
A\overline{\tau}$ and
$d^2_{t_0-\frac{\overline{\tau}}{2}}(q,q(\frac{\overline{\tau}}{2}))\leq
A\overline{\tau}$. So we have proved the assertion (3.5).

The scaling of the ancient $\kappa$-solution around
$q(\frac{\overline{\tau}}{2})$ with factor
$(\frac{\overline{\tau}}{2})^{-1}$ is
$$\widetilde{g}_{ij}(s)=\frac{2}{\overline{\tau}}g_{ij}(\cdot,t_0-s\frac{\overline{\tau}}{2})$$
for $s\in [0,+\infty)$. The assertion (3.5) implies that for all
$s\in [1,2A]$ and all $q$ with
$dist^2_{\widetilde{g}_{ij}(1)}(q,q(\frac{\overline{\tau}}{2}))\leq
A$, we have $\widetilde{R}(q,s)\leq B$ where $\widetilde{R}$ is
the scalar curvature of the rescaled metric $\widetilde{g}_{ij}$.
Then we can use Hamilton's compactness theorem (\cite{Ha5} or more
precisely Theorem 16.1 of \cite{Ha6}) and the
$\kappa$-noncollapsing assumption to obtain a sequence
$\overline{\tau}_k\rightarrow +\infty$ such that the marked
evolving manifolds
($M^4,\widetilde{g}^{(k)}_{ij}(s),q(\frac{\overline{\tau}_k}{2})$),
with
$\widetilde{g}^{(k)}_{ij}(s)=\frac{2}{\overline{\tau}_k}g_{ij}(\cdot,t_0-s\frac{\overline{\tau}_k}{2})$
and $s \in [1,+\infty)$, converge in $C^{\infty}_{loc}$ topology
to an evolving manifold
$(\overline{M}^4,\overline{g}_{ij}(s),\overline{q})$ with $s \in
[1,+\infty)$, where $\overline{g}_{ij}(s)$ satisfies
$\frac{\partial}{\partial s}\overline{g}_{ij}=2\overline{R}_{ij}$
on $\overline{M}\times [1,+\infty)$.

Denote by $\widetilde{l}_k$ the corresponding reduced distance of
$\widetilde{g}^{(k)}_{ij}(s)$. It is easy to see that
$\widetilde{l}_k(q,s)=l(q,\frac{\overline{\tau}_k}{2}s)$ for $s\in
[1,+\infty)$. After rescaling we still have $$|\nabla
\widetilde{l}_k|^2_{\widetilde{g}^{(k)}_{ij}}+\widetilde{R}^{(k)}\leq
6\widetilde{l}_k$$ and by (3.5), $\widetilde{l}_k$ are uniformly
bounded at finite distances. Thus the above gradient estimate
implies that the functions $\widetilde{l}_k$ tend (up to a
subsequence) to a function $\overline{l}$ which is a locally
Lipschitz function on $\overline{M}$.

From (3.6)-(3.8), we have $$\frac{\partial}{\partial
s}(\widetilde{l}_k)-\triangle
\widetilde{l}_k+|\nabla\widetilde{l}_k|^2-\widetilde{R}^{(k)}+\frac{2}{s}\geq
0,$$ $$2\triangle
\widetilde{l}_k-|\nabla\widetilde{l}_k|^2+\widetilde{R}^{(k)}+\frac{\widetilde{l}_k-4}{s}\leq
0,$$ which can be rewritten as $$(\frac{\partial}{\partial
s}-\triangle+\widetilde{R}^{(k)})((4\pi
s)^{-2}\exp(-\widetilde{l}_k))\leq 0, \eqno (3.14)$$
$$-(4\triangle -\widetilde{R}^{(k)})e^{-\frac{\widetilde{l}_k}{2}}+\frac{\tilde{l}_{k}-4}{s}e^{-\frac{\widetilde{l}_k}{2}}\leq 0, \eqno
(3.15)$$ in the sense of distribution. Clearly, these two
inequalities imply that the limit $\overline{l}$ satisfies
$$(\frac{\partial}{\partial s}-\triangle +\overline{R})((4\pi s)^{-2}\exp(-\overline{l}))\leq 0, \eqno
(3.16)$$ $$-(4\triangle
-\overline{R})e^{-\frac{\overline{l}}{2}}+\frac{\bar{l}-4}{s}e^{-\frac{\overline{l}}{2}}\leq
0, \eqno (3.17)$$ in the sense of distributions.

Denote by $\widetilde{V}^{(k)}(s)$ the reduced volume of the
rescaled metric $\widetilde{g}^{(k)}_{ij}(s)$. Since
$\widetilde{l}_k(q,s)=l(q,\frac{\overline{\tau}_k}{2}s)$, we see
that
$\widetilde{V}^{(k)}(s)=\widetilde{V}(\frac{\overline{\tau}_k}{2}s)$.
The monotonicity of the reduced volume $\widetilde{V}(\tau)$ (see
\cite{P1}) then implies that
$$\lim_{k\rightarrow+\infty}\widetilde{V}^{(k)}(s)=\overline{V}, \mbox{  for
}s\in[1,2]$$ for some positive constant $\overline{V}$. But we are
not sure whether the limiting $\overline{V}$ is exactly the
Perelman's reduced volume of the limiting manifold
($\overline{M}^4,\overline{g}_{ij}(s)$), because the points
$q(\frac{\overline{\tau}_k}{2})$ may diverge to infinity.
Nevertheless, we can insure that $\overline{V}$ is not less than
the Perelman's reduced volume of the limit.  Note that
$$\widetilde{V}^{(k)}(2)-\widetilde{V}^{(k)}(1)=\int^2_{1}\frac{d}{ds}(\widetilde{V}^{(k)}(s))ds$$
$$=\int^2_{1}ds\int_{M^4}(\frac{\partial}{\partial s}-\triangle+\widetilde{R}^{(k)})((4\pi
s)^{-2}\exp(-\widetilde{l}_k))dV_{\widetilde{g}^{(k)}_{ij}(s)}.$$
Thus we deduce that in the sense of distributions,
$$(\frac{\partial}{\partial s}-\triangle +\overline{R})((4\pi s)^{-2}\exp(-\overline{l}))=0, \eqno
(3.18)$$ $$-(4\triangle
-\overline{R})e^{-\frac{\overline{l}}{2}}+\frac{\bar{l}-4}{s}e^{-\frac{\overline{l}}{2}}=0,
\eqno (3.19)$$ and then the standard parabolic equation theory
implies that $\overline{l}$ is actually smooth. Here we used
(3.6)-(3.8) to show that the equality in (3.16) implies the
equality in (3.17).

Set $$\upsilon=[s(2\triangle \overline{l}-|\nabla
\bar{l}|^2+\overline{R})+\overline{l}-4]\cdot(4\pi
s)^{-2}e^{-\overline{l}}.$$
 A direct computation gives
$$(\frac{\partial}{\partial s}
-\triangle+\overline{R})\upsilon=-2s|\overline{R}_{ij}+\nabla_i\nabla_j\overline{l}-\frac{1}{2s}\overline{g}_{ij}|^2\cdot(4\pi s)^{-2}e^{-\overline{\ell}}. \eqno
(3.20)$$ Since the equation (3.18) implies $\upsilon\equiv0$, the
limit metric $\overline{g}_{ij}$ satisfies
$$\overline{R}_{ij}+\nabla_i\nabla_j\overline{l}-\frac{1}{2s}\overline{g}_{ij}=0. \eqno
(3.21)$$ Thus the limit is a gradient shrinking Ricci soliton.

To show the limiting gradient shrinking Ricci soliton to be
nonflat, we first show that constant $\overline{V}$ is strictly
less than $1$. Indeed, by considering the reduced volume
$\widetilde{V}(\tau)$ of the ancient $\kappa$-solution, we get
from Perelman's Jacobian comparison theorem \cite{P1} that
\begin{eqnarray*}
\widetilde{V}(\tau)&=&\int_{M^4}(4\pi
\tau)^{-2}e^{-l}dV_{t_0-\tau}\\
&\leq& \int_{T_pM^4}(4\pi)^{-2}e^{-|X|^2}dX
\\&=&1. \end{eqnarray*}
Recall that we assumed the nonflat ancient $\kappa$-solution is not
a gradient shrinking Ricci soliton. Thus by the monotonicity of the
reduced volume \cite{P1}, we have $\widetilde{V}(\tau)<1$ for
$\tau>0$. This implies that $\overline{V}<1$.

We now argue by contradiction. Suppose the limit
$\overline{g}_{ij}(s)$ is flat. Then by (3.21) we have
$$\nabla_i\nabla_j\overline{l}=\frac{1}{2s}\overline{g}_{ij}\ \ \mbox{ and
}\ \ \ \triangle\overline{l}=\frac{2}{s}.$$ And then by (3.19), we
get $$|\nabla\overline{l}|^2=\frac{\overline{l}}{s}.$$ Since the
function $\overline{l}$ is strictly convex, it follows that
$\sqrt{4s\overline{l}}$ is a distance function (from some point)
on the complete flat manifold $\overline{M}$. From the smoothness
of the function $\overline{l}$, we conclude that the flat manifold
$\overline{M}$ must be $\mathbb{R}^4$. In this case we would have
its reduced distance to be $\bar{l}$ and its reduced volume to be
$1$. Since $\overline{V}$ is not less than the reduced volume of
the limit, this is a contradiction. Therefore the limiting
gradient shrinking soliton $\overline{g}_{ij}$ is nonflat.

Now we consider the nonflat gradient shrinking Ricci soliton
$(\overline{M}^4,\overline{g}_{ij})$. Of course it is still
$\kappa$-noncollapsed for all scales and satisfies the restricted
isotropic curvature pinching condition (2.4). We first show that
$(\overline{M}^4,\overline{g}_{ij}(s))$ has bounded curvature at
each $s>0$. Clearly it suffices to consider $s=1$. By Lemma 3.2,
we may assume the soliton $(\overline{M}^4,\overline{g}_{ij}(1))$
has positive curvature operator everywhere. Let us argue by
contradiction. Suppose not, then we claim that for each positive
integer $k$, there exists a point $x_k$ such that
\begin{equation*}
 \left\{
\begin{split}
\bar{R}(x_{k},1)&\geq k,\\
\bar{R}(x,1)& \leq4\bar{R}(x_{k},1), \mbox{ for } x\in
B_{\bar{g}(\cdot,1)}(x_k,\frac{k}{\sqrt{\bar{R}(x_{k},1)}}).
\end{split}
\right.
\end{equation*}
Indeed, $x_k$ can be constructed as a limit of a finite sequence
$\{y_i\}$, defined as follows. Let $y_0$ be any fixed point with
$\bar{R}(y_{0},1) \geq k$. Inductively, if $y_i$ cannot be taken
as $x_k$, then there is a $y_{i+1}$ such that
\begin{equation*}
 \left\{
\begin{split}
\bar{R}(y_{i+1},1)&> 4\bar{R}(y_{i},1),\\
d_{\bar{g}(\cdot,1)}(y_i,y_{i+1}) &\leq
\frac{k}{\sqrt{\bar{R}(y_{i},1)}}.
\end{split}
\right.
\end{equation*}
Thus we have
$$\bar{R}(y_{i},1)> 4^i\bar{R}(y_{0},1)\geq 4^ik,$$
$$d_{\bar{g}(\cdot,1)}(y_i,y_{0})\leq
k\sum^i_{j=1}\frac{1}{\sqrt{4^{j-1}k}}<2\sqrt{k}.$$ Since the
soliton is smooth, the sequence $\{y_i\}$ must be finite. The last
element fits.

Note that the limiting soliton still satisfies the Li-Yau-Hamilton
inequality. Then we have
$$\bar{R}(x,s) \leq \bar{R}(x,1) \leq 4\bar{R}(x_k,1)$$
for $x\in
B_{\bar{g}(\cdot,1)}(x_k,\frac{k}{\sqrt{\bar{R}(x_{k},1)}})$ and
$1 \leq s \leq 1 + \frac{1}{\bar{R}(x_{k},1)}.$ By the
$\kappa$-noncollapsing and the Hamilton's compactness theorem
\cite{Ha5}, a sequence of
$(\bar{M}^4,\bar{R}(x_{k},1)\bar{g}(\cdot,1+\frac{(\cdot)}{\bar{R}(x_{k},1)},x_k)$
will converge to a complete smooth solution
$(\bar{\bar{M}}^4,\bar{\bar{g}})$ at least on the interval
$[0,1)$. Since $d_{\bar{g}(\cdot,1)}(x_k,x_{0})\rightarrow \infty$
and $\bar{R}(x_{k},1)\rightarrow \infty$, it follows from Lemma
3.1 that $\bar{\bar{M}}^4 = \mathbb{R}\times \mathbb{S}^3$. This
contradicts with Proposition 2.2. So we have proved that
$(\overline{M}^4,\overline{g}_{ij}(s))$ has bounded curvature at
each $s>0$.

We next show that the soliton $(\overline{M}^4,\overline{g}_{ij})$
is $\kappa'_0$-noncollapsed on all scales for some universal
positive constant $\kappa'_0$. If the soliton
$(\overline{M}^4,\overline{g}_{ij})$ has positive curvature
operator, we know from Hamilton's result \cite{Ha2} and Proposition
3.4 that either the soliton $(\overline{M}^4,\overline{g}_{ij})$ is
the round $\mathbb{S}^4$ or $\mathbb{RP}^4$ when it is compact, or
it is $\kappa'_0$-noncollapsed for all scales for some universal
positive constant $\kappa'_0$ when it is noncompact. (Furthermore,
when the soliton $(\overline{M}^4,\overline{g}_{ij})$ is the round
$\mathbb{S}^4$ or $\mathbb{RP}^4$, it follows from Hamilton's
pinching estimates in \cite{Ha2} that the original ancient
$\kappa$-solution $({M}^4,{g}_{ij}(t))$ is also the round
$\mathbb{S}^4$ or $\mathbb{RP}^4$). While if the soliton
$(\overline{M}^4,\overline{g}_{ij})$ has a nontrivial null
eigenvector somewhere at some time, we know from Lemma 3.2 that the
soliton $(\overline{M}^4,\overline{g}_{ij})$ is $\mathbb{R}\times
\mathbb{S}^{3}/\Gamma$, a metric quotient of the round neck
$\mathbb{R}\times \mathbb{S}^3$. For each $\sigma \in \Gamma$,
$(s,x)\in \mathbb{R}\times \mathbb{S}^{3}$, write
$\sigma(s,x)=(\sigma_{1}(s,x),\sigma_{2}(s,x))\in \mathbb{R}\times
\mathbb{S}^{3}.$ Since $\sigma$ sends lines to lines, and $\sigma$
sends cross spheres to cross spheres, we have
$\sigma_{1}(s,x)=\sigma_{1}(s,y)$,$ \forall\   x,y\in
\mathbb{S}^{3}$. This says that $\sigma_{1}$ reduces to a function
of $s$ alone on $\mathbb{R}$. Moreover, for any $(s,x),(s',x')\in
\mathbb{R}\times \mathbb{S}^{3}$, since $\sigma$ preserves the
distances between cross spheres $\{s\}\times \mathbb{S}^{3}$ and
$\{s'\}\times \mathbb{S}^{3}$, we have
$|\sigma_{1}(s,x)-\sigma_{1}(s',x')|=|s-s'|$. So the projection
 $\Gamma_{1}$ of $\Gamma$ to the factor $\mathbb{R}$ is an isometric
 subgroup of $\mathbb{R}$. We know that if $(\bar{M}^{4},\bar{g}_{ij})=\mathbb{R}\times \mathbb{S}^{3}/\Gamma$
  was compact, it, as an ancient solution, could not be $\kappa$-noncollapsed on all scale
  as $t\rightarrow -\infty$. Thus $(\bar{M}^{4},\bar{g}_{ij})=\mathbb{R}\times
  \mathbb{S}^{3}/\Gamma$ is noncompact. It follows that $\Gamma_{1}=\{1\}$ or
  $\mathbb{Z}_{2}$. We conclude that, in both cases, there is a $\Gamma$-invariant cross
  sphere $\mathbb{S}^{3}$ in $\mathbb{R}\times \mathbb{S}^{3}$. Denote
  it by $\{0\}\times \mathbb{S}^{3}$. $\Gamma$ acts on $\{0\}\times \mathbb{S}^{3}$
  without fixed points.
Recall that we have assumed that the ancient solution $(M^4,g_{ij})$
has positive curvature operator. Then we apply Hamilton's argument
in Theorem C4.1 of  \cite{Ha7} when $M^4$ is compact and apply the
modified argument in the proof of Lemma 3.1 when $M^4$ is noncompact
to conclude that $(\{0\}\times \mathbb{S}^3)/\Gamma$ is
incompressible in $M^4$ (i.e., $\pi_1((\{0\}\times
\mathbb{S}^3)/\Gamma)$ injects into $\pi_1(M^4)$). By Synge theorem
and the Soul theorem \cite{CE}, the fundamental group $\pi_1 (M^4)$
is either $\{1\}$ or $\mathbb{Z}_2$. This implies that $\Gamma$ is
either $\{1\}$ or $\mathbb{Z}_2$. Thus the limiting soliton
$(\overline{M}^4,\overline{g}_{ij})$ is also
$\kappa'_0$-noncollapsed on all scales for some universal positive
constant $\kappa'_0$.

We next use the $\kappa'_0$-noncollapsing of the limiting soliton
to derive a $\kappa_0$-noncollapsing for the original ancient
$\kappa$-solution. By rescaling, we may assume that $R(x,t)\leq 1$
for all $(x,t)$ satisfying $d_{t_0}(x,p)\leq 2$ and $t_0-1\leq
t\leq t_0$. We only need to bound the volume
$Vol_{t_0}(B_{t_0}(p,1))$ from below by a universal positive
constant.

Denote by $\epsilon=Vol_{t_0}(B_{t_0}(p,1))^{\frac{1}{4}}$. For any
$\upsilon\in T_pM^4$, it was known from  \cite{P1} that one can find
a $\mathcal{L}$-geodesic $\gamma(\tau)$, starting at $p$, with
$\lim_{\tau\rightarrow 0^+}\sqrt{\tau}\dot{\gamma}(\tau)=\upsilon$,
which satisfies the following $\mathcal{L}$-geodesic equation
$$\frac{d}{d\tau}(\sqrt{\tau}\dot{\gamma})-\frac{1}{2}\sqrt{\tau}\nabla R+2Ric(\sqrt{\tau}\dot{\gamma},\cdot)=0. \eqno
(3.22)$$ Note from Shi's local derivative estimate (see \cite{Sh1})
that $|\nabla R|$ is also uniformly bounded. By integrating the
$\mathcal{L}$-geodesic equation we see that as $\tau\leq \epsilon$
with the property that $\gamma(\sigma)\in B_{t_0}(p,1)$ for
$\sigma\in (0,\tau]$, there holds
$$|\sqrt{\tau}\dot{\gamma}(\tau)-\upsilon|\leq C\epsilon(|\upsilon|+1) \eqno
(3.23)$$ for some universal positive constant $C$. Here we
implicitly used the fact that the metrics $g_{ij}(t)$ are
 equivalent to each other on $B_{t_0}(p,1)\times [t_0-1,t_0]$,
which is a easy consequence of the boundedness of the curvature
there. Without loss of generally, we may assume $C\epsilon\leq
\frac{1}{4}$ and $\epsilon\leq \frac{1}{100}$. Then for $\upsilon\in
T_pM^4$ with $|\upsilon|\leq \frac{1}{4}\epsilon^{-\frac{1}{2}}$ and
for $\tau\leq \epsilon$ with the property that $\gamma(\sigma)\in
B_{t_0}(p,1)$ for $\sigma\in (0,\tau]$, we have
\begin{eqnarray*}
d_{t_0}(p,\gamma(\tau))&\leq&
\int^{\tau}_0|\dot{\gamma}(\sigma)|d\sigma\\
&<&\frac{1}{2}\epsilon^{-\frac{1}{2}}\int^{\tau}_0\frac{d\sigma}{\sqrt{\sigma}}\\
&=&1.\end{eqnarray*}
 This shows
$$\mathcal{L}\exp\{|\upsilon|\leq
\frac{1}{4}\epsilon^{-\frac{1}{2}}\}(\epsilon)\subset B_{t_0}(p,1)
\eqno (3.24)$$ where $\mathcal{L}\exp(\cdot)(\epsilon)$ denotes the
exponential map of the $\mathcal{L}$ distance with parameter
$\epsilon$ (see \cite{P1} or \cite{CaZ} for details). We decompose
the reduced volume $\widetilde{V}(\epsilon)$ as
\begin{equation}\tag{3.25}
\begin{split}
\widetilde{V}(\epsilon)&=\int_{M^4}(4\pi
\epsilon)^{-2}\exp(-l)dV_{t_0-\epsilon}\\
&\leq \int_{\mathcal{L}\exp\{|\upsilon|\leq
\frac{1}{4}\epsilon^{-\frac{1}{2}}\}(\epsilon)}+\int_{M^4\setminus\mathcal{L}\exp\{|\upsilon|\leq
\frac{1}{4}\epsilon^{-\frac{1}{2}}\}(\epsilon)}(4\pi
\epsilon)^{-2}\exp(-l)dV_{t_0-\epsilon}.
 \end{split}
\end{equation}
The first term on RHS of (3.25) can be estimated by
\begin{equation}\tag{3.26}
\begin{split}
&\ \int_{\mathcal{L}\exp\{|\upsilon|\leq
\frac{1}{4}\epsilon^{-\frac{1}{2}}\}(\epsilon)}(4\pi\epsilon)^{-2}\exp(-l)dV_{t_0-\epsilon}\\
&\leq e^{4\epsilon}\int_{B_{t_0}(p,1)}(4\pi
\epsilon)^{-2}\exp(-l)dV_{t_0}\\
&\leq e^{4\epsilon}(4\pi)^{-2}\epsilon^{-2}Vol_{t_0}(B_{t_0}(p,1)) \\
&=e^{4\epsilon}(4\pi)^{-2}{\epsilon}^{2}.
\end{split}
\end{equation}
 where we used (3.24)
and the equivalence of the evolving metric over $B_{t_0}(p,1)$.
While the second term on the RHS of (3.25) can be estimated as
follows \begin{equation}\tag{3.27}
\begin{split}
& \int_{M^4\setminus\mathcal{L}\exp\{|\upsilon|\leq
\frac{1}{4}\epsilon^{-\frac{1}{2}}\}(\epsilon)}(4\pi\epsilon)^{-2}\exp(-l)dV_{t_0-\epsilon}
\\
&\leq \int_{\{|\upsilon|>
\frac{1}{4}\epsilon^{-\frac{1}{2}}\}}(4\pi\tau)^{-2}\exp(-l)J(\tau)|_{\tau=0}d\upsilon
\end{split}
 \end{equation} by Perelman's Jacobian
comparison theorem  \cite{P1}, where $J(\tau)$ is the Jacobian of
the $\mathcal{L}$-exponential map.

For any $\upsilon\in T_pM$, we consider a $\mathcal{L}$-geodesic
$\gamma(\tau)$ starting at $p$ with $\lim_{\tau\rightarrow
0^+}\sqrt{\tau}\dot{\gamma}(\tau)=\upsilon$. To evaluate the
Jacobian of the $\mathcal{L}$ exponential map at $\tau=0$ we
choose linear independent vectors $\upsilon_1,\cdots,\upsilon_4$
in $T_pM$ and let
$$V_i(\tau)=(\mathcal{L}\exp_{\upsilon}(\tau))_*(\upsilon_i
)=\frac{d}{ds}|_{s=0}\mathcal{L}\exp_{(\upsilon+s\upsilon_i)}(\tau),
\ i=1,\cdots,4.$$ The $\mathcal{L}$-Jacobian $J(\tau)$ is given by
$$J(\tau)=|V_1(\tau)\wedge\cdots\wedge V_4(\tau)|_{g_{ij}(\tau)}/|v_1\wedge\cdots\wedge
v_4|.$$ By the $\mathcal{L}$-geodesic equation (3.22) and the
deriving of (3.23), we see that as $\tau>0$ small enough,
$$|\sqrt{\tau}\frac{d}{d\tau}\mathcal{L}\exp_{(\upsilon+s\upsilon_i)}(\tau)-(\upsilon+s\upsilon_i)|\leq
o(1)$$ for $s\in (-\epsilon,\epsilon)$ and $i=1,\cdots,4$, where
$o(1)$ tends to zero as $\tau\rightarrow0^+$ uniformly in $s$. This
implies that $$\lim_{\tau\rightarrow
0^+}\sqrt{\tau}\dot{V}_i(\tau)=\upsilon_i,\  i=1,\cdots,4,$$ so we
get $$\lim_{\tau\rightarrow 0^+}\tau^{-2}J(\tau)=1. \eqno (3.28)$$

To evaluate $l(\cdot,\tau)$ at $\tau=0$, we use (3.23) again to
get
$$l(\cdot,\tau)=\frac{1}{2\sqrt{\tau}}\int^{\tau}_0\sqrt{s}(R+|\dot{\gamma}(s)|^2)ds$$
$$\rightarrow|\upsilon|^2, \ \ \ \mbox{ as }\tau\rightarrow 0^+,$$ thus
$$l(\cdot,0)=|\upsilon|^2. \eqno (3.29)$$ Hence by combining
(3.27)-(3.29) we have \begin{equation}\tag{3.30}
\begin{split}
&\int_{M^4\setminus\mathcal{L}\exp\{|\upsilon|\leq\frac{1}{4}\epsilon^{-\frac{1}{2}}\}(\epsilon)}(4\pi
\epsilon)^{-2}\exp(-l)dV_{t_0-\epsilon} \\ &\leq
(4\pi)^{-2}\int_{\{|\upsilon|>\frac{1}{4}\epsilon^{-\frac{1}{2}}\}}\exp(-|\upsilon|^2)d\upsilon\\
&<\epsilon^2. \end{split}
\end{equation}
By summing up (3.25), (3.26) and (3.30), we obtain
$$\widetilde{V}(\epsilon)<2\epsilon^2. \eqno (3.31)$$

On the other hand we recall that there are sequences
$\tau_k\rightarrow+\infty$ and $q(\tau_k)\in M^4$ with
$l(q(\tau_k),\tau_k)\leq 2$ so that the rescalings of the ancient
$\kappa$-solution around $q(\tau_k)$ with factor $\tau^{-1}_k$
converge to a gradient shrinking Ricci soliton which is
$\kappa'_0$-noncollapsing on all scales for some universal
positive constant $\kappa'_0$. For sufficiently large $k$, we
construct a path $\gamma:[0,2\tau_k]\rightarrow M^{4}$, connecting
$p$ to any given point $q\in M^{4}$, as follows: the first half
path $\gamma|_{[0,\tau_k]}$ connects $p$ to $q(\tau_k)$ such that
$$l(q(\tau_k),\tau_k)=\frac{1}{2\sqrt{\tau_k}}\int^{\tau_k}_0\sqrt{\tau}(R+|\dot{\gamma}(\tau)|^2)d\tau\leq
3,$$ and the second half path $\gamma|_{[\tau_k,2\tau_k]}$ is a
shortest geodesic connecting $q(\tau_k)$ to $q$ with respect to
the metric $g_{ij}(t_0-\tau_k)$. Note that the rescaled metric
$\tau^{-1}_kg_{ij}(t_0-\tau)$ over the domain
$B_{t_0-\tau_k}(q(\tau_k),\sqrt{\tau_k})\times
[t_0-2\tau_k,t_0-\tau_k]$ is sufficiently close to the gradient
shrinking Ricci soliton. Then by the estimates (3.5) and the
$\kappa'_0$-noncollapsing of the shrinking soliton, we get
\begin{eqnarray*}
\widetilde{V}(2\tau_k)&=&\int_M(4\pi(2\tau_k))^{-2}\exp(-l(q,2\tau_k))dV_{t_0-2\tau_k}(q)\\
&\geq&
\int_{B_{t_0-\tau_k}(q(\tau_k),\sqrt{\tau_k})}(4\pi(2\tau_k))^{-2}\exp(-l(q,2\tau_k))dV_{t_0-2\tau_k}(q)\\
&\geq& \beta
\end{eqnarray*} for some universal positive constant $\beta$. By applying the
monotonicity of the reduced volume  \cite{P1} and (3.31), we
deduce that
$$\beta\leq \widetilde{V}(2\tau_k)\leq
\widetilde{V}(\epsilon)<2\epsilon^2.$$ This proves
$$Vol_{t_0}(B_{t_0}(p,1))\geq \kappa_0>0$$ for some universal
positive constant $\kappa_0$. Therefore we have proved the
Theorem.
$$\eqno \#$$ \vskip 0.3cm

Once the universal noncollapsing of ancient $\kappa$-solution with
restricted isotropic curvature pinching is established, we can also
strengthen the elliptic type estimates in Proposition 3.3 to the
following form.\vskip 0.3cm

$\underline{\mbox{\textbf{Proposition 3.6}}}$ \emph{ There exist a
positive constant $\eta$ and a positive function
$\omega:[0,+\infty)\rightarrow(0,+\infty)$ with the following
properties. Suppose we have a four-dimensional ancient
$\kappa$-solution $(M^4,g_{ij}(t)),-\infty<t\leq 0$, with restricted
isotropic curvature pinching for some $\kappa>0$. Then}

(i) \emph{for every $x,y\in M^4$ and $t\in (-\infty,0]$, there
holds
$$R(x,t)\leq R(y,t)\cdot \omega(R(y,t)d^2_t(x,y));$$}
 \ \ (ii) \emph{for all $x\in M^4$ and $t\in (-\infty,0]$, there hold
$$|\nabla R|(x,t)\leq \eta R^{\frac{3}{2}}(x,t) \mbox{ and }|\frac{\partial R}{\partial t}|(x,t)\leq \eta
R^2(x,t).$$ }\vskip 0.3cm

The following result generalizes Theorem 11.7 of Perelman \cite{P1}
to four-dimension.\vskip 0.3cm

$\underline{\mbox{\textbf{Corollary 3.7}}}$ \emph{The set of
four-dimensional ancient $\kappa$-solutions with restricted
isotropic curvature pinching and \textbf{positive} curvature
operator is precompact modulo scaling in the sense that for any
sequence of such solutions and marked points $(x_k,0)$ with
$R(x_k,0)=1$, we can extract a $C^{\infty}_{loc}$ converging
subsequence, and the limit is also an ancient $\kappa_0$-solution
with restricted isotropic curvature pinching.}\vskip 0.2cm

$\underline{\mbox{\textbf{Proof}}}$: Consider any sequence of
four-dimensional ancient $\kappa$-solutions with restricted
isotropic curvature pinching and positive curvature operator and
marked points $(x_k,0)$ with $R(x_k,0)=1$. By the  Proposition 3.6
(i), Li-Yau-Hamilton inequality \cite{Ha4} and Hamilton's
compactness theorem (Theorem 16.1 of \cite{Ha6}), we can extract a
$C^{\infty}_{loc}$ converging subsequence such that the limit
$(\overline{M}^4,\overline{g}_{ij}(t))$ is an ancient solution to
the Ricci flow and satisfies the restricted isotropic curvature
pinching condition (2.4), as well as is $\kappa$-noncollapsed for
all scales. Moreover, the limit still satisfies the
Li-Yau-Hamilton inequality and the assertions (i) and (ii) of
Proposition 3.6. To show the limit is an ancient
$\kappa$-solution, we remain to show the limit has bounded
curvature at the time $t=0$.

By the virtue of Lemma 3.2, we may assume the limit has positive
curvature operator everywhere. We now argue by contradiction.
Suppose the curvature of the limit (at $t=0$)
$(\bar{M}^{4},\bar{g}_{ij}(0))$ is unbounded, then there is a
sequence of points $P_{l}$ divergent to infinity at the time $t=0$
with the scalar curvature $R(P_{l},0)\rightarrow +\infty$. By
Hamilton's compactness theorem (Theorem 16.1 of \cite{Ha6}) and the
estimates in the assertions (i) and (ii) of Proposition 3.6, we know
that a subsequence of the rescaled solutions
$(\overline{M}^4,\bar{R}(P_{l},0)\overline{g}_{ij}(\cdot,\frac{t}{\bar{R}(P_{l},0)}),P_{l})$
converges in $C^{\infty}_{loc}$ to a smooth nonflat limit. And by
Lemma 3.1, the limit must be the round neck $\mathbb{R}\times
\mathbb{S}^{3}$. This contradicts Proposition 2.2.

Therefore we have proved the corollary.
$$\eqno \#$$ \vskip 0.5cm

\subsection*{3.4  Canonical neighborhood structures}

We now examine the structures of four-dimensional nonflat ancient
$\kappa$-solutions with restricted isotropic curvature pinching. As
before by Lemma 3.2, we have seen a four-dimensional nonflat ancient
$\kappa$-solution with restricted isotropic curvature pinching,
whose curvature operator has a nontrivial null vector somewhere at
some time, must be a metric quotient of the round cylinder
$\mathbb{R}\times \mathbb{S}^3$. So we only need to consider the
ancient $\kappa$-solutions with \textbf{\emph{positive}} curvature
operator. The following theorem gives their canonical neighborhood
structures. The analogous result in three-dimensional case was given
by Perelman in Section 1.5 of \cite{P2}.\vskip 0.3cm

$\underline{\mbox{\textbf{Theorem 3.8}}}$ \emph{ For every
$\varepsilon>0$ one can find positive constants
$C_1=C_1(\varepsilon)$, $C_2=C_2(\varepsilon)$ such that for each
point $(x,t)$ in every four-dimensional ancient $\kappa$-solution
(for some $\kappa>0$) with restricted isotropic curvature pinching
and with \textbf{positive} curvature operator, there is a radius
$r$, $0<r<C_1(R(x,t))^{-\frac{1}{2}}$, so that some open
neighborhood $B_t(x,r)\subset B\subset B_t(x,2r)$ falls into one of
the following three categories:}

\emph{(a) $B$ is an \textbf{evolving $\varepsilon$-neck} (in the
sense that it is the time slice at time $t$ of the parabolic
region $\{(x',t')|x'\in B, t'\in [t-\varepsilon^{-2}
R(x,t)^{-1},t]\}$ which
 is, after
scaling with factor $R(x,t)$ and shifting the time $t$ to $0$,
$\varepsilon$-close (in $C^{[\varepsilon^{-1}]}$ topology) to the
 subset $(\mathbb{I}\times \mathbb{S}^{3})\times [-\varepsilon^{-2},0]$
 of the evolving round cylinder
$\mathbb{R}\times \mathbb{S}^3$, having scalar curvature one and
length $2\varepsilon^{-1}$ to $\mathbb{I}$ at time zero, or}

\emph{(b) $B$ is an $\textbf{evolving $\varepsilon$-cap}$ (in the
sense that it is the time slice at the time $t$ of an evolving
metric on open $\mathbb{B}^4$ or $\mathbb{RP}^4\setminus
\overline{\mathbb{B}^4}$ such that the region outside some
suitable compact subset of $\mathbb{B}^4$ or
$\mathbb{RP}^4\setminus \overline{\mathbb{B}^4}$ is an evolving
$\varepsilon$-neck), or}

\emph{(c) $B$ is a compact manifold (without boundary) with
positive curvature operator
(thus it is diffeomorphic to $\mathbb{S}^4$ or $\mathbb{RP}^4$); \\
furthermore, the scalar curvature of the ancient $\kappa$-solution
in $B$ at time $t$ is between $C^{-1}_2R(x,t)$ and $C_2R(x,t)$, and
the volume of $B$ in case (a) and case (b) satisfies
$$
(C_2R(x,t))^{-2} \leq Vol_t(B).
$$. }\vskip 0.2cm

$\underline{\mbox{\textbf{Proof}}}$. If the nonflat ancient
$\kappa$-solution is noncompact, the conclusions follow
immediately from (the proof of) Proposition 3.4. We thus assume
the nonflat ancient $\kappa$-solution is compact. By Theorem 3.5
we see that such ancient $\kappa$-solution is
$\kappa_0$-noncollapsed for all scales for some universal positive
constant $\kappa_0$.

We argue by contradiction. Suppose for some $\varepsilon>0$, there
exists a sequence of compact ancient $\kappa_0$-solutions
($M^4_k,g_k$) with restricted isotropic curvature pinching and
with positive curvature operator, a sequence of points $x_k\in
M^4_k$, and sequences of positive constants $C_{1k}$ with
$C_{1k}\rightarrow+\infty$ as $k\rightarrow+\infty$ and
$C_{2k}=\omega(4C^2_{1k})$ with the function $\omega$ given in
Proposition 3.6 such that at time $t$, for every radius $r,\
0<r<C_{1k}R(x_k,t)^{-\frac{1}{2}}$, any open neighborhood $B$ with
$B_t(x_k,r)\subset B\subset B_t(x_k,2r)$ can not fall into any one
of the three categories (a), (b) and (c). Clearly, the diameter of
each $M^4_k$ at time $t$ is at least
$C_{1k}R(x_k,t)^{-\frac{1}{2}}$; otherwise one can choose suitable
$r\in (0,C_{1k}R(x_k,t)^{-\frac{1}{2}})$ and $B=M^4_k$, which
falls into the category (c), so that the scalar curvature in $B$
at $t$ is between $C^{-1}_{2k}R(x_k,t)$ and $C_{2k}R(x_k,t)$ by
using Proposition 3.6 (i). Now by scaling the ancient
$\kappa_0$-solutions along the points $(x_k,t)$ with factors
$R(x_k,t)$ and shifting the time $t$ to $0$, it follows from
Corollary 3.7  that a subsequence of these rescaled ancient
$\kappa_0$-solutions converge in $C^{\infty}_{loc}$ topology to a
noncompact nonflat ancient $\kappa_0$-solution with restricted
isotropic curvature pinching.

If the noncompact limit has a nontrivial null curvature
eigenvector somewhere, then by Lemma 3.2 we conclude that the
limit is round cylinder $\mathbb{R}\times \mathbb{S}^{3}$ or a
metric quotient $\mathbb{R}\times \mathbb{S}^3/\Gamma$. By the
same reason as in the proof of Theorem 3.5, the projection
$\Gamma_{1}$ of $\Gamma$ to the factor $\mathbb{R}$ is an
isometric subgroup of $\mathbb{R}$. Since the limit
$\mathbb{R}\times \mathbb{S}^{3}/\Gamma$ is noncompact,
$\Gamma_{1}$ must be $\{1\}$ or $\mathbb{Z}_{2}$. Thus we have a
$\Gamma$-invariant cross-sphere $\mathbb{S}^{3}$ in
$\mathbb{R}\times \mathbb{S}^{3}/\Gamma$, and $\Gamma$ acts on it
without fixed points. Denote this cross sphere by $\{0\}\times
\mathbb{S}^{3}$. Since each $(M^4_k,g_k)$ is compact and has
positive curvature operator, we know from \cite{Ha2} that each
$M^4_k$ is diffeomorphic to $\mathbb{S}^4$ or $\mathbb{RP}^4$.
Then by the proof of Theorem 3.5 and applying theorem C 4.1 of
\cite{Ha7}, we conclude that the limit is $\mathbb{R}\times
\mathbb{S}^{3}$ or $\mathbb{R}\times \mathbb{S}^{3}/\Gamma $ with
$\Gamma=\mathbb{Z}_{2}.$ If $\Gamma=\mathbb{Z}_{2}$, we claim that
$\Gamma$ must act on $\mathbb{R}\times \mathbb{S}^{3}/\Gamma $ by
flipping both $\mathbb{R}$ and $\mathbb{S}^{3}$.

Indeed, as shown before, $\Gamma_{1}=\{1\}$ or $\mathbb{Z}_{2}$.
If $\Gamma_{1}=\{1\}$, then $\mathbb{R}\times
\mathbb{S}^{3}/\Gamma=\mathbb{R}\times \mathbb{RP}^{3}$. Let
$\Gamma^{+}$ be the normal subgroup of $\Gamma$ preserving the
orientation of the cylinder, and $\pi_{1}(M_{k}^{4})^{+}$
 be the normal subgroup of $\pi_{1}(M_{k}^{4})$ preserving the orientation
 of the universal cover of $M_{k}^{4}$. Since the manifold $M^{4}_{k}$ is
diffeomorphic to $\mathbb{R}\mathbb{P}^{4}$, this induces an
absurd commutative diagram:
$$
\begin{array}{ccccccccc}
 & & \mathbb{Z}_2 & & \mathbb{Z}_2 & & 0 & &  \\
 & & \| & & \| & &\| & &\\
 0&\longrightarrow &\Gamma^+& \longrightarrow & \Gamma
 &\longrightarrow & \Gamma/\Gamma^{+}& \longrightarrow& 0\\
 & & \downarrow & & \downarrow \wr& & \downarrow & &\\
 0&\longrightarrow &\pi_{1}(M_{k}^4)^{+}& \longrightarrow & \pi_{1}(M_{k}^4)
 &\longrightarrow & \pi_{1}(M_{k}^4)/\pi_{1}(M_{k}^4)^{+}& \longrightarrow & 0\\
 & & \|& & \|& &\| & &\\
& & 0& & \mathbb{Z}_{2}& &\mathbb{Z}_{2} & &
 \end{array}
$$
where the vertical morphisms are induced by the inclusion
$\mathbb{R}\times \mathbb{S}^{3}/\Gamma\subset M $. Therefore
$\Gamma_{1}=\mathbb{Z}_{2}$. Denote by $\sigma_{1}$ be the
isometry of $\mathbb{R}\times \mathbb{S}^{3}$ acting by flipping
both $\mathbb{R}$ and $\mathbb{S}^{3}$ around $\{0\}\times S^{3}$.
Clearly, for any $\sigma\in \Gamma$ with $\sigma\neq 1$,
$\sigma\circ\sigma_{1}$ is an isometry of $\mathbb{R}\times
\mathbb{S}^{3}$ whose projection on the factor $\mathbb{R}$ is the
identity map. Then $\sigma\circ\sigma_{1}$ is only a rotation of
the factor $\mathbb{S}^{3}$ in $\mathbb{R}\times \mathbb{S}^{3}$.
Note that $\sigma\circ\sigma_{1}\mid_{\{0\}\times \mathbb{S}^{3}}$
is identity. We conclude that $\sigma=\sigma_{1}$ and the claim
holds.

 When the
limit is the round cylinder $\mathbb{R}\times \mathbb{S}^3$, a
suitable neighborhood $B$ (for suitable $r$) of $x_k$ would fall
into the category (a) for sufficiently large $k$; while when the
limit is the $\mathbb{Z}_2$ quotient of the round cylinder
$\mathbb{R}\times \mathbb{S}^3$ with the antipodal map flipping
both $\mathbb{S}^3$ and $\mathbb{R}$, a suitable neighborhood $B$
(for suitable $r$) of $x_k$ would fall into the category (b) (over
$\mathbb{R}\mathbb{P}^4\setminus \overline{\mathbb{B}^4}$) or into
the category (a) for sufficiently large $k$. This is a
contradiction.

If the noncompact limit has positive curvature operator
everywhere, then by Proposition 3.4, a suitable neighborhood $B$
(for suitable $r$) of $x_k$ would fall into the category (b) (over
$\mathbb{B}^4$) for sufficiently large $k$. We also get a
contradiction.

Finally, the statements on the curvature estimate and volume
estimate for the neighborhood $B$ follows directly from Theorem 3.6
and Proposition 3.4. Therefore we have proved the theorem.
$$\eqno \#$$ \vskip 1cm

 \centerline{\large{\textbf{4. The Structure of Solutions at the
Singular Time}}} \vskip 0.5cm
 Let
$(M^4,g_{ij}(x))$ be a four-dimensional compact Riemannian manifold
with positive isotropic curvature and let $g_{ij}(x,t)$, $x\in M^4$
and $t\in [0,T),$ be a maximal solution to the Ricci flow (1.1) with
$g_{ij}(x,0)=g_{ij}(x)$ on $M^4$. Since the initial metric
$g_{ij}(x)$ has positive scalar curvature, it is easy to see that
the maximal time $T$ must be finite and the curvature tensor becomes
unbounded as $t\rightarrow T$. According to Perelman's noncollapsing
theorem I (Theorem 4.1 of \cite{P1}), the solution $g_{ij}(x,t)$ is
$\kappa$-noncollapsed on the scale $\sqrt{T}$ for all $t\in [0,T)$
for some $\kappa>0$. Now let us take a sequence of times
$t_k\rightarrow T$, and a sequence of points $p_k\in M^4$ such that
for some positive constant $C$, $|R_m|(x,t)\leq CQ_k$ with
$Q_{k}=|Rm(p_{k},t_{k})|$ whenever $x\in M^4$ and $t\in [0,t_k]$,
called a sequence of (almost) \textbf{maximal points}. Then by
Hamilton's compactness theorem \cite{Ha5}, a sequence of the
scalings of the solution $g_{ij}(x,t)$ along the points $p_k$ with
factors $Q_k$ converges to a complete ancient $\kappa$-solution with
restricted isotropic curvature pinching. This says, for any
$\varepsilon>0$, there exists a positive number $k_0$ such that as
$k\geq k_0$, the solution in the parabolic region $\{(x,t)\in
M^4\times [0,T)\ |\ d^2_{t_k}(x,x_k)<\varepsilon^{-2}
Q_k^{-1},t_k-\varepsilon^{-2} Q_k^{-1}<t\leq t_k\}$ is, after
scaling with the factor $Q_k$, $\varepsilon$-close (in
$C^{[\varepsilon^{-1}]}$-topology) to the corresponding subset of
the ancient $\kappa$-solution with restricted isotropic curvature
pinching.

Let us describe the structure of any ancient $\kappa$-solution
(with restricted isotropic curvature pinching).
 If the curvature operator is positive
everywhere, then each point of the ancient $\kappa$-solution has a
canonical neighborhood described in Theorem 3.8. While if the
curvature operator has a nontrivial null eigenvector somewhere,
then by Hamilton's strong maximum principle and the pinching
condition (2.4) the ancient $\kappa$-solution is isometric to
$\mathbb{R}\times \mathbb{S}^3/\Gamma$, a metric quotient of the
round cylinder $\mathbb{R}\times \mathbb{S}^3$. Since it is
$\kappa$-noncollapsed for all scales, the metric quotient
$\mathbb{R}\times \mathbb{S}^3/\Gamma$ can not be compact. Suppose
we make an additional assumption that the compact four manifold
$M^4$ has no essential incompressible space form. Then by the
proofs of Theorems 3.5 and 3.8 and applying Theorem C 4.1 of
\cite{Ha7}, we have $\Gamma=\{1\}$, or $\Gamma=\mathbb{Z}_2$
acting antipodally on $\mathbb{S}^3$ and by reflection on
$\mathbb{R}$. Thus in both cases, each point of the ancient
$\kappa$ solution has also a canonical neighborhood described in
Theorem 3.8.

Hence we see that each such (almost) maximal point $(x_k,t_k)$ has a
canonical neighborhood which is either an evolving
$\varepsilon$-neck or an evolving $\varepsilon$-cap, or a compact
manifold (without boundary) with positive curvature operator. This
gives the structure of the singularities coming from a sequence of
(almost) maximal points $(x_k,t_k)$. However this argument does not
work for the singularities coming from a sequence of points
$(y_k,\tau_k)$ with $\tau_k\rightarrow T$ and
$|Rm(y_k,\tau_k)|\rightarrow+\infty$ when $|Rm(y_k,\tau_k)|$ is not
comparable with the maximum of the curvature at the time $\tau_k$,
since we can not take a limit directly. We now follow a refined
rescaling argument of Perelman (Theorem 12.1 of \cite{P1}) to obtain
a uniform canonical neighborhood structure theorem for
four-dimensional solutions at any point where its curvature is
suitable large. \vskip 0.3cm

$\underline{\mbox{\textbf{Theorem 4.1}}}$ \emph{Given
$\varepsilon>0$, $\kappa >0$, $0<\theta,\rho,\Lambda, P<+\infty$,
one can find $r_0>0$ with the following property. If $g_{ij}(x,t)$,
$t\in [0,T)$ with $T>1$, is a solution to the Ricci flow on a
four-dimensional manifold $M^{4}$ with no essential incompressible
space form, which has positive isotropic curvature, is
$\kappa$-noncollapsed on the scales $\leq \theta$ and satisfies
(2.1), (2.2) and (2.3) in Lemma 2.1, then for any point $(x_0,t_0)$
with $t_0\geq 1$ and $Q=R(x_0,t_0)\geq r^{-2}_0$, the solution in
the parabolic region $\{(x,t)\in M^4\times [0,T)|
d^2_{t_0}(x,x_0)<\varepsilon^{-2} Q^{-1},t_0-\varepsilon^{-2}
Q^{-1}<t\leq t_0\}$ is, after scaling by the factor $Q$,
$\varepsilon$-close (in $C^{[\varepsilon^{-1}]}$-topology) to the
corresponding subset of some ancient $\kappa$-solution with
restricted isotropic curvature pinching.}

 \emph{Consequently each point
$(x_0,t_0)$, with $t_0\geq 1$ and $Q=R(x_0,t_0)\geq r^{-2}_0$,
satisfies the gradient estimates
$$|\nabla R(x_0,t_0)|<2\eta R^{\frac{3}{2}}(x_0,t_0)
\mbox{ and }|\frac{\partial}{\partial t}R(x_0,t_0)|<2\eta
R^2(x_0,t_0), \eqno (4.1)$$ and has a canonical neighborhood $B$
with $B_{t_0}(x_0,r)\subset B \subset B_{t_0}(x_0,2r)$ for some $
0<r<C_1(\varepsilon)(R(x_0,t_0))^{-\frac{1}{2}}$, which is either an
evolving $\varepsilon$-neck, or an evolving $\varepsilon$-cap, or a
compact four-manifold with positive curvature operator. Here $\eta$
is the universal constant in Proposition 3.6 and $C_1(\varepsilon)$
is the positive constant in Theorem 3.8.}\vskip 0.2cm

 $\underline{\mbox{\textbf{Proof}}}$. Let
$C(\varepsilon)$ be a positive constant depending only on
$\varepsilon$ such that $C(\varepsilon)\rightarrow+\infty$ as
$\varepsilon\rightarrow 0^+$. It suffices to prove that there exists
$r_0>0$ such that for any point $(x_0,t_0)$ with $t_0\geq 1$ and
$Q=R(x_0,t_0)\geq r^{-2}_0$, the solution in the parabolic region
$\{(x,t)\in M^4\times [0,T)\ |\  d^2_{t_0}(x,x_0)<C(\varepsilon)
Q^{-1},t_0-C(\varepsilon) Q^{-1}<t\leq t_0\}$ is, after scaling by
the factor $Q$, $\varepsilon$-close to the corresponding subset of
some ancient $\kappa$-solution with restricted isotropic curvature
pinching. The constant $C(\varepsilon)$ will be determined later.

We argue by contradiction. Suppose for some $\varepsilon>0$,
$\kappa
>0$, $0<\theta, \rho, \Lambda,P<+\infty$, there exists a sequence of
solutions $(M^4_k,g^{(k)}_{ij}(\cdot,t))$ to the Ricci flow on
compact four-manifolds with no essential incompressible space form,
having positive isotropic curvature and satisfying (2.1), (2.2) and
(2.3), defined on the time interval $[0,T_k)$ with $T_k>1$, and a
sequence of positive numbers $r_k\rightarrow 0$ such that each
solution $(M^4_k,g^{(k)}_{ij}(\cdot,t))$ is $\kappa$-noncollapsed on
the scales $\leq \theta$; but there exists a sequence of points
$x_k\in M^4_k$ and times $t_k\geq 1$ with $Q_k=R_k(x_k,t_k)\geq
r^{-2}_k$ such that the solution in the parabolic region $\{(x,t)\in
M^4_k\times [0,T_k)\ |\  d^2_{t_k}(x,x_k)<C(\varepsilon)
Q^{-1}_k,t_k-C(\varepsilon) Q^{-1}_k<t\leq t_k\}$ is not, after
scaling by the factor $Q_k$, $\varepsilon$-close to the
corresponding subset of any ancient $\kappa$-solution with
restricted isotropic curvature pinching, where $R_k$ denotes the
scalar curvature of $(M^4_k,g^{(k)}_{ij}(\cdot,t))$. For each
solution $(M^4_k,g^{(k)}_{ij}(\cdot,t))$, we may adjust the point
$(x_k,t_k)$ with $t_k\geq \frac{1}{2}$ and with $Q_k=R_k(x_k,t_k)$
as large as possible so that the conclusion of the theorem fails at
$(x_k,t_k)$, but holds for any $(x,t)\in M^4_k\times
[t_k-H_kQ^{-1}_k,t_k]$ satisfying $R(x,t)\geq 2Q_k$, where
$H_k=\frac{1}{4}r^{-2}_k\rightarrow +\infty$ as
$k\rightarrow+\infty$. Indeed, suppose not, by setting
$(x^{(1)}_k,t^{(1)}_k)=(x_k,t_k)$, we can inductively choose
$(x^{(\ell)}_k,t^{(\ell)}_k)\in M^4_k\times
[t^{(\ell-1)}_k-H_k(R_k(x^{(\ell-1)}_k,t^{(\ell-1)}_k))^{-1},t^{(\ell-1)}_k]$
satisfying $R_k(x^{(\ell)}_k,t^{(\ell)}_k)\geq
2R_k(x^{(\ell-1)}_k,t^{(\ell-1)}_k)$, but the conclusion of the
theorem fails at $(x^{(\ell)}_k,t^{(\ell)}_k)$ for each
$\ell=2,3,\cdots$. Since the solution is smooth and
\begin{eqnarray*}R_k(x^{(\ell)}_k,t^{(\ell)}_k)&\geq&
2R_k(x^{(\ell-1)}_k,t^{(\ell-1)}_k)\\ &\geq&
2^{\ell-1}R_k(x_k,t_k),
\end{eqnarray*}
\begin{eqnarray*}
t^{(\ell)}_k&\geq&
t^{(\ell-1)}_k-H_k(R_k(x^{(\ell-1)}_k,t^{(\ell-1)}_k))^{-1}\\
&\geq& t_k-H_k\sum^{\ell-1}_{i=1}(2^{i-1}R_k(x_k,t_k))^{-1}\\&\geq
&t_k-2H_k(R(x_k,t_k))^{-1}\\ &\geq& \frac{1}{2} ,\end{eqnarray*}
the above choosing process must terminate in finite step and the
last element fits.

Let $(M^4_k,\widetilde{g}^{(k)}_{ij}(\cdot,t),x_k)$ be the rescaled
solutions obtained by rescaling the manifolds
$(M^4_k,g^{(k)}_{ij}(\cdot,t))$ with factors $Q_k=R_k(x_k,t_k)$ and
shifting the time $t_{k}$ to $0$. Denote by $\widetilde{R}_k$ the
rescaled scalar curvature. We will show that a subsequence of the
rescaled solutions $(M^4_k,\widetilde{g}^{(k)}_{ij}(\cdot,t),x_k)$
converges to an ancient $\kappa$-solution with restricted isotropic
curvature pinching, which is a contradiction. In the followings we
divide the argument into four steps. \vskip 0.2cm

$\underline{\mbox{\textbf{Step 1}}}$ We want to prove a local
curvature estimate in the following assertion.

$\underline{\mbox{\textbf{Claim}}}$: \emph{For each
$(\overline{x},\overline{t})$ with
$t_k-\frac{H_k}{2}Q^{-1}_k<\overline{t}\leq t_k$, we have
$R_k(x,t)\leq 4\overline{Q}_k$ whenever
$\overline{t}-c\overline{Q}^{-1}_k\leq t\leq \overline{t}$ and
$d^2_{\overline{t}}(x,\overline{x})\leq c\overline{Q}^{-1}_k$,
where $\overline{Q}_k=Q_k+R_k(\overline{x},\overline{t})$ and
$c>0$ is a small universal constant.}

To prove this, we consider any point $(x,t)\in
B_{\overline{t}}(\overline{x},(c\overline{Q}^{-1}_k)^{\frac{1}{2}})\times
[\overline{t}-c\overline{Q}^{-1}_k,\overline{t}]$ with $c>0$ to be
determined. If $R_k(x,t)\leq 2Q_k$, there is nothing to show. If
$R_k(x,t)>2Q_k$, consider a space time curve $\gamma$ that goes
straightly from $(x,t)$ to $(x,\overline{t})$ and goes from
$(x,\overline{t})$ to $(\overline{x},\overline{t})$ along a
minimizing geodesic (with respect to the metric $g^{(k)}_{ij}(\cdot,
\overline{t})$). If there is a point on $\gamma$ with the scalar
curvature $2Q_k$, let $p$ be the nearest such point to $(x,t)$; if
not, put $p=(\overline{x},\overline{t})$. On the segment of $\gamma$
from $(x,t)$ to $p$, the scalar curvature is not less than $2Q_k$.
According to the choice of the point $(x_k,t_k)$, the solution along
the segment is $\varepsilon$-close to that of some ancient
$\kappa$-solutions with restricted isotropic curvature pinching. Of
course we may assume $\varepsilon>0$ is very small. It follows from
Proposition 3.6 (ii) that $$|\nabla (R^{-\frac{1}{2}}_k)|\leq 2\eta
\mbox{ and } |\frac{\partial}{\partial t}(R^{-1}_k)|\leq 2\eta$$ on
the segment for some universal constant $\eta>0$. Then by choosing
$c>0$ (depending only on $\eta$) small enough we get the desired
curvature bound by integrating the above derivative estimate along
the segment. This proves the assertion. \vskip 0.2cm

$\underline{\mbox{\textbf{Step 2}}}$ \  We next want to show that
the curvatures of the rescaled solutions
$\widetilde{g}^{(k)}_{ij}(\cdot,t)$ at the new time $t=0$ (i.e.,
the original time $t_k$) stay uniformly bounded at bounded
distances from $x_k$.

For all $\sigma\geq 0$, set
$$M(\sigma)=\sup\{\widetilde{R}_k(x,0)\ |\ k\geq 1,x\in M^4_k \mbox{ with }d_0(x,x_k)\leq
\sigma\}$$ and $$\sigma_0=\sup\{\sigma\geq 0\ |\
M(\sigma)<+\infty\}.$$ Note that $\sigma_0>0$ by Step 1. By the
assumptions (2.1) and (2.2), it suffices to show $\sigma_0=+\infty$.
We now argue by contradiction to show $\sigma_0=+\infty$. Suppose
not, we may find (after passing to a subsequence if necessary) a
sequence of points $y_k\in M^4_k$ with $d_0(y_k,x_k)\rightarrow
\sigma_0<+\infty$ and $\widetilde{R}_k(y_k,0)\rightarrow+\infty$ as
$k\rightarrow+\infty$. Let $\gamma_k(\subset M^4_k)$ be a minimizing
geodesic segment from $x_k$ to $y_k$, $z_k\in \gamma_k$ the point on
$\gamma_k$ closest to $y_k$ at which $\widetilde{R}_k(z_k,0)=2$ and
$\beta_k$ the subsegment of $\gamma_k$ running from $y_k$ to $z_k$.
By Step 1 the length of $\beta_k$ is uniformly bounded away from
zero for all $k$. And by the assumptions (2.1) and (2.2), we have a
uniform curvature bound on the open balls
$B_0(x_k,\sigma)\subset(M^4_k,\widetilde{g}^{(k)}_{ij}(\cdot,0))$
for each fixed $\sigma<\sigma_0$. Note that the
$\kappa$-noncollapsing assumption implies the uniform injectivity
radius bound for $(M^4_k,\widetilde{g}^{(k)}_{ij}(\cdot,0))$ at the
marked points $x_k$. Then by the virtue of Hamilton's compactness
theorem 16.1 in \cite{Ha6} ( see \cite{CaZ} for the details on
generalizing Hamilton's compactness theorem to finite balls), we can
extract a subsequence of the marked
($B_0(x_k,\sigma_0),\widetilde{g}^{(k)}_{ij}(\cdot,0),x_k)$ which
converges in $C^{\infty}_{loc}$ topology to a marked (noncomplete)
manifold ($B_{\infty},\widetilde{g}^{\infty}_{ij},x_{\infty}$), so
that the segments $\gamma_k$ converge to a geodesic segment (missing
an endpoint) $\gamma_{\infty}\subset B_{\infty}$ emanating from
$x_{\infty}$, and $\beta_k$ converge to a subsegment
$\beta_{\infty}$ of $\gamma_{\infty}$. Let $\overline{B}_{\infty}$
denote the completion of
$(B_{\infty},\widetilde{g}^{(\infty)}_{ij})$, and $y_{\infty}\in
\overline{B}_{\infty}$ the limit point of $\gamma_{\infty}$.

Denote by $\widetilde{R}_{\infty}$ the scalar curvature of
$(B_{\infty},\widetilde{g}^{(\infty)}_{ij})$. Since the rescaled
scalar curvatures of $\widetilde{R}_k$ along $\beta_k$ are at
least 2, it follows from the choice of the points $(x_k,t_k)$ that
for any $q_0\in \beta_{\infty}$, the manifold
$(B_{\infty},\widetilde{g}^{(\infty)}_{ij})$ in $\{q\in
B_{\infty}|dist^2_{\widetilde{g}^{(\infty)}_{ij}}(q,q_0)<C(\varepsilon)(\widetilde{R}_{\infty}(q_0))^{-1}\}$
is $2\varepsilon$-close to the corresponding subset of (a time
slice of) of some ancient $\kappa$-solution with restricted
isotropic curvature pinching. From the argument in the second
paragraph of this section, we know that such an ancient
$\kappa$-solution with restricted isotropic curvature pinching at
each point $(x,t)$ has a radius $r$,
$0<r<C_1(2\varepsilon)R(x,t)^{-\frac{1}{2}}$, such that its
canonical neighborhood $B$ with $ B_t(x,r)\subset B\subset
B_{t}(x,2r)$,  is either an evolving $2\varepsilon$-neck, or an
evolving $2\varepsilon$-cap, or a compact manifold (without
boundary) with positive curvature operator, moreover the scalar
curvature on the ball is between $C_2(2\varepsilon)^{-1}R(x,t)$
and $C_2(2\varepsilon)R(x,t)$, where $C_1(2\varepsilon)$ and
$C_2(2\varepsilon)$ are the positive constants in Theorem 3.8. We
now choose
$C(\varepsilon)=\max\{2C_1(2\varepsilon)^{2},\varepsilon^{-2}\}$.
By the local curvature estimate in Step 1, we see that the scalar
curvature $\widetilde{R}_{\infty}$ becomes unbounded along
$\gamma_{\infty}$ going to $y_{\infty}$. This implies that the
canonical neighborhood around $q_0$ can not be a compact manifold
(without boundary) with positive curvature operator. Note that
$\gamma_{\infty}$ is shortest since it is the limit of a sequence
of shortest geodesics. Without loss of generality, we may assume
$\varepsilon$ is suitably small. These imply that as $q_0$
sufficiently close to $y_{\infty}$, the canonical neighborhood
around $q_0$ can not be a $2\varepsilon$-cap.
 Thus we conclude that each $q_0\in
\gamma_{\infty}$, which is sufficiently close to $y_{\infty}$, is
the center of a $2\varepsilon$-neck.

Denote by $$U=\bigcup_{q_0\in
\gamma_{\infty}}B(q_0,24\pi(\widetilde{R}_{\infty}(q_0))^{-\frac{1}{2}})\
\ \ (\subset(B_{\infty},\widetilde{g}^{(\infty)}_{ij}))$$ where
$B(q_0,24\pi(\widetilde{R}_{\infty}(q_0))^{-\frac{1}{2}})$ is the
ball centered at $q_0$ of radius
$24\pi(\widetilde{R}_{\infty}(q_0))^{-\frac{1}{2}}$.

Clearly, it follows from the assumptions (2.1), (2.2) and (2.3)
that $U$ has nonnegative curvature operator. Since the metric
$\widetilde{g}^{(\infty)}_{ij}$ is cylindrical at any point
$q_0\in \gamma_{\infty}$ which is sufficiently close to
$y_{\infty}$, we see that the metric space $\overline{U}=U\cup
\{y_{\infty}\}$ by adding the point $y_{\infty}$, is locally
complete and strictly intrinsic near $y_{\infty}$. Here strictly
intrinsic means that the distance between any two points can be
realized by shortest geodesics. Furthermore $y_{\infty}$ cannot be
an interior point of any geodesic segment in $\overline{U}$. This
implies that the curvature of $\overline{U}$ at $y_{\infty}$ is
nonnegative in Alexandrov sense. Note that for any very small
radius $\sigma$, the geodesic sphere $\partial
B(y_{\infty},\sigma)$ is an almost round sphere of radius $\leq 3
\varepsilon\sigma \pi$. By \cite{BGP} or \cite{CC} we have a
four-dimensional tangent cone at $y_{\infty}$ with aperture $\leq
20\varepsilon$. Moreover, by \cite{BGP} or \cite{CC}, any
four-dimensional tangent cone $C_{y_{\infty}}\overline{U}$ at
$y_{\infty}$ must be a metric cone. For each tangent cone, pick
$z\in C_{y_{\infty}}\overline{U}$ such that the distance between
the vertex $y_{\infty}$ and $z$ is one. Then the ball
$B(z,\frac{1}{2})\subset C_{y_{\infty}}\overline{U}$ is the
Gromov-Hausdorff limit of the scalings of a sequence of balls
$B_0(z_k,s_k)\subset(M^4_k,\widetilde{g}^{(k)}_{ij}(\cdot,0))$ by
some factors $a_k$, where $s_k\rightarrow 0^+$. Since the tangent
cone is four-dimensional and has aperture $\leq 20\varepsilon$,
the factors $a_k$ must be comparable with
$\tilde{R}_{k}(z_{k},0)$. By using the local curvature estimate in
Step 1, we actually have the convergence in the $C^{\infty}_{loc}$
topology for the solutions $\widetilde{g}^{(k)}_{ij}(\cdot,t)$
over the balls $B_0(z_k,s_k)$ and over some time interval $t\in
[-\delta, 0]$ for some sufficiently small $\delta>0$. The limiting
$B(z,\frac{1}{2})\subset C_{y_{\infty}}\overline{U}$ is a piece of
the nonnegative (operator) curved and nonflat metric cone. On the
other hand, since the radial direction of the cone is flat, by
Hamilton's strong maximum principle \cite{Ha2} and the pinching
condition (2.4) as in the proof of Lemma 3.2, the limiting
$B(z,\frac{1}{2})$ would be a piece of $\mathbb{R}\times
\mathbb{S}^3$ or $\mathbb{R}\times \mathbb{S}^3/\Gamma$ (a metric
quotient). This is a contradiction. So we have proved that the
curvatures of the rescaled metrics
$\widetilde{g}^{(k)}_{ij}(\cdot,0)$ stay uniformly bounded at
bounded distances from $x_k$.

By the local curvature estimate in Step 1, we can locally extend the
above curvature control backward in time a little. Then by the
$\kappa$-noncollapsing assumption and Shi's derivative estimates
\cite{Sh1}, we can take a $C^{\infty}_{loc}$ limit from the sequence
of marked rescaled solutions
$(M^4_k,\widetilde{g}^{(k)}_{ij}(\cdot,t),x_k)$. The limit, denoted
by
($M^4_{\infty},\widetilde{g}^{(\infty)}_{ij}(\cdot,t),x_{\infty}$),
is $\kappa$-noncollapsing on all scales, is defined on a space-time
open subset of $M^4_{\infty}\times (-\infty,0]$ containing the time
slice $M^4_{\infty}\times \{0\}$, and satisfies the restricted
isotropic curvature pinching condition (2.4) by the assumptions
(2.1), (2.2) and (2.3). \vskip 0.2cm

$\underline{\mbox{\textbf{Step 3}}}$  We further claim that the
limit ($M^4_{\infty},\widetilde{g}^{(\infty)}_{ij}(\cdot,t)$) at
the time slice $t=0$ has bounded curvature.

We have known that the curvature operator of the limit
($M^4_{\infty},\widetilde{g}^{(\infty)}_{ij}(\cdot,t)$) is
nonnegative everywhere. If the curvature operator has a nontrivial
null eigenvector somewhere, we can argue as in the proof of Lemma
3.2 by using Hamilton's strong maximum principle \cite{Ha2} and
the restricted isotropic curvature pinching condition (2.4) to
deduce that the universal cover of the limit is isometric to the
standard $\mathbb{R}\times \mathbb{S}^3$. Thus the curvature of
the limit is bounded in this case.

Assume that the curvature operator of the limit
($M^4_{\infty},\widetilde{g}^{(\infty)}_{ij}(\cdot,t)$) at the time
slice $t=0$ is positive everywhere. Suppose there exists a sequence
of points $p_j\in M^4_{\infty}$ such that their scalar curvatures
$\widetilde{R}_{\infty}(p_j,0)\rightarrow +\infty$ as $j\rightarrow
+\infty$. By the local curvature estimate in Step 1 and the
assertion of the above Step 2 (for the marked points $p_{j}$) as
well as the $\kappa$-noncollapsed assumption, a subsequence of the
rescaled and marked manifolds
$(M^4_{\infty},\widetilde{R}_{\infty}(p_j,0)\widetilde{g}^{(\infty)}_{ij}(\cdot),p_j)$
converges in $C^{\infty}_{loc}$ topology to a smooth nonflat limit
$Y$. Then by Lemma 3.1 we conclude that $Y$ is isometric to
$\mathbb{R}\times \mathbb{S}^3$ with the standard metric. This
contradicts with Proposition 2.2. So the curvature of the limit
($M^4_{\infty},\widetilde{g}^{(\infty)}_{ij}(\cdot,t)$) at the time
slice $t=0$ must be bounded. \vskip 0.2cm

$\underline{\mbox{\textbf{Step 4}}}$  Finally we want to extend
the limit backward in time to $-\infty$.

By the local curvature estimate in Step 1, we now know that the
limit ($M^4_{\infty},\widetilde{g}^{(\infty)}_{ij}(\cdot,t)$) is
defined on $[-a,0]$ for some $a>0$.

 Denote by
\begin{eqnarray*} t' = \inf \{\ \ \tilde{t}&|&\mbox {we can take a
smooth limit  on } (\tilde{t},0] \mbox {(with bounded curvature at}
\\ & &\mbox {each time slice) from a subsequence of the rescaled
solutions } \tilde{g}_k \}.\end{eqnarray*} We first claim that there
is a subsequence of the rescaled solutions $\tilde{g}_k$ which
converges in $C^{\infty}_{loc}$ topology to a smooth limit
$(M_{\infty},\tilde{g}_{\infty}(\cdot,t))$ on the maximal time
interval $(t',0]$.

Indeed, let $t_k$ be a sequence of negative numbers such that
$t_k\rightarrow t'$ and there exist smooth limits
$(M_{\infty},\tilde{g}_{\infty}^{k}(\cdot,t))$ defined on
$(t_k,0]$. For each $k$, the limit has nonnegative and bounded
curvature operator at each time slice. Moreover by the Claim in
Step 1, the limit has bounded curvature on each subinterval
$[-b,0]\subset(t_k,0]$. Denote by $\tilde{Q}$ the scalar curvature
upper bound of the limit at time zero (where $\tilde{Q}$ is the
same for all $k$). Then we can apply
 Li-Yau-Hamilton inequality \cite{Ha4} to get
 $$
 \tilde{R}_{\infty}^{k}(x,t)\leq
\tilde{Q} (\frac{-t_k}{t-t_k}),
 $$
 where $\tilde{R}_{\infty}^{k}(x,t)$ are the scalar curvatures
 of the limits $(M_{\infty},\tilde{g}_{\infty}^{k}(\cdot,t))$.
 Hence by the definition of convergence and the above curvature estimates,
  we can find a subsequence of
the rescaled solutions $\tilde{g}_k$ which converges in
$C^{\infty}_{loc}$ topology to a smooth limit
$(M_{\infty},\tilde{g}_{\infty}(\cdot,t))$ on the maximal time
interval $(t',0]$.

 We next claim that $t'=-\infty$.

Suppose not, then the curvature of the limit
($M^4_{\infty},\widetilde{g}^{(\infty)}_{ij}(\cdot,t)$) becomes
unbounded as $t\rightarrow t'>-\infty$. Since the minimum of the
scalar curvature is nondecreasing in time and
$\widetilde{R}_{\infty}(x_{\infty},0)=1$, we see that there is a
$y_{\infty}\in M^4_{\infty}$ such that
$$0<\widetilde{R}_{\infty}(y_{\infty},t'+\frac{c}{3})<\frac{3}{2}$$
where $c>0$ is the universal constant in the assertion of Step 1.
By using Step 1 again we see that the limit
($M^4_{\infty},\widetilde{g}^{(\infty)}_{ij}(\cdot,t)$) in a small
neighborhood of the point $y_{\infty}$ at the time slice
$t=t'+\frac{c}{3}$ can be extended backward to the time interval
$[t'-\frac{c}{3},t'+\frac{c}{3}]$. We remark that the distances at
the time $t$ and the time $0$ are roughly equivalent in the
following sense $$d_t(x,y)\geq d_0(x,y)\geq d_t(x,y)-const. \eqno
(4.2)$$ for any $x,y\in M^4_{\infty}$ and $t\in (t',0]$. Indeed
from the Li-Yau-Hamilton inequality \cite{Ha4} we have the
estimate
$$\frac{\partial}{\partial t}\widetilde{R}_{\infty}(x,t)\geq
-\widetilde{R}_{\infty}(x,t)\cdot(t-t')^{-1}, \mbox{ for }(x,t)\in M^4_{\infty}\times
(t',0].$$ If $\widetilde{Q}$ denotes the supermum of the scalar
curvature $\widetilde{R}_{\infty}$ at $t=0$, then
$$\widetilde{R}_{\infty}(x,t)\leq \widetilde{Q}(\frac{-t'}{t-t'}), \mbox{ on }M^4_{\infty}\times
(t',0].$$ By applying Lemma 8.3 (b) of \cite{P1}, we have
$$d_t(x,y)\leq d_0(x,y)+30(-t')\sqrt{\widetilde{Q}}$$ for any
$x,y\in M^4_{\infty}$ and $t\in (t',0]$. On the other hand, since
the curvature operator of the limit
$\widetilde{g}^{\infty}_{ij}(\cdot,t)$ is nonnegative, we have
$$d_t(x,y)\geq d_0(x,y)$$ for any $x,y\in M^4_{\infty}$ and $t\in
(t',0]$. Thus we obtain the estimate (4.2).

The estimate (4.2) insures that the limit around the point
$y_{\infty}$ at any time $t\in (t',0]$ is exactly the original
limit around $x_{\infty}$ at the time $t=0$. Consider the rescaled
sequence of ($M^4_k,\widetilde{g}^{(k)}_{ij}(\cdot,t)$) with the
marked points replaced by the associated sequence $y_k\rightarrow
y_{\infty}$. By applying the same arguments as the above Step 2
and Step 3 to the new marked sequence
($M^4_k,\widetilde{g}^{(k)}_{ij}(\cdot,t),y_k$), we conclude the
original limit
($M^4_{\infty},\widetilde{g}^{(\infty)}_{ij}(\cdot,t)$) is
actually well defined on the time slice $M^4_{\infty}\times
\{t'\}$ and also has uniformly bounded curvature for all $t\in
[t',0]$. By taking a subsequence from the original subsequence
 and combining Step 1, we can extend the limit backward to a larger interval
 $[t^{''},0]\supsetneq (t',0]$. This is a contradiction with the definition of $t'$.

Therefore we have proved a subsequence of the rescaled solutions
($M^4_k$,\\
$\widetilde{g}^{(k)}_{ij}(\cdot,t),x_k$) converges to an ancient
$\kappa$-solution with restricted isotropic curvature pinching.
This is a contradiction.  We finish the proof of the theorem.
$$\eqno \#$$ \vskip 0.3cm

From now on, we always assume that the initial datum is a compact
four-manifold $M^4$ with no essential incompressible space form and
with positive isotropic curvature. Let $g_{ij}(x,t)$, $x\in M^4$ and
$t\in [0,T)$, be a maximal solution to the Ricci flow with
$T<+\infty$. Without loss of generality, after a scaling on the
initial metric, we may assume $T>1$. It was shown in \cite{Ha7} that
the solution $g_{ij}(x,t)$ remains positive isotropic curvature. By
Lemma 2.1, there hold (2.1), (2.2) and (2.3) for some positive
constants $0<\rho, \Lambda, P<+\infty$ (depending only on the
initial datum). And by Perelman's no local collapsed theorem I
\cite{P1} the solution is $\kappa$-noncollapsed on the scale
$\sqrt{T}$ for some $\kappa>0$ (depending only on the initial
datum). Then for any sufficiently small $\varepsilon>0$, we can find
$r_0>0$ with the property described in Theorem 4.1.

Let $\Omega$ denote the set of all points in $M^4$, where
curvature stays bounded as $t\rightarrow T$. The estimates (4.1)
imply that $\Omega$ is open and $R(x,t)\rightarrow+\infty$ as
$t\rightarrow T$ for each $x\in M^4\backslash \Omega$. If $\Omega$
is empty, then the solution becomes extinct at time $T$ and the
manifold is either diffeomorphic to $\mathbb{S}^4$ or
$\mathbb{R}\mathbb{P}^4$, or entirely covered by evolving
$\varepsilon$-necks or evolving $\varepsilon$-caps shortly before
the maximal time $T$, so $M^4$ is diffeomorphic to $\mathbb{S}^4$,
or $\mathbb{R}\mathbb{P}^4$, or $\mathbb{R}\mathbb{P}^4\#
\mathbb{R}\mathbb{P}^4$ or $\mathbb{S}^3\times \mathbb{S}^1$, or
$\mathbb{S}^3\widetilde{\times} \mathbb{S}^1$. The reason is as
follows. We only need to consider the situation that the manifold
$M^{4}$ is entirely covered by evolving $\varepsilon$-necks and
evolving $\varepsilon$-caps shortly before the maximal time $T$.
If $M^{4}$
 contains a cap $C$, then there is a cap or a neck adjacent to the
 neck like end of $C$. The former case implies that $M^{4}$ is
 diffeomorphic to $\mathbb{S}^4$, $\mathbb{RP}^4$, or
$\mathbb{RP}^4\#\mathbb{RP}^4$. In the latter case, we get a new
longer cap and continue the procedure. Finally, we must end up
with a cap, producing a $\mathbb{S}^4$, $\mathbb{RP}^4$, or
$\mathbb{RP}^4\#\mathbb{RP}^4$. If $M^{4}$ contains no caps, we
start with a neck $N$, consider the other necks adjacent to the
boundary of $N$, this gives a longer neck and we continue the
procedure. After a finite number of steps, the neck must repeat
itself. By considering the orientation of $M^{4}$, we conclude
that $M^{4}$ is diffeomorphic to $\mathbb{S}^3\times \mathbb{S}^1$
or $\mathbb{S}^3\widetilde{\times} \mathbb{S}^1$.

 We can now
assume that $\Omega$ is not empty. By using the local derivative
estimates of Shi \cite{Sh1} (or see \cite{Ha6}), we see that as
$t\rightarrow T$, the solution $g_{ij}(\cdot,t)$ has a smooth
limit $\overline{g}_{ij}(\cdot)$ on $\Omega$. Let
$\overline{R}(x)$ denote the scalar curvature of
$\overline{g}_{ij}$. By the positive isotropic curvature
assumption on the initial metric, we know that the metric
$\overline{g}_{ij}(\cdot)$ also has positive isotropic curvature;
in particular, $\overline{R}(x)$ is positive. For any
$\sigma<r_0$, let us consider the set $$\Omega_{\sigma}=\{x\in
\Omega \ |\ \overline{R}(x)\leq \sigma^{-2}\}.$$ Note that for any
fixed $x\in \partial \Omega$, as $x_j\in \Omega$ and
$x_j\rightarrow x$ with respect to the initial metric
$g_{ij}(\cdot,0)$, we have $\overline{R}(x_j)\rightarrow +\infty$.
In fact, if there was a subsequence $x_{j_k}$ so that the limit
$\lim_{k\rightarrow \infty}\overline{R}(x_{j_k})$ exists and is
finite, then it would follow from the gradient estimates (4.1)
that $\overline{R}$ is uniformly bounded in some small
neighborhood of $x\in \partial \Omega$ (with respect to the
induced topology of the initial metric $g_{ij}(\cdot,0)$); this is
a contradiction. From this observation and the compactness of the
initial manifold, we see that $\Omega_{\sigma}$ is compact (with
respect to the metric $\overline{g}_{ij}(\cdot)$).

For the further discussion, we follow \cite{P2} to introduce the
following terminologies. Denote by $\mathbb{I}$ a (finite or
infinite) interval.

Recall that an $\textbf{$\varepsilon$-neck}$ (of radius $r$) is an
open set with a Riemannian metric, which is, after scaling the
metric with factor $r^{-2}$, $\varepsilon$-close (in
$C^{[\varepsilon^{-1}]}$ topology) to the standard neck
$\mathbb{S}^3\times \mathbb{I}$ with the product metric, where
$\mathbb{S}^3$ has constant scalar curvature one and $\mathbb{I}$
has length $2\varepsilon^{-1}$. A metric on $\mathbb{S}^3\times
\mathbb{I}$, such that each point is contained in some
$\varepsilon$-neck, is called an \textbf{$\varepsilon$-tube}, or an
\textbf{$\varepsilon$-horn}, or a \textbf{double
$\varepsilon$-horn}, if the scalar curvature stays bounded on both
ends, or stays bounded on one end and tends to infinity on the
other, or tends to infinity on both ends, respectively. A metric on
$\mathbb{B}^4$ or $\mathbb{RP}^4\backslash \overline{\mathbb{B}^4}$,
such that each point outside some compact subset is contained in an
$\varepsilon$-neck, is called an \textbf{$\varepsilon$-cap} or a
\textbf{capped $\varepsilon$-horn}, if the scalar curvature stays
bounded or tends to infinity on the end, respectively.

Now take any $\varepsilon$-neck in $(\Omega, \overline{g}_{ij})$
and consider a point $x$ on one of its boundary components. If
$x\in \Omega\backslash \Omega_{\sigma}$, then there is either an
$\varepsilon$-cap or an $\varepsilon$-neck, adjacent to the
initial $\varepsilon$-neck. In the latter case we can take a point
on the boundary of the second $\varepsilon$-neck and continue.
This procedure can either terminate when we get into
$\Omega_{\sigma}$ or an $\varepsilon$-cap, or go on infinitely,
producing an $\varepsilon$-horn. The same procedure can be
repeated for the other boundary component of the initial
$\varepsilon$-neck. Therefore, we conclude that each
$\varepsilon$-neck of $(\Omega,\overline{g}_{ij})$ is contained in
a subset of $\Omega$ of one of the following types:
$$
\begin{array}{lllll}
\mbox{(a) an $\varepsilon$-tube with boundary components in
$\Omega_{\sigma}$, or}\\
\mbox{(b) an $\varepsilon$-cap with
boundary in $\Omega_{\sigma}$, or}\\
\mbox{(c) an $\varepsilon$-horn with boundary in
$\Omega_{\sigma}$, or \ \ \ \ \ \ \ \ \ \ \ \ \ \ \ \ \ \ \ \ \ \ \ \ \ \ \ \ \ (4.3)}\\
\mbox{(d) a capped $\varepsilon$-horn, or}\\
\mbox{(e) a double $\varepsilon$-horn.}
 \end{array}$$
Similarly, each $\varepsilon$-cap of $(\Omega,\overline{g}_{ij})$
is contained in a subset of $\Omega$ of either type (b) or type
(d).

It is clear that there is a definite lower bound (depending on
$\sigma$) for the volume of subsets of types (a), (b), (c), so
there can be only finite number of them. Thus we conclude that
there is only a finite number of components of $\Omega$,
containing points of $\Omega_{\sigma}$, and every such component
has a finite number of ends, each being an $\varepsilon$-horn.
While by taking into account that $\Omega$ has no compact
components, every component of $\Omega$, containing no points of
$\Omega_{\sigma}$, is either a capped $\varepsilon$-horn, or a
double $\varepsilon$-horn. Nevertheless, if we look at the
solution for a slightly earlier time $t$, each $\varepsilon$-neck
or $\varepsilon$-cap of $(M,g_{ij}(\cdot,t))$ is contained in a
subset of types (a) and (b); while the $\varepsilon$-horns, capped
$\varepsilon$-horns and double $\varepsilon$-horns, observed at
the maximal time $T$, are connected together to form
$\varepsilon$-tubes and $\varepsilon$-caps at the slightly earlier
time $t$.

Hence, by looking at the solution for times just before $T$, we see
that the topology of $M^4$ can be reconstructed as follows: take the
all components $\Omega_j,1\leq j\leq k$, of $\Omega$ which contains
points of $\Omega_{\sigma}$, truncate their $\varepsilon$-horns, and
glue a finite collection of tubes $\mathbb{S}^{3}\times \mathbb{I}$
and caps $\mathbb{B}^4$ or $\mathbb{RP}^4\backslash
\overline{\mathbb{B}^4}$ to the boundary components of truncated
$\Omega_j$. Thus $M^4$ is diffeomorphic to a connected sum of
$\overline{\Omega}_j$, $1\leq j\leq k$, with a finite number of
$\mathbb{S}^3\times \mathbb{S}^1$ or
$\mathbb{S}^3\widetilde{\times}\mathbb{S}^1$ (which correspond to
glue a tube to two boundary components of the same $\Omega_j$), and
a finite number of $\mathbb{RP}^4$. Here $\overline{\Omega}_j$
denotes $\Omega_j$ with each $\varepsilon$-horn one point
compactified. (One might wonder why we do not also cut other
$\varepsilon$-tubes or $\varepsilon$-caps so that we can remove more
volumes; we will explain it a bit later.)

More geometrically, one can get $\overline{\Omega}_j$ in the
following way: in every $\varepsilon$-horn of $\Omega_j$ one can
find an $\varepsilon$-neck, cut it along the middle three-sphere,
remove the horn-shaped end, and glue back a cap (i.e., a
differentiable four-ball). Thus to understand the topology of
$M^4$, one only need to understand the topologies of the compact
four-manifolds $\overline{\Omega}_j, 1\leq j\leq k$.

Recall that the four-manifold $M^4$ has no essential
incompressible space form, we now claim that each
$\overline{\Omega}_j$
 still has no essential incompressible space form. Clearly, we only need to check the assertion that if $N$ is an essential incompressible
 space form  in $\overline{\Omega}_j$, then $N$ will be
 also incompressible in $M^4$. After moving $N$ slightly, we can choose $N$ such
 that $N\subset \Omega_{j}$. Then $N$ can be regarded as a
 submanifold in $M^4$ (unaffected by the surgery). We now argue by contradiction. Suppose $\gamma\subset N$ is a homotopically nontrivial curve which bounds a
 disk $D$ in $M^4$. We want to modify the map of disk $D$ so that $\gamma$
 bound a new disk in $\Omega_{j}$, which will gives the desired contradiction. Let $E_{1},E_{2},\cdots, E_{m}$ be all the $\varepsilon$-horn ends of $\Omega_{j}$, $S_{1},S_{2},\cdots, S_{m}\subset
 \Omega_{j}$ be the corresponding cross spheres lying
 inside the $\varepsilon$-horn ends
 $E_{j}$ respectively.  Let us perturb the spheres $S_{1},S_{2},\cdots, S_{m}$
  slightly so that they meet $D$ transversely in a finite number of
  simple closed curves (we only consider those $S_{j}$ with $S_{j}\cap D\neq \phi$).
  After removing those curves which are contained in larger ones
  in $D$, we are left with a finite number of disjoint simple
  closed curves, denoted by $C_{1},C_{2},\cdots,C_{l}$. We denote
  the enclosed disks of $C_{1},C_{2},\cdots,C_{l}$ in $D$ by $D_{1},D_{2},\cdots,D_{l}$. Since
  $\mathbb{S}^{3}$ is simply-connected, each intersection curves in
  $S_{1},S_{2},\cdots,S_{m}$ can be shrunk to a point. So by filling the holes
$D_{1},D_{2},\cdots,D_{l}$, we obtain a new continuous map from
$D$ to $M^4$ such that
  the image of $D_{1}\cup D_{2}\cdots\cup D_{l}$ is contained in $S_{1}\cup S_{2}\cdots
  \cup S_{m}\subset \Omega_{j}$. On the other hand, since $D\backslash(D_{1}\cup D_{2}\cdots\cup D_{l})$
   is connected, $\gamma( \mbox{the image of}\  \partial D)\subset N$, we know that the image of $D\backslash(D_{1}\cup D_{2}\cdots\cup
   D_{l})$ must be contained in $\Omega_{j}$. Therefore, $\gamma$
   bounds a new disk in $\Omega_{j}$. This proves that after the surgery, each $\overline{\Omega}_j$
 still has no essential incompressible space form.

 As shown by
Hamilton in Section D of \cite{Ha7}, provided $\varepsilon>0$ small
enough, one can perform the above surgery procedure carefully so
that the compact four-manifolds $\overline{\Omega}_j,1\leq j\leq k$,
also have positive isotropic curvature. Naturally, one can evolve
each $\overline{\Omega}_j$ by the Ricci flow again and carry out the
same surgery procedure to produce a finite collection of new compact
four-manifolds with no essential incompressible space form and with
positive isotropic curvature. By repeating this procedure
indefinitely, it will be likely to give us the long time existence
of a kind of ``weak" solution to Ricci flow.\vskip 1cm

 \centerline{\large{\textbf{5. Ricci Flow with Surgery for
Four-manifolds}}} \vskip 0.5cm We begin with an abstract definition
of the solution to the Ricci flow with surgery which is adapted from
\cite{P2}.\vskip 0.3cm

$\underline{\mbox{\textbf{Definition 5.1}}}$ \emph{ Suppose we have
a collection of compact four-dimensional smooth solutions
$g^{(k)}_{ij}(t)$ to the Ricci flow on $M^4_k\times [t^-_k,t^+_k)$
with no essential incompressible space form and with positive
isotropic curvature, which go singular as $t\rightarrow t^+_k$ and
where each manifold $M^4_k$ may be disconnected with only a finite
number of connected components. Let
$(\Omega_k,\overline{g}^{(k)}_{ij})$ be the limits of the
corresponding solutions $g^{(k)}_{ij}(t)$ as $t\rightarrow t^+_k$.
Suppose also that for each $k$ we have $t^-_k=t^+_{k-1}$, and
$(\Omega_{k-1},\overline{g}^{(k-1)}_{ij})$ and
$(M^4_k,g^{(k)}_{ij}(t^-_k))$ contain compact (possibly
disconnected) four-dimensional submanifolds with smooth boundary
which are isometric. Then by identifying these isometric
submanifolds, we say it is a solution to the \textbf{Ricci flow with
surgery} on the time interval which is the union of all
$[t^-_k,t^+_{k})$, and say the times $t^+_k$ are \textbf{surgery
times}.}\vskip0.2cm
 The procedure described in the last paragraph of the previous
 section gives us a solution to the Ricci flow with surgery. However, in
 order to understand the topology of the initial manifold from the solution to
the Ricci flow with surgery, one encounters the following two
difficulties:

(i) how to prevent the surgery times from accumulation?

(ii) how to get the long time behavior of the solution to the
Ricci flow with surgery?

In view of this, it is natural to consider those solutions having
"good" properties. Let $\varepsilon$ be a fixed small positive
number. We will only consider those solutions to the Ricci flow with
surgery which satisfy the following a priori assumptions (with
accuracy $\varepsilon$):\vskip 0.3cm

$\underline{\mbox{\textbf{Pinching assumption}}}$: \emph{There exist
positive constants $\rho,\Lambda,P<+\infty$ such that there hold
$$a_1+\rho>0 \mbox{ and } c_1+\rho>0, \eqno (5.1)$$ $$\max\{a_3,b_3,c_3\}
\leq \Lambda(a_1+\rho) \mbox{ and }\max\{a_3,b_3,c_3\}\leq \Lambda(c_1+\rho), \eqno
(5.2)$$ and $$\frac{b_3}{\sqrt{(a_1+\rho)(c_1+\rho)}}\leq
1+\frac{\Lambda
e^{Pt}}{\max\{\log\sqrt{(a_1+\rho)(c_1+\rho)},2\}}, \eqno (5.3)$$
 everywhere.} \vskip 0.3cm

$\underline{\mbox{\textbf{Canonical neighborhood assumption (with
accuracy $\varepsilon$)}}}$: \emph{For the given $\varepsilon>0$,
there exist two constants $C_1(\varepsilon)$, $C_2(\varepsilon)$ and
a non-increasing positive function $r$ on $[0,+\infty)$ such that
for every point $(x,t)$ where the scalar curvature $R(x,t)$ is at
least $r^{-2}(t)$, there is an open neighborhood $B$,
$B_t(x,\sigma)\subset   B\subset B_t(x,2\sigma)$ with
$0<\sigma<C_1(\varepsilon)R(x,t)^{-\frac{1}{2}}$,
 which falls into
one of the following three categories:}

\emph{(a) $B$ is a \textbf{strong $\varepsilon$-neck} (in the
sense that $B$ is an $\varepsilon$-neck and it is the slice at
time $t$ of the parabolic neighborhood $\{(x',t')\ |\ x'\in B,
t'\in[t-R(x,t)^{-1},t]\}$, where the solution is well defined on
the whole parabolic neighborhood and is, after scaling with factor
$R(x,t)$ and shifting the time to zero, $\varepsilon$-close (in
$C^{[\varepsilon^{-1}]}$ topology) to the corresponding subset
 of the evolving
standard round cylinder $\mathbb{S}^3\times \mathbb{R}$ with
scalar curvature $1$ at the time zero), or}

 \emph{(b) $B$ is an \textbf{$\varepsilon$-cap}, or}

\emph{(c) $B$ is a compact four-manifold with positive curvature
operator;\\ furthermore, the scalar curvature in $B$ at time
 $t$ is between $C^{-1}_2R(x,t)$ and $C_2R(x,t)$, and satisfies
 the gradient estimate $$|\nabla R|<\eta R^{\frac{3}{2}}\mbox{ and }|
 \frac{\partial R}{\partial t}|<\eta R^2, \eqno
 (5.4)$$ and the volume of $B$ in case (a) and case (b) satisfies}
$$
(C_2R(x,t))^{-2} \leq Vol_t(B).
$$ \emph{Here $C_1$ and $C_2$ are some positive constants
 depending only on $\varepsilon$, and $\eta$ is a universal
 positive constant.} \vskip 0.3cm

 Clearly, we may always assume the above $C_1$ and $C_2$ are twice bigger
 than the corresponding constants $C_1(\frac{\varepsilon}{2})$ and
 $C_2(\frac{\varepsilon}{2})$ in Theorem 3.8 with the accuracy
 $\frac{\varepsilon}{2}$.

The main purpose of this section is to construct a long-time
 solution to the Ricci flow with surgery which starts with an
 arbitrarily given compact four-manifold with no essential
 incompressible space form and with positive isotropic curvature,
 so that the a priori assumptions are satisfied and there are only a finite number
 of surgery times at each finite time interval. The construction
 will be given by an induction argument.

 Firstly, for an arbitrarily given compact four-manifold
 $(M^4,g_{ij}(x))$ with no essential incompressible space form and
 with positive isotropic curvature, the Ricci flow with it as
 initial data has a maximal solution $g_{ij}(x,t)$ on $[0,T_0)$ with
 $T_0<+\infty$. Without loss of generality,
  after a scaling
 on the initial metric, we may assume $T_0>1$.
 It follows from Lemma 2.1 and Theorem 4.1 that the
 a priori assumptions above hold for the smooth solution on
 $[0,T_0)$.

 Suppose that we have a solution to the Ricci flow with surgery,
 with the given compact four-manifold $(M^4,g_{ij}(x))$ as initial
 datum,
 which is defined on $[0,T)$ with $T<+\infty$, going singular at
 the time $T$, satisfies the a priori assumptions and has only a
 finite number of surgery times on $[0,T)$. Let $\Omega$ denote
 the set of all points where the curvature stays bounded as
 $t\rightarrow T$. As shown before, the gradient estimate (5.4) in
 the canonical neighborhood assumption implies that $\Omega$ is
 open and that $R(x,t)\rightarrow +\infty$ as $t\rightarrow T$ for
 $x$ lying outside $\Omega$. Moreover, as $t\rightarrow T$, the
 solution $g_{ij}(x,t)$ has a smooth limit $\overline{g}_{ij}(x)$
 on $\Omega$.

 For $\delta>0$ to be chosen much smaller than $\varepsilon$, we
 let $\sigma=\delta r(T)$ where $r(t)$ is the positive
 nonincreasing function in the definition of the canonical
 neighborhood assumption. We consider the corresponding compact
 set $$\Omega_{\sigma}=\{x\in \Omega \ |\ \overline{R}(x)\leq
 \sigma^{-2}\}$$ where $\overline{R}(x)$ is the scalar curvature
 of $\overline{g}_{ij}$. If $\Omega_{\sigma}$ is empty, the
 manifold (near the maximal time $T$) is entirely covered by
 $\varepsilon$-tubes, $\varepsilon$-caps and compact components
 with positive curvature operator. Clearly, the number of
  compact components is finite. Then in this case the manifold
  (near the maximal time $T$) is
 diffeomorphic to the union of a finite number of $\mathbb{S}^4$,
 or $\mathbb{RP}^4$, or $\mathbb{S}^3\times
 \mathbb{S}^1$, or
 $\mathbb{S}^3\widetilde{\times }\mathbb{S}^1$, or a connected sum of them. Thus when
 $\Omega_{\sigma}$ is empty, the procedure stops here, and we say
 that the \textbf{solution becomes extinct}. We now assume
 $\Omega_{\sigma}$ is not empty. Every point $x\in
 \Omega\backslash \Omega_{\sigma}$ lies in one of subsets of
 listing in (4.3), or in a compact component with positive curvature operator or
 in a compact component which is contained in $\Omega\backslash \Omega_{\sigma}$ and is diffeomorphic to
 $\mathbb{S}^4$,
 or $\mathbb{RP}^4$, or $\mathbb{S}^3\times
 \mathbb{S}^1$, or
 $\mathbb{S}^3\widetilde{\times }\mathbb{S}^1$. Note again that the number of
  compact components is finite. Let us throw away all the compact
  components lying $\Omega\backslash \Omega_{\sigma}$ or with positive curvature
  operator, and then consider the all components $\Omega_j,1\leq j\leq
 k$, of $\Omega$ which contains points of $\Omega_{\sigma}$. (We will consider the components
 of $\Omega\backslash \Omega_{\sigma}$ consisting of
 capped $\varepsilon$-horns and double $\varepsilon$-horns later). We
 could perform Hamilton's surgerical procedure in Section D of
 \cite{Ha7} at every horn of $\Omega_j, 1\leq j\leq k$, so that the
 positive isotropic curvature condition and the pinching
 assumption is preserved.

 Note that if we perform the surgeries at the necks with certain fixed
accuracy  $\varepsilon$ on the high curvature region at each surgery
time, then it is possible that the errors of surgeries may
accumulate to a certain amount so that for some later time we can
not recognize the structure of very high curvature region. This
prevents us to carry out the process in finite time with finite
steps. Hence in order to maintain the a priori assumptions
\textbf{\emph{with the same accuracy}} after surgery, we need to
find sufficient ``fine" necks in the $\varepsilon$-horns and to glue
sufficient ``fine" caps in the procedure of surgery. Note that
$\delta>0$ will be chosen much smaller than $\varepsilon>0$. The
following lemma gives us the ``fine" necks in the
$\varepsilon$-horns. (The corresponding result in three-dimension
 is Lemma 4.3 in \cite{P2}).

 Now we explain that why we only perform the surgeries in
 the horns with boundary in $\Omega_{\sigma}$. At the first sight,
 we should also
cut off all those $\varepsilon$-tubes and $\varepsilon$-caps in the
surgery procedure. But in general, we are not able to find a
``finer" neck in an $\varepsilon$-tube or in $\varepsilon$-cap, and
such surgeries at ``rough" $\varepsilon$-necks will certainly loss
some accuracy. This is the reason why
 we will only perform the surgeries in the
$\varepsilon$-horns with boundary in $\Omega_{\sigma}$.\vskip 0.3cm

$\underline{\mbox{\textbf{Lemma 5.2}}}$\emph{  Given
$0<\varepsilon<\frac{1}{100}, 0<\delta<\varepsilon$ and
$0<T<+\infty$, there exists a radius $0<h< \delta \sigma$, depending
only on $\delta, r(T)$ and the pinching assumption, such that if we
have a solution to the Ricci flow with surgery,
 with a compact four-manifold $(M^4,g_{ij}(x))$ with no essential incompressible space form and
 with positive isotropic curvature as initial
 data,
 defined on $[0,T)$, going singular at
 the time $T$, satisfies the a priori assumptions and has only a
 finite number of surgery times on $[0,T)$, then for
each point $x$ with $h(x)=\overline{R}^{-\frac{1}{2}}(x)\leq h$ in
an $\varepsilon$-horn of $(\Omega,\overline{g}_{ij})$ with boundary
in $\Omega_{\sigma}$, the neighborhood
$B_T(x,\delta^{-1}h(x))=\{y\in
\Omega|dist_{\overline{g}_{ij}}(y,x)\leq \delta^{-1}h(x)\}$ is a
strong $\delta$-neck (i.e.,
 $\{(y,t)\ |\ y\in
B_T(x,\delta^{-1}h(x)),t\in[T-h^2(x),T]\}$ is, after scaling with
factor $h^{-2}(x)$, $\delta$-close (in $C^{[\delta^{-1}]}$ topology)
to the corresponding subset of the evolving standard round cylinder
$\mathbb{S}^3 \times \mathbb{R}$ over the time interval $[-1,0]$
with scalar curvature $1$ at the time zero).} \vskip 0.2cm

$\underline{\mbox{\textbf{Proof.}}}$ We argue as in \cite{P2} by
contradiction. Suppose that there exists a sequence of solution
$g^{(k)}_{ij}(\cdot,t), k=1,2,\cdots$, to the Ricci flow with
surgery, satisfying the a priori assumptions, defined on $[0,T)$
with limits $(\Omega^k,\overline{g}^{(k)}_{ij}), k=1,2,\cdots$, as
$t\rightarrow T$, and exist points $x_k$, lying inside an
$\varepsilon$-horn of $\Omega^k$, which contains the points of
$\Omega^k_{\sigma}$, and having $h(x_k)\rightarrow 0$ as
$k\rightarrow +\infty$ such that the neighborhood
$B_T(x_k,\delta^{-1}h(x_k))$ are not strong $\delta$-necks.

Let $\widetilde{g}^{(k)}_{ij}(\cdot,t)$ be the rescaled solutions
by the factor $\overline{R}(x_k)=h^{-2}(x_k)$ around $(x_k,T)$. We
will show that a sequence of $\widetilde{g}^{(k)}_{ij}(\cdot,t)$
converges to the evolving round $\mathbb{R}\times \mathbb{S}^3$,
which gives the desired contradiction.

Note that $\widetilde{g}^{(k)}_{ij}(\cdot,t), k=1, 2, \cdots,$ are
modified by surgery. We can not apply Hamilton's compactness theorem
directly since it states only for smooth solutions. For each
(unrescaled) surgical solution $\widetilde{g}^{(k)}_{ij}(\cdot,t)$,
we pick a point $z_k$, with $\bar{R}(z_k) =
2C_2^2(\varepsilon)\sigma^{-2},$ in the $\varepsilon$-horn of
$(\Omega^k,\bar{g}^{(k)}_{ij})$ with boundary in $\Omega_\sigma^k$,
where $C_2(\varepsilon)$ is the positive constant in the canonical
neighborhood assumption. From the definition of $\varepsilon$-horn
and the canonical neighborhood assumption, we know that each point
$x$ lying inside the $\varepsilon$-horn of
$(\Omega^k,\bar{g}^{(k)}_{ij})$ with
$d_{\bar{g}^{(k)}_{ij}}(x,\Omega_\sigma^k) \geq
d_{\bar{g}^{(k)}_{ij}}(z_k,\Omega_\sigma^k)$ has a strong
$\varepsilon$-neck as its canonical neighborhood. Since
$h(x_k)\rightarrow0$, each $x_k$ lies deeply inside the
$\varepsilon$-horn. Thus for each positive $A < +\infty$, the
rescaled (surgical) solutions $\widetilde{g}^{(k)}_{ij}(\cdot,t)$
with the marked origins $x_k$ over the geodesic balls
$B_{\widetilde{g}^{(k)}_{ij}(\cdot,0)}(x_k,A)$, centered at $x_k$ of
radii $A$ (with respect to the metrics
$\widetilde{g}^{(k)}_{ij}(\cdot,0)$), will be smooth on some uniform
(size) small time intervals for all sufficiently large $k$, if the
curvatures of the rescaled solutions $\widetilde{g}^{(k)}_{ij}$ at
$t=0$ in $B_{\widetilde{g}^{(k)}_{ij}(\cdot,0)}(x_k,A)$ are
uniformly bounded. In such situation, the Hamilton's compactness
theorem is applicable. Then we can now apply the same argument in
Step 2 of the proof Theorem 4.1 to conclude that the curvatures of
the rescaled solutions $\widetilde{g}^{(k)}_{ij}(\cdot,t)$ at the
time $T$ stay uniformly bounded at bounded distances from $x_k$;
otherwise we get a piece of a non-flat nonnegative curved metric
cone as a blow-up limit, which would contradict with Hamilton strong
maximum principle \cite{Ha2}. Hence as before we can get a
$C^{\infty}_{loc}$ limit $\widetilde{g}^{(\infty)}_{ij}(\cdot,t)$,
defined on a space-time set which is relatively open in the half
space-time $\{t\leq T\}$ and contains the time slice $\{t=T\}$, from
the rescaled solutions $\widetilde{g}^{(k)}_{ij}(\cdot,t)$.

By the pinching assumption, the limit is a complete manifold with
the restricted pinching condition (2.4) and with nonnegative
curvature operator. Since $x_k$ was contained in an
$\varepsilon$-horn with boundary in $\Omega^k_{\sigma}$, and
$h(x_k)/\sigma\rightarrow 0$, the limiting manifold has two ends.
Thus by Toponogov splitting theorem, it admits a (maybe not round
at this moment) metric splitting $\mathbb{R}\times \mathbb{S}^3$
because $x_k$ was the center of a strong $\varepsilon$-neck. We
further apply the restricted isotropic curvature pinching
condition (2.4) and contracted second Bianchi identity as before
to conclude that the factor $\mathbb{S}^3$ must be round at time
$0$. By combining with the canonical neighborhood assumption, we
see that the limit is defined on the time interval $[-1,0]$. By
Toponogov splitting theorem, the splitting $\mathbb{R}\times
\mathbb{S}^3$ is at each time $t\in [-1,0]$; so the limiting
solution is just the standard evolving round cylinder. This is a
contradiction. We finish the proof of Lemma 5.2.
$$\eqno \#$$ \vskip 0.3cm

The property in the above lemma that the radius $h$ depends only on
$\delta$, the time $T$ and the pinching assumption, independent of
the surgical solution, is crucial; otherwise we will not be able to
cut off enough volume at each surgery to guarantee the number of
surgeries being finite in each finite time interval.\vskip 0.3cm

{\textbf{Remark}}.  The proof of Lemma 5.2 actually proves a more
stronger result:\emph{ for any $\delta>0$, there exists a radius
$0<h< \delta \sigma$, depending only on $\delta, r(T)$ and the
pinching assumption, such that for each point $x$ with
$h(x)=\overline{R}^{-\frac{1}{2}}(x)\leq h$ in an $\varepsilon$-horn
of $(\Omega,\overline{g}_{ij})$ with boundary in $\Omega_{\sigma}$,
$\{(y,t)\ |\ y\in
B_T(x,\delta^{-1}h(x)),t\in[T-\delta^{-2}h^2(x),T]\}$ is, after
scaling with factor $h^{-2}(x)$, $\delta$-close (in
$C^{[\delta^{-1}]}$ topology) to the corresponding subset of the
evolving standard round cylinder $\mathbb{S}^3 \times \mathbb{R}$
over the time interval $[-\delta^{-2},0]$ with scalar curvature $1$
at the time zero.}
 This fact will be used in the proof of
Proposition 5.4.

The reason is as follows. Let us use the notation in the proof the
Lemma 5.2 and argue by contradiction. Note that the scalar curvature
of the limit at time $t=-1$ is $\frac{1}{1-\frac{2}{3}(-1)}$. Since
$h(x^k)/\rho\rightarrow0$, each point in the limiting manifold at
time $t=-1$ has also a strong $\varepsilon$-neck as its canonical
neighborhood. Thus the limit is defined at least on the time
interval $[-2,0]$. Inductively, suppose the limit is defined on the
time interval $[-m,0]$ with bounded curvature for some positive
integer $-m$, then by the isotropic pinching condition, Toponogov
splitting theorem and evolution equation of the scalar curvature on
the round $\mathbb{R}\times \mathbb{S}^3$, we see that
$R=\frac{1}{1+\frac{2}{3}m}$ at time $-m$. Since
$h(x^k)/\rho\rightarrow0$, each point in the limiting manifold at
time $t=-m$ has also a strong $\varepsilon$-neck as its canonical
neighborhood, we see that the limit is defined at least on the time
interval $[-(m+1),0]$ with bounded curvature. So by induction we
prove that the limit exists on the ancient time interval $(-\infty,
0]$. Therefore the limit is the evolving round cylinder
$\mathbb{S}^3\times\mathbb{R}$ over the time interval $(-\infty,
0]$, which gives the desired contradiction.

\vskip0.3cm

 To specialize our surgery, we now fix a standard capped infinite cylinder
 for $n=4$ as follows.
 Consider the semi-infinite standard round cylinder $N_0 =
 \mathbb{S}^3 \times (-\infty,4)$ with the metric $g_0$
 of scalar curvature 1. Denote by $z$ the coordinate of
 the second factor $(-\infty,4)$. Let $f$ be a
 smooth nondecreasing convex function on
$(-\infty,4)$ defined by
$$
   \left\{
   \begin{array}{lll}
    f(z) = 0, \ \ \ z\leq 0,
         \\[3mm]
    f(z) = ce^{-\frac{D}{z}}, \ \ \ z \in (0,3], \\[3mm]
    f(z) \mbox{ is strictly convex on } z \in [3,3.9], \\[3mm]
    f(z) = -\frac{1}{2}\log(16-z^2), \ \ \ z \in [3.9,4),
\end{array}
\right.
$$
where the small (positive) constant $c$ and big (positive)
constant $D$ will be determined later. Let us replace the standard
metric $g_0$ on the portion $\mathbb{S}^3 \times [0,4)$ of the
semi-infinite cylinder by $\hat{g} = e^{-2f}g_0$. Then the
resulting metric $\hat{g}$ will be smooth on $\mathbb{R}^4$
obtained by adding a point to $\mathbb{S}^3 \times (-\infty,4)$ at
$z=4$. We denote by $C(c,D) = (\mathbb{R}^4,\hat{g})$. Clearly,
$C(c,D)$ is a standard capped infinite cylinder.

 We next use a compact portion of the
standard capped infinite cylinder $C(c,D)$ and the $\delta$-neck
obtained in Lemma 5.2 to perform the following surgery  due to
Hamilton \cite{Ha7}.

Consider the solution metric $\bar{g}$ at the maximal time
$T<+\infty$. Take an $\varepsilon$-horn with boundary in
$\Omega_\rho$. By Lemma 5.2, there exists a $\delta$-neck $N$ of
radius $0<h<\delta \rho$ in the $\varepsilon$-horn. By definition,
$(N,h^{-2}\bar{g})$ is $\delta$-close (in $C^{[\delta^{-1}]}$
topology) to the standard round neck $\mathbb{S}^3\times
\mathbb{I}$ of scalar curvature 1 with
$\mathbb{I}=(-\delta^{-1},\delta^{-1})$. The parameter $z \in
\mathbb{I}$ induces a function on the $\delta$-neck $N$.

Let us cut the $\delta$-neck $N$ along the middle (topological)
three-sphere $N\bigcap\{z=0\}$. Without loss of generality, we may
assume that the right hand half portion $N\bigcap\{z\geq 0\}$ is
contained in the horn-shaped end. Let $\varphi$ be a smooth bump
function with $\varphi = 1$ for $z\leq2$, and $\varphi = 0$ for
$z\geq3$. Construct a new metric $\tilde{g}$ on a (topological)
four-ball $\mathbb{B}^4$ as follows
$$
   \tilde{g} = \left\{
   \begin{array}{lll}
    \bar{g}, \ \ \ z= 0,
         \\[4mm]
    e^{-2f}\bar{g}, \ \ \ z \in [0,2], \\[4mm]
    \varphi e^{-2f}\bar{g} + (1-\varphi)e^{-2f}h^2g_0, \ \ \ z \in [2,3], \\[4mm]
    h^2e^{-2f}g_0, \ \ \ z\in [3,4].
\end{array}
\right.
$$
The surgery is to replace the horn-shaped end by the cap
$(\mathbb{B}^4,\tilde{g})$. The following lemma, due to Hamilton
\cite{Ha7}, determines the constants $c$ and $D$ in the
$\delta$-cutoff surgery so that the pinching assumption is preserved
under the surgery.

\vskip 0.3cm \noindent \textbf{Lemma 5.3} ( Hamilton \cite{Ha7}
D3.1) (Justification of the pinching assumption)

\emph{ There are universal positive constants $\delta_0$, $c_0$
and $D_0$ such that for any $\tilde{T}$ there is a constant
$h_0>0$ depending on the initial metric and $\tilde{T}$ such that
if we take a $\delta$-cutoff surgery at a $\delta$-neck of radius
$h$ at time $T\leq\tilde{T}$ with $\delta < \delta_0$ and $h^{-2}
\geq h_0^{-2}$, then we can choose $c=c_0$ and $D=D_0$ in the
definition of $f(z)$ such that after the surgery, there still
holds the pinching condition (2.1) (2.2) (2.3):
$$
a_1+\rho>0 \mbox{ and }c_1+\rho>0, $$
$$\max\{a_3,b_3,c_3\}\leq \Lambda (a_1+\rho) \mbox{ and
}\max\{a_3,b_3,c_3\}\leq \Lambda (c_1+\rho), $$ and
$$\frac{b_3}{\sqrt{(a_1+\rho)(c_1+\rho)}}\leq
1+\frac{\Lambda e^{Pt}}{\max\{\log\sqrt{(a_1+\rho)(c_1+\rho)},2\}}
$$ at all points at time $T$. Moreover,
after the surgery, any metric ball of radius
$\delta^{-\frac{1}{2}}h$ with center near the tip (i.e. the origin
of the attached cap) is, after scaling with factor $h^{-2}$,
$\delta^{\frac{1}{2}}$-close the corresponding ball of the standard
capped infinite cylinder $C(c_0,D_0)$.}$$\eqno \#$$\vskip0.2cm
 We
call the above procedure as a $\textbf{$\delta$-cutoff surgery}$.
Since there are only finite number of horns with their other ends
connected to $\Omega_{\sigma}$, we only need to \emph{perform a
finite number of such $\delta$-cutoff surgeries at the time $T$}.
Besides those horns, there could be capped horns, double horns and
compact components lying $\Omega\setminus\Omega_{\sigma}$ or with
positive curvature operator. As explained before, capped horns and
double horns are connected with horns to form tubes or capped tubes
at any time slightly before $T$. Thus when we truncated the horns at
the $\delta$-cutoff surgeries, we actually had removed these
together with the horn-shaped ends away. So \emph{we can regard the
capped horns and double horns (of $\Omega\setminus\Omega_{\sigma}$)
to be extinct and throw them away at the time $T$}. Remember that
\emph{we have thrown away all the compact components lying in
$\Omega\setminus\Omega_{\sigma}$ or with positive curvature
operator}. Each of such compact components is diffeomorphic to
$\mathbb{S}^4$,
 or $\mathbb{RP}^4$, or $\mathbb{S}^3\times
 \mathbb{S}^1$, or
 $\mathbb{S}^3\widetilde{\times }\mathbb{S}^1$, and the number
of compact components is finite. Thus we actually throw a finite
number of $\mathbb{S}^4$, $\mathbb{RP}^4$, $\mathbb{S}^3\times
 \mathbb{S}^1$ or
 $\mathbb{S}^3\widetilde{\times }\mathbb{S}^1$ at the time $T$ also.
(Note that we allow that the manifold may be disconnected before and
after the surgeries). Let us agree to \textbf{\emph{declare extinct
every compact component with positive curvature operator or lying in
$\Omega\setminus\Omega_{\sigma}$}}; in particular, that allows to
exclude the components with positive curvature operator from the
list of canonical neighborhoods.\vskip 0.3cm

\emph{Summarily, our surgery at the time $T$ consists of the
following four procedures:}

 (1) \emph{perform $\delta$-cutoff surgeries for all $\varepsilon$-horns which
have the other ends connected to $\Omega_\sigma$,}

(2) \emph{declare extinct every compact component which has positive
curvature operator,}

(3) \emph{throw away all capped horns and double horns lying in
$\Omega \setminus\Omega_\sigma$,}

(4) \emph{declare extinct every compact components lying in $\Omega
\setminus\Omega_\sigma$.}

\vskip 0.5cm
 After the surgery at the time $T$, the pinching assumption still holds for the surgically modified
 manifolds. With this (maybe disconnected) surgically modified manifold as initial data, we now continue
 our solution until it becomes
singular for the next time $T'(>T)$. Therefore we have extended the
solution to the Ricci flow with surgery, originally defined on
$[0,T)$, to the new time interval $[0,T')$ (with $T'>T$). Moreover,
as long as $0<\delta\leq \delta_0$, the solution with
$\delta$-cutoff surgeries on the new time interval $[0,T')$ still
has positive isotropic curvature and no essential incompressible
space form, and from \cite{Ha7} and Lemma 5.3 it still satisfies the
pinching assumption.

 Denote the minimum of the
scalar curvature at time $t$ by $R_{\min}(t)>0$. Since the
$\delta$-cutoff surgeries occur at the points lying deeply in the
$\varepsilon$-horns, the minimum of the scalar curvature
$R_{min}(t)$ of the solution at each time-slice is achieved in the
region unaffected by the surgeries. Thus we know from the
evolution equation of the scalar curvature that
$$\frac{d}{dt}R_{min}(t)\geq \frac{1}{2}R^2_{min}(t).$$ By
integrating this inequality, we conclude that the maximal time $T$
 of any solution to the Ricci flow with $\delta$-cutoff surgeries
must be bounded by $2/R_{\min}(0)<+\infty$.
 Let
$\tilde{T}=2/R_{\min}(0)$ in Lemma 5.3, then there is a constant
$h_{0}$ determined by $\tilde{T}$. Set
$\bar{\delta}=\frac{1}{2}R_{\min}(0)^{\frac{1}{2}}h_{0}.$ We know
that if we perform the $\delta$-cutoff surgery with
$\delta<\min\{\bar{\delta},\delta_0\}$, then the pinching
assumptions (5.1),(5.2),(5.3) are satisfied for the solution to the
Ricci flow with $\delta$-cutoff surgery. Next we make further
restrictions on $\delta$ to justify the canonical neighborhood
assumption. Clearly, we only need to check the following assertion.
\vskip 0.3cm

$\underline{\mbox{\textbf{Proposition 5.4}}}$ (Justification of
the canonical neighborhood assumption) \emph{ Given a compact
four-manifold with positive isotropic curvature and no essential
incompressible space form and given $\varepsilon > 0$, there exist
decreasing sequences
  $\varepsilon>\widetilde{r}_j>0$, $\kappa_j>0$, $\min
\{\varepsilon^2,\delta_0,\bar{\delta}\}>\widetilde{\delta}_j>0$,
$j=1,2,\cdots$, with the following property.  Define a positive
function $\widetilde{\delta}(t)$ on $[0,+\infty)$ by
$\widetilde{\delta}(t) = \widetilde{\delta}_j$ when $t \in
[(j-1)\varepsilon^2,j\varepsilon^2)$. Suppose we have a solution to
the Ricci flow with surgery,
 with the given
four-manifold as initial datum defined on the time interval $[0,T)$
and with a finite number of $\delta$-cutoff surgeries such that any
$\delta$-cutoff surgery at a time $t \in (0,T)$ with
$\delta=\delta(t)$ satisfies $0<\delta(t) \leq
\widetilde{\delta}(t).$ Then on each the time interval
$[(j-1)\varepsilon^2,j\varepsilon^2]\bigcap [0,T)$, the solution
satisfies the $\kappa_j$-noncollapsing condition on all scales less
than $\varepsilon$ and the canonical neighborhood assumption (with
accuracy $\varepsilon$) with $r=\widetilde{r}_j$.} \vskip 0.2cm

Here and in the followings, we call a (four-dimensional) solution
$g_{ij}(t), 0 \leq t <T$, to the Ricci flow with surgery is
\textbf{$\kappa$-noncollapsed} at a point $(x_0,t_0)$ on the scales
less than $\rho$ (for some $\kappa>0,\rho>0$) if it satisfies the
following property: whenever $r < \rho$ and
$$ |Rm(x,t)| \leq r^{-2}$$  for all those $(x,t) \in P(x_0,t_0,r,-r^2)=\{ (x',t') \
|\ x'\in B_{t'}(x_0,r), t' \in [t_0 -r^2,t_0] \}$, for which the
solution is defined, we have $$Vol_{t_0}(B_{t_0}(x_0,r)) \geq \kappa
r^4.$$ Before we give the proof of the proposition, we need to check
$\kappa$-noncollapsing condition. \vskip 0.3cm

$\underline{\mbox{\textbf{Lemma 5.5}}}$ \emph{For a given compact
four-manifold with positive isotropic curvature and no essential
incompressible space form and given $\varepsilon>0$, suppose we have
constructed the sequences, satisfying the above proposition for
$1\leq j\leq \ell$. Then there exists $\kappa>0$, such that for any
$r$, $0<r<\varepsilon$, one can find $\widetilde{\delta}$ with
$0<\widetilde{\delta}<\min\{\varepsilon^2,\delta_0,\bar{\delta}\}$,
which depends on $r$, $\varepsilon$ and may also depend on the
already constructed sequences, with the following property. Suppose
we have a solution, with the given four-manifold as initial data, to
the Ricci flow with surgery defined on a time interval $[0,T]$ with
$\ell \varepsilon^{2}\leq T<(\ell+1)\varepsilon^2$ such that the
assumptions and conclusions of Proposition 5.4 hold on $[0,\ell
\varepsilon^2)$, the canonical neighborhood assumption (with
accuracy $\varepsilon$) with  $r$ holds on $[\ell \varepsilon^2,T]$,
and each $\delta(t)$-cutoff surgery in the time interval
$t\in[(\ell-1)\varepsilon^2,T]$ has $0<\delta(t)<\tilde{\delta}$.
Then the solution is $\kappa$-noncollapsed on $[0,T]$ for all scales
less than $\varepsilon$.}\vskip 0.2cm

$\underline{\mbox{\textbf{Proof}}}$. Consider a parabolic
neighborhood $P(x_0,t_0,r_0,-r^2_0)=\{(x,t)|x\in B_t(x_0,r_0),t\in
[t_0-r^2_0,t_0]\}$, with $\ell\varepsilon^2\leq t_0\leq T$, and $
0<r_0\leq \varepsilon$, where the solution satisfies $|Rm|\leq
r^{-2}_0$ whenever it is defined. We will prove that
$Vol_{t_0}(B_{t_0}(x_0,r_0))\geq  \kappa r^4_0$.

Let $\eta$ be the universal positive constant in the definition of
the canonical neighborhood assumption. Without loss of generality,
we always assume $\eta \geq 10$. Firstly, we want to show that one
may assume $r_0 \geq \frac{1}{2\eta}r$.

Obviously, the curvature satisfies the estimate
$$|Rm(x,t)|\leq 20r^{-2}_0,$$
for those $(x,t) \in
P(x_0,t_0,\frac{1}{2\eta}r_0,-\frac{1}{8\eta}r^2_0) = \{(x,t)\ |\
x\in
B_t(x_0,\frac{1}{2\eta}r_0),t\in[t_0-\frac{1}{8\eta}r^2_0,t_0]\}$,
for which the solution is defined. When $r_0< \frac{1}{2\eta}r$, we
can enlarge $r_0$ to some $r'_0\in[r_0,r]$ so that
$$|Rm|\leq 20r'^{-2}_0$$
 on $P(x_0,t_0,\frac{1}{2\eta}r'_0,-\frac{1}{8\eta}r'^2_0)$ (whenever
it is defined), and either the equality holds somewhere or $r'_0=r$.

In the case that the equality holds somewhere, it follows from the
pinching assumption that we have
$$ R > 10r'^{-2}_0$$
somewhere in
$P(x_0,t_0,\frac{1}{2\eta}r'_0,-\frac{1}{8\eta}r'^2_0)$. Here,
without loss of generality, we have assumed $r$ is suitably small.
Then by the gradient estimates in the definition of the canonical
neighborhood assumption, we know $$R(x_0,t_0) > r'^{-2}_0 \geq
r^{-2}.$$ Hence the desired noncollapsing estimate in this case
follows directly from the canonical neighborhood assumption. (Recall
that we have excluded every component which has positive sectional
curvature in the surgery procedure and then we have excluded them
from the list of canonical neighborhoods. Here we also used the
standard volume comparison when the canonical neighborhood is an
$\varepsilon$-cap).

While in the case that $r'_0 = r$, we have the curvature bound
$$|Rm(x,t)|\leq (\frac{1}{2\eta}r)^{-2},$$
for those $(x,t) \in
P(x_0,t_0,\frac{1}{2\eta}r,-(\frac{1}{2\eta}r)^{2}) = \{(x,t)\ |\
x\in
B_t(x_0,\frac{1}{2\eta}r),t\in[t_0-(\frac{1}{2\eta}r)^2,t_0]\}$, for
which the solution is defined. It follows from the standard volume
comparison that we only need to verify the noncollapsing estimate
for $r_0 = \frac{1}{2\eta}r$. Thus we have reduced the proof to the
case $r_0\geq \frac{1}{2\eta}r$.

The reduced distance from $(x_0,t_0)$ is
$$l(q,\tau)=\frac{1}{2\sqrt{\tau}}\inf\{\int^{\tau}_0\sqrt{s}
(R(\gamma(s),t_0-s)+|\dot{\gamma}(s)|^2_{g_{ij}(t_0-s)})ds\mid\gamma(0)=x_{0},
\gamma(\tau)=q\}$$ where $\tau=t_0-t$ with $t<t_0$. Firstly, we need
to check that the minimum of the reduced distance is achieved by
curves unaffected by surgery. According to Perelman \cite{P2}, we
call a space-time curve in the solution track is
$\textbf{admissible}$ if it stays in the space-time region
unaffected by surgery, and we call a space-time curve in the
solution track is a $\textbf{barely admissible curve}$ if it is on
the boundary of the set of admissible curves. The following
assertion gives a big lower bound for the reduced lengths of barely
admissible curves. \vskip 0.2cm

$\underline{\mbox{\textbf{Claim 1}}}$. \emph{For any $L<+\infty$ one
can find
$\widetilde{\delta}=\widetilde{\delta}(L,r,\widetilde{r}_{\ell},\varepsilon)>0$
with the following property. Suppose that we have a curve $\gamma$,
parametrized by $t\in [T_0,t_0], (\ell-1)\varepsilon^2\leq T_0<t_0$,
such that $\gamma(t_0)=x_0$, $T_0$ is a surgery time and
$\gamma(T_0)$ lies in a $4h$-collar of the middle three-sphere of a
$\delta$-neck with the radius $h$ obtained in Lemma 5.2, where the
$\delta$-cutoff surgery was taken. Suppose also each
$\delta(t)$-cutoff surgery in the time interval
$t\in[(\ell-1)\varepsilon^2,T]$ has $0<\delta(t)<\tilde{\delta}$.
Then we have an estimate
$$\int^{t_0-T_0}_0\sqrt{\tau}(R(\gamma(t_0-\tau),t_0-\tau)+
|\dot{\gamma}(t_0-\tau)|^2_{g_{ij}(t_0-\tau)})d\tau\geq L, \eqno
(5.5)$$ where $\tau=t_0-t\in [0,t_0-T_0]$. } \vskip 0.2cm

Before we can verify this assertion, we need to do some premilary
works.

Let $O$ be the point near $\gamma(T_0)$ which corresponds to the
center of the (rotationally symmetric) capped infinite round
cylinder. Recall from Lemma 5.3 that a metric ball of radius
$\delta^{-\frac{1}{2}}h$ at time $T_0$ centered at $O$ is, after
scaling with factor $h^{-2}$, $\delta^{\frac{1}{2}}$-close (in
$C^{[\delta^{-\frac{1}{2}}]}$ topology) to the corresponding ball in
the capped infinite round cylinder. We need to consider the
solutions to the Ricci flow with the capped infinite round cylinder
(with scalar curvature 1 outsider some compact set) as initial data
and we require the solutions have also bounded curvature; we call
such a solution a \textbf{standard solution} as in \cite{P2}. From
Shi \cite{Sh1}, we know such a solution exists. The uniqueness of
the Ricci flow for compact manifolds is well-known (see for example,
Section 6 of \cite{Ha6}). In \cite{CZ2}, we prove a uniqueness
theorem which states that if the initial data is a complete
noncompact Riemannian manifold with bounded curvature, then the
solution to the Ricci flow in the class of complete solutions with
bounded curvature is unique. Thus the standard solution with a
capped infinite round cylinder as initial data is unique. In the
appendix, we will show that the standard solution exists on the time
interval $[0,\frac{3}{2})$ and has nonnegative curvature operator,
and its scalar curvature satisfies
$$R(x,t)\geq
\frac{C^{-1}}{\frac{3}{2}-t}, \eqno (5.6)$$ everywhere for some
positive constant $C$.

For any $0<\theta<\frac{3}{2}$, let $Q$ be the maximum of the
scalar curvature of the standard solution in the time interval
$[0,\theta]$ and let $\triangle
t=(T_1-T_0)/N<\varepsilon\eta^{-1}Q^{-1}h^2$ with
$T_1=\min\{t_0,T_0+\theta h^2\}$ and $\eta$ given in the canonical
neighborhood assumption. Set $t_k=T_0+k\triangle t,k=1,\cdots,N$.

Note that the ball $B_{T_0}(O,A_0h)$ at time $T_0$ with
$A_0=\delta^{-\frac{1}{2}}$ is, after scaling with factor $h^{-2}$,
$\delta^{\frac{1}{2}}$-close to the corresponding ball in the capped
infinite round cylinder. Assume first that for each point in
$B_{T_0}(O,A_0h)$, the solution is defined on $[T_0,t_1]$. By the
gradient estimate (5.4) in the canonical neighborhood assumption and
the choice of $\triangle t$ we have a uniform curvature bound on
this set for $h^{-2}$-scaled metric. Then by the uniqueness theorem
in \cite{CZ2}, if $\delta^{\frac{1}{2}}\rightarrow 0$ (i.e.,
$A_0=\delta^{-\frac{1}{2}}\rightarrow +\infty$), the solution with
$h^{-2}$-scaled metric will converge to the standard solution in
$C^{\infty}_{loc}$ topology. Therefore we can find $A_1$, depending
only on $A_0$ and tending to infinity with $A_0$, such that the
solution in the parabolic region $P(O,T_0,A_1h,t_1-T_0)=\{(x,t)|x\in
B_t(O,A_1h),t\in [T_0,T_0+(t_1-T_0)]\}$ is, after scaling with
factor $h^{-2}$ and shifting time $T_0$ to zero, $A^{-1}_1$-close to
the corresponding subset in the standard solution. In particular,
the scalar curvature on this subset does not exceed $2Qh^{-2}$. Now
if each point in $B_{T_0}(O,A_1h)$ the solution is defined on
$[T_0,t_2]$, then we can repeat the procedure, defining $A_2$, such
that the solution in the parabolic region
$P(O,T_0,A_2h,t_2-T_0)=\{(x,t)|x\in B_t(p,A_2h),t\in
[T_0,T_0+(t_2-T_0)]\}$ is, after scaling with factor $h^{-2}$ and
shifting time $T_0$ to zero, $A^{-1}_2$-close to the corresponding
subset in the standard solution. Again, the scalar curvature on this
subset still does not exceed $2Qh^{-2}$. Continuing this way, we
eventually define $A_N$. Note that $N$ is depending only on
$\theta$. Thus for arbitrarily given $A>0$ (to be determined), we
can choose $\widetilde{\delta}(A,\theta, \varepsilon)>0$ such that
as $\delta<\widetilde{\delta}(A,\theta, \varepsilon)$, and assuming
that for each point in $B_{T_0}(O,A_{(N-1)}h)$ the solution is
defined on $[T_0,T_1]$, we have $A_0>A_1>\cdots>A_N>A$, and the
solution in $P(O,T_0,Ah,T_1-T_0)=\{(x,t)|x\in B_t(O,Ah),t\in
[T_0,T_1]\}$ is, after scaling with factor $h^{-2}$ and shifting
time $T_0$ to zero, $A^{-1}$-close to the corresponding subset in
the standard solution.

Now assume that for some $k$ $(1\leq k\leq N-1)$ and a surgery time
$t^+\in(t_k,t_{k+1}]$(or $ t^{+}\in (T_{0},t_{1}]$) such that on
$B_{T_0}(O,A_kh)$ the solution is defined on $[T_0,t^+)$, but for
some point of this ball it is not defined past $t^+$. Clearly the
above argument also shows that the parabolic region
$P(O,T_0,A_{k+1}h,t^+-T_0)=\{(x,t)|x\in B_t(x,A_{k+1}h),t\in
[T_0,t^+)\}$ is, after scaling with factor $h^{-2}$ and shifting
time $T_0$ to zero, $A^{-1}_{k+1}$-close to the corresponding subset
in the standard solution. In particular, as the time tends to $t^+$,
the ball $B_{T_0}(O,A_{k+1}h)$ keeps on looking like a cap. Since
the scalar curvature on the set $B_{T_0}(O,A_{k}h)\times [T_0,t_k]$
does not exceed $2Qh^{-2}$, it follows from the pinching assumption,
the gradient estimates in the canonical neighborhood assumption and
the evolution equation of the metric that the diameter of the set
$B_{T_0}(O,A_{k}h)$ at any time $t\in [T_0,t^+)$ is bounded from
above by $4\delta^{-\frac{1}{2}}h$. These imply that no point of the
ball $B_{T_0}(O,A_kh)$ at any time near $t^+$ can be the center of a
$\delta$-neck for any
$0<\delta<\widetilde{\delta}(A,\theta,\varepsilon)$ with
$\widetilde{\delta}(A,\theta,\varepsilon)>0$ small enough, since
$4\delta^{-\frac{1}{2}}h<< \delta^{-1}h$. However, the solution
disappears somewhere in the ball $B_{T_0}(O,A_kh)$ at the time $t^+$
because of a $\delta$-cutoff surgery and the surgery is always done
along the middle three-sphere of a $\delta$-neck. So the set
$B_{T_0}(O,A_kh)$ at the time $t^+$ is a part of a capped horn.
(Recall that we have declared extinct every compact component with
positive curvature operator or lying in $\Omega \setminus
\Omega_{\sigma}$). And then for each point of $B_{T_0}(O,A_kh)$ the
solution terminates at $t^+$.

The above observations will give us the following consequence.\vskip
0.2cm

$\underline{\mbox{\textbf{Claim 2.}}}$ \emph{For any
$\widetilde{L}<+\infty$, one can find $A=A(\widetilde{L})<+\infty$
and $\theta=\theta(\widetilde{L}),0<\theta<\frac{3}{2}$, with the
following property. Suppose $\gamma$ is a smooth curve in the set
$B_{T_0}(O,Ah)$, parametrized by $t\in [T_0,T_{\gamma}]$, such
that $\gamma(T_0)\in B_{T_0}(O,\frac{1}{2}Ah)$ and either
$T_{\gamma}=T_1$ and the solution on $B_{T_0}(O,Ah)$ exists up to
the time interval $[T_0,T_1]$ with $T_1=\min\{t_0,T_0+\theta
h^2\}<t_0$, or $T_{\gamma}<T_1$ and $\gamma(T_{\gamma})\in
\partial B_{T_0}(O,Ah)$. Then as
$\delta<\widetilde{\delta}(A,\theta,\varepsilon)$ chosen before,
there holds
$$\int^{T_{\gamma}}_{T_0}(R(\gamma(t),t)+|\dot{\gamma}(t)|^2_{g_{ij}(t)})dt >\widetilde{L}. \eqno (5.7)$$}

Indeed, we know from the estimate (5.6) that on the standard
solution,
\begin{eqnarray*}
\int_{0}^{\theta}R dt&\geq& const.
\int_{0}^{\theta}(\frac{3}{2}-t)^{-1}dt\\
&=&-const.\log(1-\frac{2\theta}{3}).
\end{eqnarray*}
By choosing $\theta=\theta(\widetilde{L})$ sufficiently close to
$\frac{3}{2}$, we have the desired estimate on the standard
solution.

If $T_{\gamma}=T_1<t_0$ and the solution on $B_{T_0}(O,Ah)$ exists
up to the time interval $[T_0,T_1]$, the solution in the parabolic
region $P(O,T_0,Ah,T_1-T_0)=\{(x,t)|x\in B_t(O,Ah), t\in
[T_0,T_1]\}$ is, after scaling with factor $h^{-2}$ and shifting
time $T_0$ to zero, $A^{-1}$-close to the corresponding subset in
the standard solution. Then we have
$$\int^{T_{\gamma}}_{T_0}(R(\gamma(t),t)+|\dot{\gamma}(t)|^2_{g_{ij}(t)})dt\geq const.
\int^{\theta}_0(\frac{3}{2}-t)^{-1}dt$$ $$=-const.
\log(1-\frac{2\theta}{3}),$$ which gives the desired estimate in
this case.

While if $T_{\gamma}<T_1$ and $\gamma(T_{\gamma})\in \partial
B_{T_0}(O,Ah)$, we see that the solution on $B_{T_0}(O,A_0h)$
exists up to the time interval $[T_0,T_{\gamma}]$ and is, after
scaling, $A^{-1}$-close to corresponding set in the standard
solution. Let $\theta=\theta(\widetilde{L})$ be chosen as above
and set $Q=Q(\widetilde{L})$ to be the maximum of the scalar
curvature of the standard solution in the time interval
$[0,\theta]$. On the standard solution, we can choose
$A=A(\widetilde{L})$ so large that for each $t\in [0,\theta]$,
 \begin{eqnarray*}d_t(O,\partial
B_0(O,A))&\geq& d_0(O,\partial B_0(O,A))-4(Q+1)t\\
&\geq& A-4(Q+1)\theta\\&\geq& \frac{3}{5}A,
\end{eqnarray*} and $$d_t(O,\partial B_0(O,\frac{A}{2}))\leq
\frac{A}{2},$$ where we used Lemma 8.3 of \cite{P1} in the first
inequality. Now our solution in the subset $B_{T_0}(O,Ah)$ up to
the time interval $[T_0,T_{\gamma}]$ is (after scaling)
$A^{-1}$-close to the corresponding subset in the standard
solution. This implies $$\frac{1}{5}Ah\leq
\int^{T_{\gamma}}_{T_0}|\dot{\gamma}(t)|_{g_{ij}(t)}\leq
(\int^{T_{\gamma}}_{T_0}|\dot{\gamma}(t)|^2_{g_{ij}(t)}dt)^{\frac{1}{2}}(T_{\gamma}-T_0)^{\frac{1}{2}}$$
and then
$$\int^{T_{\gamma}}_{T_0}(R(\gamma(t),t)+|\dot{\gamma}(t)|^2_{g_{ij}(t)})dt \geq
\frac{A^2}{25\theta}>\widetilde{L},$$ by choosing
$A=A(\widetilde{L})$ large enough. This proves the Claim 2.

We now use the above Claim 2 to verify Claim 1. Since $r_0\geq
\frac{1}{2\eta}r$ and $|R_m|\leq r^{-2}_0$ on
$P(x_0,t_0,r_0,-r^2_0)=\{(x,t)|x\in B_t(x_0,r_0),t\in
[t_0-r^2_0,t_0]\}$ (whenever it is defined), we can require
$\widetilde{\delta}>0$, depending on $r$ and $\widetilde{r}_{\ell}$,
so that $\gamma(T_0)$ does not lie in the region
$P(x_0,t_0,r_0,-r^2_0)$. Let $\triangle t$ be maximal such that
$\gamma|_{[t_0-\triangle t,t_0]}\subset P(x_0,t_0,r_0,-\triangle t)$
(i.e., $t=t_0-\triangle t$ is the first time for $\gamma$ escaping
the parabolic region $P(x_0,t_0,r_0,-r^2_0)$). Obviously we may
assume that
$$\int^{\triangle
t}_0\sqrt{\tau}(R(\gamma(t_0-\tau),t_0-\tau)+|\dot{\gamma}(t_0-\tau)|^2_{g_{ij}(t_0-\tau)})d\tau<L.$$
If $\triangle t<r^2_0$, it follows from the curvature bound
$|Rm|\leq r^{-2}_0$ on $P(x_0,t_0,r_0,-r^2_0)$ and the Ricci flow
equation that $$\int^{\triangle
t}_0|\dot{\gamma}(t_0-\tau)|d\tau\geq cr_0$$ for some universal
positive constant $c$. On the other hand, by Cauchy-Schwartz
inequality, we have
 \begin{eqnarray*}\int^{\triangle
t}_0|\dot{\gamma}(t_0-\tau)|d\tau &\leq& (\int^{\triangle
t}_0\sqrt{\tau}(R+|\dot{\gamma}|^2)d\tau)^{\frac{1}{2}}\cdot(\int^{\triangle
t}_0\frac{1}{\sqrt{\tau}}d\tau)^{\frac{1}{2}}\\ &\leq& 2
L^{\frac{1}{2}}(\triangle t)^{\frac{1}{4}}
\end{eqnarray*}
 which yields
$$(\triangle t)^{\frac{1}{2}}\geq \frac{c^2r^2_0}{4L}.$$ Thus we
always have $$(\triangle t)^{\frac{1}{2}}\geq
\min\{r_0,\frac{c^2r^2_0}{4L}\}.$$ Then
\begin{eqnarray*}\int^{t_0-T_0}_0\sqrt{\tau}(R+|\dot{\gamma}|^2)d\tau&\geq&
(\triangle t)^{\frac{1}{2}}\int^{t_0-T_0}_{\triangle
t}(R+|\dot{\gamma}|^2)d\tau\\ &\geq&
\min\{r_0,\frac{c^2r^2_0}{4L}\}\int^{t_0-T_0}_{\triangle
t}(R+|\dot{\gamma}|^2)d\tau. \end{eqnarray*} By applying Claim 2, we
can require the above $\widetilde{\delta}$ further to find
$\widetilde{\delta}=\widetilde{\delta}(L,r,\widetilde{r}_{\ell})>0$
so small that as $0<\delta<\widetilde{\delta}$, there holds
$$\int^{t_0-T_0}_{\triangle t}(R+|\dot{\gamma}|^2)d\tau\geq L(\min\{r_0,\frac{c^2r^2_0}{4L}\})^{-1}.$$ Hence we have verified the
desired assertion (5.5).

Now choose $L=100$ in (5.5), then it follows from Claim 1 that there
exists $\widetilde{\delta}>0$, depending on $r$ and
$\widetilde{r}_{\ell}$, such that as each $\delta$-cutoff surgery at
the time interval $t\in [(\ell-1)\varepsilon^2,T]$ has
$\delta<\widetilde{\delta}$, every barely admissible curve $\gamma$
with endpoints $(x_0,t_0)$ and $(x,t)$, where $t\in
[(\ell-1)\varepsilon^2,t_0)$, has
$$L(\gamma)=\int^{t_0-t}_0\sqrt{\tau}(R(\gamma(\tau),t_0-\tau)+|\dot{\gamma}(\tau)|^2_{g_{ij}(t_0-\tau)})d\tau\geq
100,$$ which implies the reduced distance from $(x_0,t_0)$ to
$(x,t)$ satisfies $$l\geq 25\varepsilon^{-1}. \eqno (5.8)$$

We also observe that the absolute value of $l(x_0,\tau)$ is very
small as $\tau$ closes to zero. We can then apply a maximum
principle argument as in Section 7.1 of \cite{P1} to conclude
$$l_{\min}(\tau)=\min\{l(x,\tau)|\mbox{ $x$ lies on the
solution manifold at time $t_0-\tau$}\}\leq 2,$$ for $\tau\in
(0,t_0-(\ell-1)\varepsilon^2],$ because barely admissible curves
do not carry minimum. In particular, there exists a minimizing
curve $\gamma$ of $l_{\min}(t_0-(\ell-1)\varepsilon^2)$, defined
on $\tau\in [0,t_0-(\ell-1)\varepsilon^2]$ with $\gamma(0)=x_0$,
such that
$$L(\gamma)\leq 2\cdot(2\sqrt{2}\varepsilon)<10\varepsilon. \eqno
(5.9)$$ Consequently, there exists a point
$(\overline{x},\overline{t})$ on the minimizing curve $\gamma$
with $\overline{t}\in
[(\ell-1)\varepsilon^2+\frac{1}{4}\varepsilon^2,(\ell-1)\varepsilon^2+\frac{3}{4}\varepsilon^2]$
such that $$R(\overline{x},\overline{t})\leq
50\widetilde{r}^{-2}_{\ell}. \eqno (5.10)$$ Otherwise, we would
have
\begin{eqnarray*}
L(\gamma)&\geq& \int^{t_0-(\ell-1)\varepsilon^2-\frac{1}{4}
\varepsilon^2}_{t_0-(\ell-1)\varepsilon^2-\frac{3}{4}\varepsilon^2}\sqrt{\tau}R(\gamma(\tau),t_0-\tau)d\tau
\\&\geq& 50 \widetilde{r}^{-2}_{\ell}\cdot
\frac{2}{3}(\frac{1}{2}\varepsilon^2)^{\frac{3}{2}}\\&>&
10\varepsilon
\end{eqnarray*} since $0<\widetilde{r}_{\ell}<\varepsilon$; this contradicts (5.9).

Next we want to get a lower bound for the reduced volume of a ball
around $\overline{x}$ of radius about $\widetilde{r}_{\ell}$ at
some time-slice slightly before $\overline{t}$. Since the solution
satisfies the canonical neighborhood assumption on the time
interval $[(\ell-1)\varepsilon^2,\ell \varepsilon^2)$, it follows
from the gradient estimate (5.4) that $$R(x,t)\leq
400\widetilde{r}^{-2}_{\ell} \eqno (5.11)$$ for those $(x,t)\in
P(\overline{x},\overline{t},\frac{1}{16}\eta^{-1}\widetilde{r}_{\ell},-\frac{1}{64}\eta^{-1}\widetilde{r}^2_{\ell})$
for which the solution is defined. And since the points where
occur the $\delta$-cutoff surgeries in the time interval
$[(\ell-1)\varepsilon^2,\ell \varepsilon^2)$ have their scalar
curvature at least $\delta^{-2}\widetilde{r}^{-2}_{\ell}$, the
solution is defined on the whole parabolic region
$P(\overline{x},\overline{t},\frac{1}{16}\eta^{-1}\widetilde{r}_{\ell},-\frac{1}{64}\eta^{-1}\widetilde{r}^2_{\ell})$
(this says, this parabolic region is unaffected by surgery). Thus
by combining (5.9) and (5.11), the reduced distance from
$(x_0,t_0)$ to each point of the ball
$B_{\overline{t}-\frac{1}{64}\eta^{-1}\widetilde{r}^2_{\ell}}(\overline{x},\frac{1}{16}\eta^{-1}\widetilde{r}_{\ell})$
is uniformly bounded by some universal constant. Let us define the
reduced volume of the ball
$B_{\overline{t}-\frac{1}{64}\eta^{-1}\widetilde{r}^2_{\ell}}(\overline{x},\frac{1}{16}\eta^{-1}\widetilde{r}_{\ell})$
by
$$\widetilde{V}_{t_0-\overline{t}+\frac{1}{64}\eta^{-1}\widetilde{r}^2_{\ell}}(B_{\overline{t}-\frac{1}{64}
\eta^{-1}\widetilde{r}^2_{\ell}}(\overline{x},\frac{1}{16}\eta^{-1}\widetilde{r}_{\ell}))$$
$$=\int_{B_{\overline{t}-\frac{1}{64}\eta^{-1}\widetilde{r}^2_{\ell}}
(\overline{x},\frac{1}{16}\eta^{-1}\widetilde{r}_{\ell})}(4\pi(t_0-\overline{t}+\frac{1}
{64}\eta^{-1}\widetilde{r}^2_{\ell}))^{-2}\exp(-l(q,t_0-\overline{t}+\frac{1}{64}\eta^{-1}\widetilde{r}^2_{\ell}))
dV_{\overline{t}-\frac{1}{64}\eta^{-1}\widetilde{r}^2_{\ell}}(q).$$
Hence by the $\kappa_{\ell}$-noncollapsing assumption on the time
interval $[(\ell-1)\varepsilon^2,\ell\varepsilon^2)$, we conclude
that the reduced volume of the ball
$B_{\overline{t}-\frac{1}{64}\eta^{-1}\widetilde{r}^2_{\ell}}(\overline{x},\frac{1}{16}\eta^{-1}\widetilde{r}_{\ell})$
is bounded from below by a positive constant depending only on
$\kappa_{\ell}$ and $\widetilde{r}_{\ell}$.

Finally we want to get a lower bound estimate for the volume of the
ball $B_{t_0}(x_0,r_0)$. We have seen the reduced distance from
$(x_0,t_0)$ to each point of the ball
$B_{\overline{t}-\frac{1}{64}\eta^{-1}\widetilde{r}^2_{\ell}}(\overline{x},\frac{1}{16}\eta^{-1}\widetilde{r}_{\ell})$
is uniformly bounded by some universal constant. Without loss of
generality, we may assume $\varepsilon>0$ is very small. Then it
follows from (5.8) that the points in the ball
$B_{\overline{t}-\frac{1}{64}\eta^{-1}\widetilde{r}^2_{\ell}}(\overline{x},\frac{1}{16}\eta^{-1}\widetilde{r}_{\ell})$
can be connected to $(x_0,t_0)$ by shortest $\mathcal{L}$-geodesics,
and all of these $\mathcal{L}$-geodesics are admissible (i.e., they
stay in the region unaffected by surgery). The union of all shortest
$\mathcal{L}$-geodesics from $(x_0,t_0)$ to the ball
$B_{\overline{t}-\frac{1}{64}\eta^{-1}\widetilde{r}^2_{\ell}}(\overline{x},\frac{1}{16}\eta^{-1}\widetilde{r}_{\ell})$,
denoted by
$CB_{\overline{t}-\frac{1}{64}\eta^{-1}\widetilde{r}^2_{\ell}}(\overline{x},\frac{1}{16}\eta^{-1}\widetilde{r}_{\ell})$,
forms a cone-like subset in space-time with the vertex $(x_0,t_0)$.
Denote $B(t)$ by the intersection of
$CB_{\overline{t}-\frac{1}{64}\eta^{-1}\widetilde{r}^2_{\ell}}(\overline{x},\frac{1}{16}\eta^{-1}\widetilde{r}_{\ell})$
with the time-slice at $t$. The reduced volume of the subset $B(t)$
is defined by
$$\widetilde{V}_{t_0-t}(B(t))=\int_{B(t)}(4\pi(t_0-t))^{-2}\exp(-l(q,t_0-t))dV_t(q).$$
Since the cone-like subset
$CB_{\overline{t}-\frac{1}{64}\eta^{-1}\widetilde{r}^2_{\ell}}(\overline{x},\frac{1}{16}\eta^{-1}\widetilde{r}_{\ell})$
lies entirely in the region unaffected by surgery, we can apply
Perelman's Jacobian comparison \cite{P1} to conclude that
\begin{equation}\tag{5.12}
\begin{split}\widetilde{V}_{t_0-t}(B(t))&\geq
\widetilde{V}_{t_0-\overline{t}+\frac{1}{64}\eta^{-1}
\widetilde{r}^2_{\ell}}(B_{\overline{t}-\frac{1}{64}\eta^{-1}\widetilde{r}^2_{\ell}}
(\overline{x},\frac{1}{16}\eta^{-1}\widetilde{r}_{\ell}))
\\&\geq c(\kappa_{\ell},\widetilde{r}_{\ell})
\end{split}
\end{equation} for all $t\in
[\overline{t}-\frac{1}{64}\eta^{-1}\widetilde{r}^2_{\ell},t_0]$,
where $c(\kappa_{\ell},\widetilde{r}_{\ell})$ is some positive
constant depending only on $\kappa_{\ell}$ and
$\widetilde{r}_{\ell}$.

Denote by
$\xi=r^{-1}_0V_0\ell_{t_0}(B_{t_0}(x_0,r_0))^{\frac{1}{4}}$. Our
purpose is to give a positive lower bound for $\xi$. Without loss
of generality, we may assume $\xi<\frac{1}{4}$, thus $0<\xi
r^2_0<t_0-\bar{t}+\frac{1}{64}\eta^{-1}\widetilde{r}^2_{\ell}$.
And denote by $\widetilde{B}(t_0-\xi r^2_0)$ the subset of the
points at the time-slice $\{t=t_0-\xi r^2_0\}$ where every point
can be connected to $(x_0,t_0)$ by an admissible shortest
$\mathcal{L}$-geodesic. Clearly $B(t_0-\xi r^2_0)\subset
\widetilde{B}(t_0-\xi r^2_0)$.

Since $r_0\geq \frac{1}{2\eta}r$ and
$\widetilde{\delta}=\widetilde{\delta}(r,\widetilde{r}_{\ell},\varepsilon)$
sufficiently small, the  region $P(x_0,t_0,r_0,-r^{2}_0)$ is
unaffected by surgery. Then by the exactly same argument as
deriving (3.24) in the proof of Theorem 3.5, we see that there
exists a universal positive constant $\xi_0$ such that as $0<\xi
\leq \xi_0$, there holds $$\mathcal{L}\exp_{\{|v|\leq
\frac{1}{4}\xi^{-\frac{1}{2}}\}}(\xi r^2_0)\subset
B_{t_0}(x_0,r_0). \eqno (5.13)$$ The reduced volume
$\widetilde{B}(t_0-\xi r^2_0)$ is given by
\begin{equation}\tag{5.14}
\begin{split}
&\ \ \widetilde{V}_{\xi r^2_0}(\widetilde{B}(t_0-\xi
r^2_0))\\&=\int_{\widetilde{B}(t_0-\xi r^2_0)}(4\pi \xi
r^2_0)^{-2}\exp(-l(q,\xi r^2_0))dV_{t_0-\xi r^2_0}(q)\\
&= \int_{\widetilde{B}(t_0-\xi r^2_0)\cap \mathcal{L}\exp_{\{|v|\leq
\frac{1}{4}\xi^{-\frac{1}{2}}\}}(\xi r^2_0)}(4\pi \xi
r^2_0)^{-2}\exp(-l(q,\xi r^2_0))dV_{t_0-\xi
r^2_0}(q) \\
&\ \ \ +\int_{\widetilde{B}(t_0-\xi r^2_0)\setminus
\mathcal{L}\exp_{\{|v|\leq\frac{1}{4}\xi^{-\frac{1}{2}}\}}(\xi
r^2_0)}(4\pi \xi r^2_0)^{-2}\exp(-l(q,\xi r^2_0))dV_{t_0-\xi
r^2_0}(q).
\end{split}
\end{equation} By (5.13), the first term on the RHS of
(5.14) can be estimated by
\begin{equation}\tag{5.15}
\begin{split}& \int_{\widetilde{B}(t_0-\xi r^2_0)\cap
\mathcal{L}\exp_{\{|v|\leq \frac{1}{4}\xi^{-\frac{1}{2}}\}}(\xi
r^2_0)}(4\pi \xi r^2_0)^{-2}\exp(-l(q,\xi r^2_0))dV_{t_0-\xi
r^2_0}(q)\\&\leq e^{4\xi}\int_{B_{t_0}(x_0,r_0)}(4\pi \xi
r^2_0)^{-2}\exp(-l)dV_{t_0}(q)\\&\leq e^{4\xi}(4\pi)^{-2}\xi^2.
\end{split}
\end{equation} And the second term on
the RHS of (5.14) can be estimated by
\begin{equation}\tag{5.16}
\begin{split}
&\ \ \int_{\widetilde{B}(t_0-\xi r^2_0)\setminus
\mathcal{L}\exp_{\{|v|\leq\frac{1}{4}\xi^{-\frac{1}{2}}\}}(\xi
r^2_0)}(4\pi \xi r^2_0)^{-2}\exp(-l(q,\xi r^2_0))dV_{t_0-\xi
r^2_0}(q)\\&\leq \int_{\{|v|>\frac{1}{4}\xi^{-\frac{1}{2}}\}}(4\pi
\tau)^{-2}\exp(-l)J(\tau)|_{\tau=0}dv\\&=
(4\pi)^{-2}\int_{\{|v|>\frac{1}{4}\xi^{-\frac{1}{2}}\}}\exp(-|v|^2)dv,
\end{split}
\end{equation}
 by using Perelman's Jacobian comparison theorem \cite{P1}
(as deriving (3.30) in the proof of Theorem 3.5). Hence the
combination of (5.12), (5.14), (5.15) and (5.16) bounds $\xi$ from
blow by a positive constant depending only on $\kappa_{\ell}$ and
$\widetilde{r}_{\ell}$.

Therefore we have completed the proof of the lemma.
$$\eqno \#$$ \vskip 0.3cm

Now we can prove the proposition.\vskip 0.2cm

$\underline{\mbox{\textbf{Proof of Proposition 5.4}}}$.

 The proof
of the proposition is by induction: having constructed our sequences
for $1\leq j\leq \ell$, we make one more step, defining
$\widetilde{r}_{\ell+1}$, $\kappa_{\ell+1}$,
$\widetilde{\delta}_{\ell+1}$, and redefining
$\widetilde{\delta}_\ell=\widetilde{\delta}_{\ell+1}$. In views of
the previous lemma, we only need to define $\widetilde{r}_{\ell+1}$
and $\widetilde{\delta}_{\ell+1}$.

In Theorem 4.1 we have obtained the canonical neighborhood structure
for smooth solutions. When adapting the arguments in the proof of
Theorem 4.1 to the present surgical solutions, we will encounter two
new difficulties. The first new difficulty is how to take a limit
for the surgerically modified solutions. The idea to overcome the
first difficulty consists of two parts. The first part, due to
Perelman \cite{P2}, is to choose $\widetilde{\delta}_{\ell}$ and
$\widetilde{\delta}_{\ell+1}$ small enough to push the surgical
regions to infinity in space. (This is the reason why we need to
redefine $\widetilde{\delta}_{\ell} = \widetilde{\delta}_{\ell+1}$.)
The second part is to show that solutions are smooth on some uniform
small time intervals (on compact subsets) so that we can apply
Hamilton's compactness theorem, since we only have curvature bounds;
otherwise Shi's interior derivative estimate may not be applicable.
In fact, the second part idea is more crucial. That is just
concerned with the question whether the surgery times accumulate or
not. Unfortunately, as written down in the third paragraph of
section 5.4 of \cite{P2}, the second part was not addressed. The
second new difficulty is that, when extending the limiting
surgically modified solution backward in time, it is possible to
meet the surgical regions in finite time. This also indicates the
surgery times may accumulate. The idea to overcome this difficulty
is somewhat similar to the above second part idea for the first
difficulty. We will use the canonical neighborhood charaterization
of the standard solution in Corollary A.2 in Appendix to exclude
this possibility.

We now start to prove the proposition by contradiction. Suppose for
sequence of positive numbers $r^{\alpha}$ and
$\widetilde{\delta}^{\alpha\beta}$, satisfying
$r^{\alpha}\rightarrow0$ as $\alpha\rightarrow\infty$ and
$\widetilde{\delta}^{\alpha\beta}\leq\frac{1}{\alpha\beta}(\rightarrow0)$,
there exist sequences of solutions $g^{\alpha\beta}_{ij}$ to the
Ricci flow with surgery, where each of them has only a finite number
of cutoff surgeries and has the given compact four-manifold as
initial datum, so that the following two assertions hold:

(i) each $\delta$-cutoff at a time
$t\in[(\ell-1)\varepsilon^2,(\ell+1)\varepsilon^2]$ satisfies
$\delta \leq \widetilde{\delta}^{\alpha\beta}$; and

(ii) the solutions satisfy the statement of the proposition on
$[0,\ell\varepsilon^2]$, but violate the canonical neighborhood
assumption (with accuracy $\varepsilon$) with $r=r^{\alpha}$ on
$[\ell\varepsilon^2,(\ell+1)\varepsilon^2]$.

For each solution $g^{\alpha\beta}_{ij}$, we choose $\bar{t}$
(depending on $\alpha, \beta$) to be the nearly first time for which
the canonical neighborhood assumption (with accuracy $\varepsilon$)
is violated. More precisely, we choose
$\bar{t}\in[\ell\varepsilon^2,(\ell+1)\varepsilon^2]$ so that the
canonical neighborhood assumption with $r=r^{\alpha}$ and with
accuracy parameter $\varepsilon$ is violated at some
$(\bar{x},\bar{t})$, however the canonical neighborhood assumption
with accuracy parameter $2\varepsilon$ holds on
$t\in[\ell\varepsilon^2,\bar{t}]$. After passing to subsequences, we
may assume each $\widetilde{\delta}^{\alpha\beta}$ is less than the
$\widetilde{\delta}$ in Lemma 5.5 with $r=r^{\alpha}$ when $\alpha$
is fixed. Then by Lemma 5.5 we have uniform $\kappa$-noncollapsing
for all scales less than $\varepsilon$ on $[0,\bar{t}]$ with some
$\kappa>0$ independent of $\alpha,\beta$.

Slightly abusing notation, we will often drop the indices
$\alpha,\beta$.

Let $\widetilde{g}^{\alpha\beta}_{ij}$ be the rescaled solutions
along $(\bar{x},\bar{t})$ with factors $R(\bar{x},\bar{t})(\geq
r^{-2}\rightarrow+\infty)$ and shift the times $\bar{t}$ to zero.
We hope to take a limit of the rescaled solutions for subsequences
of $\alpha,\beta\rightarrow\infty$ and show the limit is an
ancient $\kappa$-solution, which will give the desired
contradiction. We divide the following arguments into six steps.

\vskip 0.2cm \noindent{\bf Step 1.}  Let $(y,\hat{t})$ be a point on
the rescaled solution $\widetilde{g}^{\alpha\beta}_{ij}$ with
$\widetilde{R}(y,\hat{t})\leq A$ $(A\geq 1)$ and $\hat{t}\in
[-(\bar{t} - (\ell-1)\varepsilon^{2})R(\bar{x},\bar{t}),0]$, then we
have estimate
$$\widetilde{R}(x,t)\leq10A \eqno (5.17)$$ for those
$(x,t)$ in the parabolic neighborhood
$P(y,\hat{t},\frac{1}{2}\eta^{-1}A^{-\frac{1}{2}},-\frac{1}{8}\eta^{-1}A^{-1})$
$\triangleq\{(x',t')\ |\ x'\in
\widetilde{B}_{t'}(y,\frac{1}{2}\eta^{-1}A^{-\frac{1}{2}}),t'\in[\hat{t}-\frac{1}{8}\eta^{-1}A^{-1},\hat{t}]\}$,
for which the rescaled solution is defined.

Indeed, as in the first step of the proof of Theorem 4.1, this
follows directly from the gradient estimates (5.4) in the
canonical neighborhood assumption with parameter $2\varepsilon$.

\vskip 0.2cm \noindent{\bf Step 2.} In this step, we will prove
three time extending results.\vskip 0.2cm

\noindent \textbf{Assertion 1.}\ \
    \emph{ For arbitrarily fixed $\alpha$, $0<A<+\infty$,  $1 \leq C<+\infty$
     and $0 \leq B < \frac{1}{2}\varepsilon^2(r^{\alpha})^{-2}-
\frac{1}{8}\eta^{-1}C^{-1}$,
     there is a $\beta_0=\beta_0(\varepsilon, A, B, C)$
     (independent of $\alpha$)
such that if $\beta\geq\beta_0$ and the rescaled solution
$\widetilde{g}^{\alpha\beta}_{ij}$ on the ball
$\widetilde{B}_{0}(\bar{x},A)$ is defined on a time interval
$[-b,0]$ with $0 \leq b \leq B$ and the scalar curvature satisfies
$$\widetilde{R}(x,t)\leq C , \ \ \mbox{ on }
\widetilde{B}_{0}(\bar{x},A) \times [-b,0],$$ then the rescaled
solution $\widetilde{g}^{\alpha\beta}_{ij}$ on the ball
$\widetilde{B}_{0}(\bar{x},A)$ is also defined on the extended time
interval $ [-b-\frac{1}{8}\eta^{-1}C^{-1},0]$.} \vskip 0.2cm

Before the proof, we need a simple \textbf{observation}: once a
space point in the Ricci flow with surgery is removed by surgery at
some time, then it never appears for later time; if a space point at
some time $t$ can not be defined before the time $t$ , then either
the point lies in a gluing cap of the surgery at time $t$ or the
time $t$ is the initial time of the Ricci flow.

\vskip 0.2cm \noindent \textbf{Proof of Assertion 1.}\ \ Firstly we
claim that there exists $\beta_0=\beta_0(  \varepsilon, A, B, C)$
such that as $\beta\geq\beta_0$, the rescaled solution
$\widetilde{g}^{\alpha\beta}_{ij}$ on the ball
$\widetilde{B}_{0}(\bar{x},A)$ can be defined before the time $-b$
(i.e., there are no surgeries interfering in
$\widetilde{B}_{0}(\bar{x},A)\times [-b-\epsilon',-b]$ for some
$\epsilon'>0$).

 We argue by contradiction. Suppose not, then there is some
point $\tilde{x}\in \widetilde{B}_{0}(\bar{x},A)$ such that the
rescaled solution $\widetilde{g}^{\alpha\beta}_{ij}$ at $\tilde{x}$
can not be defined before the time $-b$. By the above observation,
there is a surgery at the time $-b$ such that the point $\tilde{x}$
lies in the instant gluing cap.

 Let $\tilde{h}$  $(=R(\bar{x},\bar{t})^{\frac{1}{2}}h$)
be the cut-off radius at the time $-b$ for the rescaled solution.
Clearly, there is a universal constant $D$ such that
$${D}^{-1}\tilde{h}\leq\widetilde{R}(\tilde{x},-b)^{-\frac{1}{2}}\leq
D\tilde{h}.$$

 By Lemma 5.3 and looking at the rescaled solution at the time $-b$, the gluing cap and the adjacent
$\delta$-neck, of radius $\tilde{h}$, constitute a
${(\widetilde{\delta}^{\alpha\beta})}^{\frac{1}{2}}$-cap
$\mathcal{K}$. For any fixed small positive constant
$\delta^{\prime}$ (much smaller than $\varepsilon$), we see
$$\widetilde{B}_{(-b)}(\tilde{x},{(\delta^{\prime})}^{-1}\widetilde{R}(\tilde{x},
-b)^{-\frac{1}{2}})\subset\mathcal{K}$$ as $\beta$ large enough.
We first verify the following

\vskip 0.2cm \noindent \textbf{Claim 1.}\emph{ For any small
constants $0<\tilde{\theta}<\frac{3}{2}$, $\delta'>0$, there exists
a $\beta(\delta',\varepsilon, \tilde{\theta})> 0$ such that as
$\beta \geq \beta(\delta',\varepsilon, \tilde{\theta})$, we have}

(i) \emph{the rescaled solution $\widetilde{g}^{\alpha\beta}_{ij}$
over
$\widetilde{B}_{(-b)}(\tilde{x},{(\delta^{\prime})}^{-1}\tilde{h})$
 is defined on the time interval $[-b,0]\cap[-b,-b+(\frac{3}{2}-\tilde{\theta})\tilde{h}^{2}]$;}

 (ii)  \emph{the ball
$\widetilde{B}_{(-b)}(\tilde{x},{(\delta^{\prime})}^{-1}\tilde{h})$
in the ${(\widetilde{\delta}^{\alpha\beta})}^{\frac{1}{2}}$-cap
$\mathcal{K}$ evolved by the Ricci flow  on the time interval
$[-b,0]\cap[-b,-b+(\frac{3}{2}-\tilde{\theta})\tilde{h}^{2}]$ is,
after scaling with factor $\tilde{h}^{-2}$,
${\delta}^{\prime}$-close ( in $C^{[\delta'^{-1}]}$ topology) to the
corresponding subset of the standard solution.}\vskip 0.2cm

This claim is somewhat known in the first claim in the proof of
Lemma 5.5. Indeed, suppose there is a surgery at some time
$\tilde{\tilde{t}}\in
[-b,0]\cap(-b,-b+(\frac{3}{2}-\tilde{\theta})\tilde{h}^{2}]$ which
removes some point  $\tilde{\tilde{x}}\in \widetilde
B_{(-b)}(\tilde{x},{(\delta^{\prime})}^{-1}\tilde{h})$. We assume
$\tilde{\tilde{t}}\in (-b,0]$ be the first time with that
property.

Then by the proof of the first claim in Lemma 5.5, there is a
$\bar{\delta}=\bar{\delta}(\delta',\varepsilon,\tilde{\theta})$ such
that if $\widetilde{\delta}^{\alpha\beta}<\bar{\delta}$, then the
ball
$\widetilde{B}_{(-b)}(\tilde{x},{(\delta^{\prime})}^{-1}\tilde{h})$
in the ${(\widetilde{\delta}^{\alpha\beta})}^{\frac{1}{2}}$-cap
$\mathcal{K}$ evolved by the Ricci flow  on the time interval
$[-b,\tilde{\tilde{t}})$ is, after scaling with factor
$\tilde{h}^{-2}$, ${\delta}^{\prime}$-close to the corresponding
subset of the standard solution. Note that the metrics for times in
$[-b,\tilde{\tilde{t}})$ on $\widetilde
B_{(-b)}(\tilde{x},{(\delta^{\prime})}^{-1}\tilde{h})$ are
equivalent. By the proof of the first claim in Lemma 5.5, the
solution on $\widetilde
B_{(-b)}(\tilde{x},{(\delta^{\prime})}^{-1}\tilde{h})$ keeps looking
like a cap for $t\in [-b,\tilde{\tilde{t}})$. On the other hand, by
definition, the surgery is always performed along the middle
three-sphere of a $\delta$-neck with $\delta <
\widetilde{\delta}^{\alpha\beta}$. Then as $\beta$ large, all the
points in $\widetilde
B_{(-b)}(\tilde{x},{(\delta^{\prime})}^{-1}\tilde{h})$ are removed
(as a part of a capped horn) at the time $\tilde{\tilde{t}}$. But
$\tilde{x}$ (near the tip of the cap) exists past the time
$\tilde{\tilde{t}}$.  This is a contradiction. Hence we have proved
that
$\widetilde{B}_{(-b)}(\tilde{x},{(\delta^{\prime}})^{-1}\tilde{h})$
 is defined on the time interval
 $[-b,0]\cap[-b,-b+(\frac{3}{2}-\tilde{\theta})\tilde{h}^{2}]$.

 The $\delta'$-closeness of the solution on $\widetilde{B}_{(-b)}(\tilde{x},
 {(\delta^{\prime})}^{-1}h)\times ([-b,0]\cap[-b,-b+(\frac{3}{2}-\tilde{\theta})\tilde{h}^{2}])$
with the corresponding subset of the standard solution follows by
the uniqueness theorem and the canonical neighborhood assumption
with parameter $2\varepsilon$ as in the proof of the first claim
in Lemma 5.5. Then we have proved Claim 1.

\vskip 0.2cm
 We next verify the following

 \vskip 0.2cm \noindent \textbf{Claim 2.} \emph{ There is $\tilde{\theta}=\tilde{\theta}(CB)$,
 $0<\tilde{\theta}<\frac{3}{2}$,
  such that $b\leq
(\frac{3}{2}-\tilde{\theta})\tilde{h}^{2}$ as $\beta$ large.} \vskip
0.2cm

Note from Theorem A.1 in Appendix, there is a universal constant
$D'>0$ such that the standard solution ( of dimension four)
satisfies the following curvature estimate
$$R(y,s) \geq \frac{2D'}{\frac{3}{2}-s}.$$
We choose $\tilde{\theta}= {3D'}/{2(2D'+2CB)}$. Then as $\beta$
large enough, the rescaled solution satisfies
$$
\widetilde{R}(x,t)\geq
\frac{D'}{\frac{3}{2}-(t+b)\tilde{h}^{-2}}\tilde{h}^{-2} \eqno
(5.18)
$$
on
$\widetilde{B}_{(-b)}(\tilde{x},{(\delta^{\prime})}^{-1}\tilde{h})\times
([-b,0]\cap[-b,-b+(\frac{3}{2}-\tilde{\theta})\tilde{h}^{2}])$.

 Suppose $b\geq
(\frac{3}{2}-\tilde{\theta})\tilde{h}^{2}$. Then by combining with
the assumption $\widetilde{R}(\tilde{x},t)\leq {C}$ for
$t=(\frac{3}{2}-\tilde{\theta})\tilde{h}^{2}-b$, we have
$$C\geq \frac{D'}{\frac{3}{2}-(t+b)\tilde{h}^{-2}}\tilde{h}^{-2},$$
and then
$$ \tilde{\theta}\geq\frac{\frac{3D'}{2CB}}{1+\frac{D'}{CB}}.
$$ This is a
contradiction. Hence we have proved Claim 2. \vskip 0.2cm

The combination of the above two claims shows that there is a
positive constant $0 <\tilde{\theta}=\tilde{\theta}(CB)<\frac{3}{2}$
such that for any $\delta'>0$, there is a positive
$\beta(\delta',\varepsilon, \tilde{\theta})$ such that as $\beta
\geq \beta(\delta',\varepsilon, \tilde{\theta})$, we have $b\leq
(\frac{3}{2}-\tilde{\theta})\tilde{h}^{2}$ and the rescaled solution
in the ball
$\widetilde{B}_{(-b)}(\tilde{x},{(\delta^{\prime}})^{-1}\tilde{h})$
on the time interval $[-b,0]$ is, after scaling with factor
$\tilde{h}^{-2}$, ${\delta}^{\prime}$-close ( in
$C^{[(\delta')^{-1}]}$ topology) to the corresponding subset of the
standard solution.

By (5.18) and the assumption $\widetilde{R}\leq {C}$ on
$\widetilde{B}_{0}(\bar{x},A) \times [-b,0],$ we know that the
cut-off radius $\tilde{h}$ at the time $-b$ for the rescaled
solution satisfies
$$\tilde{h} \geq \sqrt{\frac{2D'}{3C}}.$$

Let $\delta'>0$ be much smaller than $\varepsilon$ and
$\min\{A^{-1},A\}$. Since $\tilde{d}_0(\tilde{x},\bar{x})\leq A$,
it follows that there is constant $C(\tilde{\theta})$ depending
only on $\tilde{\theta}$ such that
$\tilde{d}_{(-b)}(\tilde{x},\bar{x})\leq C(\tilde{\theta})A \ll
(\delta')^{-1}\tilde{h}$. We now apply Corollary A.2 in Appendix
with the accuracy parameter ${\varepsilon}/{2}$. Let
$C({\varepsilon}/{2})$
 be the positive constant in Corollary A.2. Without loss of
 generality, we may assume the positive constant
 $C_1(\varepsilon)$ in the canonical neighborhood assumption is larger than $4C({\varepsilon}/{2})$.
 As $\delta'>0$ is much smaller that $\varepsilon$ and $\min\{A^{-1},A\}$,
 the point $\bar{x}$ at the time $\bar{t}$ has a neighborhood which is either a
 $\frac{3}{4}\varepsilon$-cap or a $\frac{3}{4}\varepsilon$-neck.

 Since the
canonical neighborhood assumption with accuracy parameter
$\varepsilon$ is violated at $(\bar{x},\bar{t})$, the neighborhood
of the point $\bar{x}$ at the new time zero for the rescaled
solution must be a $\frac{3}{4}\varepsilon$-neck. By Corollary A.2
(b), we know the neighborhood is the slice at the time zero of the
parabolic neighborhood
$$P(\bar{x},0,\frac{4}{3}\varepsilon^{-1}\widetilde{R}(\bar{x},0)^{-\frac{1}{2}},-\min
\{\widetilde{R}(\bar{x},0)^{-1},b\})$$ (with
$\widetilde{R}(\bar{x},0)=1$) which is
$\frac{3}{4}\varepsilon$-close (in
$C^{[\frac{4}{3}\varepsilon^{-1}]}$ topology) to the corresponding
subset of the evolving standard cylinder $\mathbb{S}^3 \times
\mathbb{R}$ over the time interval $[-\min \{b,1\},0]$ with scalar
curvature $1$ at the time zero. If $b \geq 1$, the
$\frac{3}{4}\varepsilon$-neck is strong, which is a contradiction.
While if $b<1$, the $\frac{3}{4}\varepsilon$-neck at time $-b$ is
contained in the union of the gluing cap and the adjacent
 $\delta$-neck where the $\delta$-cutoff surgery was taken. Since
 $\varepsilon$ is small (say $\varepsilon< 1/100$), it is clear that the point $\bar{x}$ at time $-b$ is the center of an $\varepsilon$-neck
 which is entirely
 contained in the adjacent
 $\delta$-neck. By the remark after Lemma 5.2, the adjacent $\delta$-neck
 approximates an ancient $\kappa$-solution. This implies the point
 $\bar{x}$ at the time $\bar{t}$ has a strong $\varepsilon$-neck,
 which is also a contradiction.

 Hence we have proved that there exists $\beta_0=\beta_0(  \varepsilon,
A, B, C)$ such that as $\beta\geq\beta_0$, the rescaled solution
on the ball $\widetilde{B}_{0}(\bar{x},A)$ can be defined before
the time $-b$.

Let $[t_{A}^{\alpha\beta},0]\supset[-b,0]$ be the largest time
interval so that the rescaled solution
$\widetilde{g}^{\alpha\beta}_{ij}$ can be defined on
$\widetilde{B}_{0}(\bar{x},A)\times[t_{A}^{\alpha\beta},0]$. We
finally claim that $t_{A}^{\alpha\beta}\le
-b-\frac{1}{8}\eta^{-1}C^{-1}$ as $\beta$ large enough.

 Indeed, suppose not, by the gradient estimates as in Step 1,
we have the curvature estimate $$\widetilde{R}(x,t)\leq {10C}
$$ on $\widetilde{B}_{0}(\bar{x},A)\times
[t_{A}^{\alpha\beta},-b]$. Hence we have the curvature estimate
$$\widetilde{R}(x,t)\leq {10C} $$ on
$\widetilde{B}_{0}(\bar{x},A)\times [t_{A}^{\alpha\beta},0]$. By
the above argument there is a $\beta_0=\beta_0( \varepsilon, A, B
+ \frac{1}{8}\eta^{-1}C^{-1}, 10C)$ such that as
$\beta\geq\beta_0$, the solution in the ball
$\widetilde{B}_{0}(\bar{x},A)$ can be defined before the time
$t_{A}^{\alpha\beta}$. This is a contradiction.

Therefore we have proved Assertion 1.\vskip 0.2cm

 \noindent \textbf{Assertion 2.}\ \
     \emph{For arbitrarily
     fixed $\alpha$,  $0<A<+\infty$, $1 \leq C<+\infty$  and $0 < B
<\frac{1}{2}\varepsilon^2 (r^{\alpha})^{-2}-
\frac{1}{50}\eta^{-1}$,
     there is a $\beta_0=\beta_0(\varepsilon, A,B,C)$ (independent of $\alpha$)
such that if $\beta\geq\beta_0$ and the rescaled solution
$\widetilde{g}^{\alpha\beta}_{ij}$ on the ball
$\widetilde{B}_{0}(\bar{x},A)$
 is defined on a time interval
$[-b+\epsilon',0]$ with $0 < b \leq B$ and $0<\epsilon'
<\frac{1}{50}\eta^{-1}$ and the scalar curvature satisfies
$$\widetilde{R}(x,t)\leq {C} \ \ \mbox{ on } \
\widetilde{B}_{0}(\bar{x},A)\times [-b+\epsilon',0],$$
 and there is a point $y\in \widetilde{B}_{0}(\bar{x},A)$ such
that $\widetilde{R}(y,-b+\epsilon')\leq \frac{3}{2}$, then the
rescaled solution $\widetilde{g}^{\alpha\beta}_{ij}$ at $y$ is
also defined on the extended time interval
$[-b-\frac{1}{50}\eta^{-1},0]$ and satisfies the estimate
$$\widetilde{R}(y,t)\leq 15 $$ for $t\in
[-b-\frac{1}{50}\eta^{-1},-b+\epsilon']$.} \vskip 0.2cm
 \noindent \textbf{Proof of Assertion 2.}  We imitate the proof of Assertion 1. If the
rescaled solution $\widetilde{g}^{\alpha\beta}_{ij}$ at $y$ can not
be defined for some time in
$[-b-\frac{1}{50}\eta^{-1},-b+\epsilon')$, then there is a surgery
at some time $\tilde{\tilde{t}} \in
[-b-\frac{1}{50}\eta^{-1},-b+\epsilon']$ such that $y$ lies in the
instant gluing cap.  Let $\tilde{h}$
$(=R(\bar{x},\bar{t})^{\frac{1}{2}}h$) be the cutoff radius at the
time $\tilde{\tilde{t}}$ for the rescaled solution. Clearly, there
is a universal constant $D>1$ such that
${D}^{-1}\tilde{h}\leq\widetilde{R}(y,\tilde{\tilde{t}})^{-\frac{1}{2}}\leq
D\tilde{h}$. By the gradient estimates as in Step 1, the cutoff
radius satisfies
$$\tilde{h} \geq D^{-1}15^{-\frac{1}{2}}.$$

As in Claim 1 (i) in the proof of Assertion 1, for any small
constants $0<\tilde{\theta}<\frac{3}{2}$, $\delta'>0$, there
exists a $\beta(\delta',\varepsilon, \tilde{\theta})> 0$ such that
as $\beta \geq \beta(\delta',\varepsilon, \tilde{\theta})$, there
is
 no surgery interfering in
$\widetilde{B}_{\tilde{\tilde{t}}}(y,(\delta')^{-1}\tilde{h})\times
([\tilde{\tilde{t}},(\frac{3}{2}-\tilde{\theta})\tilde{h}^{2}
+\tilde{\tilde{t}}]\cap(\tilde{\tilde{t}},0])$. Without loss of
generality, we may assume that the universal constant $\eta$ is much
larger than $D$. Then we have
$(\frac{3}{2}-\tilde{\theta})\tilde{h}^{2}
+\tilde{\tilde{t}}>-b+\frac{1}{50}\eta^{-1}$. As in Claim 2, we can
use the curvature bound assumption to choose
$\tilde{\theta}=\tilde{\theta}(B,C)$ such that
$(\frac{3}{2}-\tilde{\theta})\tilde{h}^{2} +\tilde{\tilde{t}}\geq
0$; otherwise
$$
C\geq\frac{D'}{\tilde{\theta}\tilde{h}^{2}}
$$
for some universal constant $D'$, and  $$|\tilde{\tilde{t}}+b|\leq
\frac{1}{50}\eta^{-1},$$ which implies $$
\tilde{\theta}\geq\frac{\frac{3D'}{2C(B+\frac{1}{50}\eta^{-1})}}{1+\frac{D'}{C(B+\frac{1}{50}\eta^{-1})}}.
$$ This is a contradiction
if we choose
$\tilde{\theta}={3D'}/{2(2D'+2C(B+\frac{1}{50}\eta^{-1}))}$.

So there is a positive constant $0
<\tilde{\theta}=\tilde{\theta}(B,C)<\frac{3}{2}$ such that for any
$\delta'>0$, there is a positive $\beta(\delta',\varepsilon,
\tilde{\theta})$ such that as $\beta \geq
\beta(\delta',\varepsilon, \tilde{\theta})$, we have
$-\tilde{\tilde{t}}\leq (\frac{3}{2}-\tilde{\theta})\tilde{h}^{2}$
and the solution in the ball
$\widetilde{B}_{\tilde{\tilde{t}}}(\tilde{x},{(\delta^{\prime})}^{-1}\tilde{h})$
on the time interval $[\tilde{\tilde{t}},0]$ is, after scaling
with factor $\tilde{h}^{-2}$, ${\delta}^{\prime}$-close (in
$C^{[\delta'^{-1}]}$ topology) to the corresponding subset of the
standard solution.

Then exactly as in the proof of Assertion 1, by using the canonical
neighborhood structure of the standard solution in Corollary A.2,
this gives the desired contradiction with the hypothesis that the
canonical neighborhood assumption with accuracy parameter
$\varepsilon$ is violated at $(\bar{x},\bar{t})$, as $\beta$
sufficiently large.

The curvature estimate at the point $y$ follows from Step 1.
Therefore we complete the proof of Assertion 2.\vskip 0.2cm

Note that the standard solution satisfies $R(x_1,t)\leq D''
R(x_2,t)$ for any $t\in [0,\frac{1}{2}]$ and any two points
$x_1,x_2$, where $D''\geq 1$ is a universal constant.
 \vskip 0.2cm \noindent \textbf{Assertion 3.}\
     \emph{For arbitrarily
     fixed $\alpha$,  $0<A<+\infty$, $1 \leq C<+\infty$ ,
     there is a $\beta_0=\beta_0(\varepsilon, AC^{\frac{1}{2}})$
such that if any point $(y_0,t_0)$  with  $0 \leq
     -t_0
<\frac{1}{2}\varepsilon^2 (r^{\alpha})^{-2}-
\frac{1}{8}\eta^{-1}C^{-1}$ of the rescaled solution
$\widetilde{g}^{\alpha\beta}_{ij}$ for $\beta\geq\beta_0$ satisfies
$\widetilde{R}(y_0,t_0)\leq C$ , then either the rescaled solution
at $y_0$ can be defined at least on
$[t_0-\frac{1}{16}\eta^{-1}C^{-1},t_0]$ and the rescaled scalar
curvature satisfies
$$\widetilde{R}(y_0,t)\leq 10 {C} \ \ \mbox{ for } t\in
[t_0-\frac{1}{16}\eta^{-1}C^{-1},t_0],
$$
or we have $$\widetilde{R}(x_1,t_0)\leq 2D''
\widetilde{R}(x_2,t_0)$$ for any two points $x_1, x_2\in
\widetilde{B}_{t_0}(y_0,A)$, where $D''$ is the above universal
constant.} \vskip 0.2cm

\noindent \textbf{Proof of Assertion 3.} Suppose the rescaled
solution $\widetilde{g}^{\alpha\beta}_{ij}$ at $y_0$ can not be
defined for some $t\in [t_0-\frac{1}{16}\eta^{-1}C^{-1},t_0)$, then
there is a surgery at some time $\tilde{t} \in
[t_0-\frac{1}{16}\eta^{-1}C^{-1},t_0]$ such that  $y_0 $ lies in the
instant gluing cap. Let $\tilde{h}$
$(=R(\bar{x},\bar{t})^{\frac{1}{2}}h$) be the cutoff radius at the
time $\tilde{t}$ for the rescaled solution
$\widetilde{g}^{\alpha\beta}_{ij}$.  By the gradient estimates as in
Step 1, the cutoff radius satisfies
$$\tilde{h} \geq D^{-1}10^{-\frac{1}{2}}C^{-\frac{1}{2}},$$
where $D$ is the universal constant in the proof of the Assertion
1.
 Since we assume $\eta$ is
suitable larger than $D$ as before, we have
$\frac{1}{2}\tilde{h}^{2} +\tilde{t}>t_0$. As in Claim 1 (ii) in
Assertion 1, for arbitrarily small $\delta'>0$, we know that as
$\beta$ large enough the rescaled solution on the ball
$\widetilde{B}_{\tilde{t}}(y_0,{(\delta^{\prime})}^{-1}\tilde{h})$
on the time interval $[\tilde{t},t_0]$ is, after scaling with
factor $\tilde{h}^{-2}$, ${\delta}^{\prime}$-close (in
$C^{[(\delta')^{-1}]}$ topology) to the corresponding subset of
the standard solution. Since  $(\delta')^{-1}\tilde{h}\gg A$ as
$\beta$ large enough, Assertion 3 follows from the curvature
estimate of standard solution in the time interval
 $ [0,\frac{1}{2}]$.

 \vskip 0.2cm \noindent{\bf Step 3.}
For any subsequence $(\alpha_m,\beta_m)$ of $(\alpha,\beta)$ with
$r^{\alpha_m}\rightarrow 0$ and $\delta^{\alpha_m\beta_m}\rightarrow
0$ as $m\rightarrow \infty$, we next argue as in the second step of
the proof of Theorem 4.1 to show that the curvatures of the rescaled
solutions $\tilde{g}^{\alpha_m\beta_m}$ at new times zero (after
shifting) stay uniformly bounded at bounded distances from $\bar{x}$
for all sufficiently large $m$. More precisely, we will prove the
following assertion:

 \vskip 0.1cm \noindent \textbf{Assertion 4.}\ \emph{Given any subsequence of the rescaled
solutions $\tilde{g}^{\alpha_m\beta_m}_{ij}$ with
$r^{\alpha_m}\rightarrow 0$ and $\delta^{\alpha_m\beta_m}\rightarrow
0$ as $m\rightarrow \infty$, then for any $L>0$, there are constants
$C(L)>0$  and $m(L)$ such that the rescaled solutions
$\tilde{g}^{\alpha_m\beta_m}_{ij}$ satisfy}

 (i) \emph{
$\tilde{R}(x,0)\leq C(L)$ for all points $x$ with
$\tilde{d}_{0}(x,\bar{x})\leq L$ and all $m \geq 1$;}

 (ii) \emph{ the rescaled
solutions over the ball $\tilde{B}_{0}(\bar{x},L)$ are defined  at
least on the time interval $[-\frac{1}{16}\eta^{-1}C(L)^{-1},0]$ for
all $m\geq m(L)$.} \vskip 0.2cm

\noindent \textbf{Proof of Assertion 4.} For all $\rho>0$, set
$$M(\rho)=\sup\{\tilde{R}(x,0)\ |\ m\geq 1 \mbox{ and } \tilde{d}_0(x,\bar{x})\leq\rho \ \mbox{ in the rescaled
solutions
 } \ \tilde{g}^{\alpha_m\beta_m}_{ij}\}$$ and

$$\rho_0=\sup\{\rho > 0\ |\ M(\rho)<+\infty\}.$$ Note that the
estimate (5.17) implies that $\rho_0>0$. For (i), it suffices to
prove $\rho_0=+\infty$.

We argue by contradiction. Suppose $\rho_0<+\infty$. Then there
are a sequence of points $y$ in the rescaled solutions
$\tilde{g}^{\alpha_m\beta_m}_{ij}$ with
$\tilde{d}_0(\bar{x},y)\rightarrow\rho_0<+\infty$ and
$\tilde{R}(y,0)\rightarrow +\infty$. Denote by $\gamma$ a
minimizing geodesic segment from $\bar{x}$ to $y$ and denote by
$\tilde{B}_0(\bar{x},\rho_0)$ the geodesic open ball centered at
$\bar{x}$ of radius $\rho_0$ on the rescaled solution
$\tilde{g}^{\alpha_m\beta_m}_{ij}$.

 First, we
claim that for any $0<\rho<\rho_0$ with $\rho$ near $\rho_0$,  the
rescaled solutions on the balls $\tilde{B}_{0}(\bar{x},\rho)$ are
defined on the time interval
$[-\frac{1}{16}\eta^{-1}M(\rho)^{-1},0]$ for all large $m$.
Indeed, this follows from Assertion 3 or Assertion 1. For the
later purpose in Step 6, we now present an argument by using
Assertion 3.  If the claim is not true, then there is a surgery at
some time $\tilde{t} \in [-\frac{1}{16}\eta^{-1}M(\rho)^{-1},0]$
such that some point $\tilde{y}\in \tilde{B}_{0}(\bar{x},\rho) $
lies in the instant gluing cap.
 We can choose sufficiently small
$\delta'>0$ such that $2\rho_0<(\delta')^{-\frac{1}{2}}\tilde{h}$,
where $\tilde{h} \geq
D^{-1}20^{-\frac{1}{2}}M(\rho)^{-\frac{1}{2}}$ are the cutoff
radius of the rescaled solutions at $\tilde{t}$. By applying
Assertion 3 with $(\tilde{y},0)=(y_0,t_0)$, we see that there is a
$m(\rho_0,M(\rho))>0$ such that as  $m\geq m(\rho_0,M(\rho))$,
$$\widetilde{R}(x,0)\leq 2D''
$$ for all $x\in \widetilde{B}_{0}(\bar{x},\rho)$. This is a
contradiction as $\rho\rightarrow\rho_0$.

Since for each fixed $0<\rho<\rho_0$ with $\rho$ near $\rho_0$, the
rescaled solutions on the ball $\tilde{B}_{0}(\bar{x},\rho)$ are
defined on the time interval
$[-\frac{1}{16}\eta^{-1}M(\rho)^{-1},0]$ for all large $m$, by Step
1 and Shi's derivative estimate, we know that the covariant
derivatives and higher order derivatives of the curvatures on
$\tilde{B}_{0}(\bar{x},\rho - \frac{(\rho_0
-\rho)}{2})\times[-\frac{1}{32}\eta^{-1}M(\rho)^{-1},0]$ are also
uniformly bounded.

 By the uniform
$\kappa$-noncollapsing and the virtue of Hamilton's compactness
theorem 16.1 in \cite{Ha6} (see \cite{CaZ} for the details on
generalizing Hamilton's compactness theorem to finite balls),
after passing to a subsequence, we can assume that the marked
sequence
$(\tilde{B}_0(\bar{x},\rho_0),\widetilde{g}^{\alpha_m\beta_m}_{ij},\bar{x})$
converges in $C^{\infty}_{loc}$ topology to a marked (noncomplete)
manifold ($B_{\infty},\widetilde{g}^{\infty}_{ij},\bar{x})$ and
the geodesic segments $\gamma$ converge to a geodesic segment
(missing an endpoint) $\gamma_{\infty}\subset B_{\infty}$
emanating from $\bar{x}$.

Clearly, the limit has restricted isotropic curvature pinching
(2.4) by the pinching assumption. Consider a tubular neighborhood
along $\gamma_{\infty}$ defined by
$$V=\bigcup_{q_0\in\gamma_{\infty}}B_{\infty}(q_0,4\pi(\widetilde{R}_{\infty}(q_0))^{-\frac{1}{2}}),$$
where $\widetilde{R}_{\infty}$ denotes the scalar curvature of the
limit and
$B_{\infty}(q_0,4\pi(\widetilde{R}_{\infty}(q_0))^{-\frac{1}{2}})$
is the ball centered at $q_0\in B_{\infty}$ with the radius
$4\pi(\widetilde{R}_{\infty}(q_0))^{-\frac{1}{2}}$. Let
$\bar{B}_{\infty}$ denote the completion of
$(B_{\infty},\widetilde{g}^{\infty}_{ij})$, and
$y_{\infty}\in\bar{B}_{\infty}$ the limit point of
$\gamma_{\infty}$. Exactly as in the second step of the proof of
Theorem 4.1, it follows from the canonical neighborhood assumption
with accuracy parameter $2\varepsilon$ that the limiting metric
$\widetilde{g}^{\infty}_{ij}$ is cylindrical at any point
$q_0\in\gamma_{\infty}$ which is sufficiently close to
$y_{\infty}$ and then the metric space $\bar{V}=V\cup
\{y_{\infty}\}$ by adding the point $y_{\infty}$ has nonnegative
curvature in Alexandrov sense. Consequently we have a
four-dimensional non-flat tangent cone $C_{y_{\infty}}\bar{V}$ at
$y_{\infty}$ which is a metric cone with aperture $\leq
20\varepsilon$.

On the other hand, note that by the canonical neighborhood
assumption, the canonical $2\varepsilon$-neck neighborhoods are
strong. Thus at each point $q\in V$ near $y_{\infty}$, the
limiting metric $\widetilde{g}^{\infty}_{ij}$ actually exists on
the whole parabolic neighborhood $$V \bigcap
P(q,0,\frac{1}{3}\eta^{-1}(\widetilde{R}_{\infty}(q))^{-\frac{1}{2}},
-\frac{1}{10}\eta^{-1}(\widetilde{R}_{\infty}(q))^{-1}),$$  and is
a smooth solution of the Ricci flow there. Pick $z\in
C_{y_{\infty}}\bar{V}$ with distance one from the vertex
$y_{\infty}$ and it is nonflat around $z$. By definition the ball
$B(z,\frac{1}{2})\subset C_{y_{\infty}}\bar{V}$ is the
Gromov-Hausdorff convergent limit of the scalings of a sequence of
balls $B_{\infty}(z_k,\sigma_k)
(\subset(V,\widetilde{g}^{\infty}_{ij}))$ where
$\sigma_k\rightarrow0$. Since the estimate (5.17) survives on
$(V,\widetilde{g}^{\infty}_{ij})$ for all $A < +\infty$, and the
tangent cone is four-dimensional and nonflat around $z$, we see
that this convergence is actually in $C^{\infty}_{loc}$ topology
and over some ancient time interval. Since the limiting
$B_{\infty}(z,\frac{1}{2})(\subset C_{y_{\infty}}\bar{V})$ is a
piece of nonnegatively (operator) curved nonflat metric cone, we
get a contradiction with Hamilton's strong maximum principle
\cite{Ha2} as before. So we have proved $\rho_0=\infty$. This
proves (i).

By the same proof of Assertion 1 in Step 2, we can further show
that for any $L$, the rescaled solutions on the balls
$\tilde{B}_{0}(\bar{x},L)$ are defined at least on the time
interval $[-\frac{1}{16}\eta^{-1}C(L)^{-1},0]$ for all
sufficiently large $m$. This proves (ii).

\vskip 0.2cm \noindent{\bf Step 4.} For any subsequence
$(\alpha_m,\beta_m)$ of $(\alpha,\beta)$ with
$r^{\alpha_m}\rightarrow 0$ and
$\widetilde{\delta}^{\alpha_m\beta_m}\rightarrow 0$ as $m\rightarrow
\infty$, by Step 3, the $\kappa$-noncollapsing and Hamilton's
compactness theorem, we can extract a $C^{\infty}_{loc}$ convergent
subsequence of $\tilde{g}^{\alpha_m\beta_m}_{ij}$ over some space
time open subsets containing $t=0$. We now want to show \textbf{any}
such limit  has bounded curvature at $t=0$. We prove by
contradiction. Suppose not, then there is a sequence of points
$z_{k}$ divergent to infinity in the limiting metric at time zero
with curvature divergent to infinity. Since the curvature at $z_k$
is large (comparable to one), $z_k$ has canonical neighborhood which
is a $2\varepsilon$-cap or strong $2\varepsilon$-neck. Note that the
boundary of $2\varepsilon$-cap lies in some $2\varepsilon$-neck. So
we get a sequence of $2\varepsilon$-necks with radius going to zero.
Note also that the limit has nonnegative sectional curvature.
Without loss of the generality, we may assume $2\varepsilon <
\varepsilon_0$, where $\varepsilon_0$ is the positive constant in
Proposition 2.2. Thus this arrives a contradiction with Proposition
2.2.

\vskip 0.2cm \noindent{\bf Step 5.} In this step, we will choose
some subsequence $(\alpha_m,\beta_m)$ of $(\alpha,\beta)$ so that
we can extract a  complete smooth limit on a time interval
$[-a,0]$ for some $a>0$ from the rescaled solutions
$\widetilde{g}^{\alpha_m \beta_m}_{ij}$ of the Ricci flow with
surgery.

Choose $\alpha_m,\beta_m\rightarrow\infty$ so that
$r^{\alpha_m}\rightarrow0$,
$\widetilde{\delta}^{\alpha_{m}\beta_{m}}\rightarrow 0$, and
Assertion 1, 2, 3
 hold with $\alpha=\alpha_m, \beta=\beta_m$ for all $A \in \{{p}/{q} \ |\ p, q=1, 2 \cdots, m\}$, and $B,C \in
\{1,2,\cdots,m\}$. By Step 3, we may assume the rescaled solutions
$\widetilde{g}^{\alpha_m \beta_m}_{ij}$ converge in
$C^{\infty}_{loc}$ topology at the time $t=0$. Since the curvature
of the limit at $t=0$ is bounded by Step 4, it follows from
Assertion 1 in Step 2 and the choice of the subsequence
$(\alpha_m,\beta_m)$ that the limiting
$(M_{\infty},\widetilde{g}^{\infty}_{ij}(\cdot,t))$ is defined at
least on a backward time interval $[-a,0]$ for some positive
constant $a$ and is a smooth solution to the Ricci flow there.

\vskip 0.2cm \noindent{\bf Step 6.} We further want to extend the
limit of Step 5 backward in time to infinity to get an ancient
$\kappa$-solution. Let $\widetilde{g}^{\alpha_m\beta_m}_{ij}$ be the
convergent sequence obtained in the above Step 5.

 Denote by
\begin{eqnarray*} t_{\max} = \sup \{\ t'&|&\mbox {we can take a
smooth limit on } (-t',0] \mbox{ (with bounded}  \\
& & \mbox{ curvature at each time slice) from a subsequence of}\\&
& \mbox{ the rescaled solutions
}\widetilde{g}^{\alpha_m\beta_m}_{ij} \}.\end{eqnarray*} We first
claim that there is a subsequence of the rescaled solutions
$\widetilde{g}^{\alpha_m\beta_m}_{ij}$ which converges in
$C^{\infty}_{loc}$ topology to a smooth limit
$(M_{\infty},\widetilde{g}^{\infty}_{ij}(\cdot,t))$ on the maximal
time interval $(-t_{\max},0]$.

Indeed, let $t_k$ be a sequence of positive numbers such that
$t_k\rightarrow t_{\max}$ and there exist smooth limits
$(M_{\infty},\widetilde{g}^{\infty}_{k}(\cdot,t))$ defined on
$(-t_k,0]$. For each $k$, the limit has nonnegative curvature
operator and has bounded curvature at each time slice. Moreover by
the gradient estimate in canonical neighborhood assumption with
accuracy parameter $2\varepsilon$, the limit has bounded curvature
on each subinterval $[-b,0]\subset(-t_k,0]$. Denote by
$\widetilde{Q}$ the scalar curvature upper bound of the limit at
time zero ($\widetilde{Q}$ independent of $k$). Then we can apply
 Li-Yau-Hamilton inequality \cite{Ha4} to get
 $$
 \widetilde{R}^{\infty}_{k}(x,t)\leq
 \frac{t_k}{t+t_k}\widetilde{Q},
 $$
 where $\widetilde{R}^{\infty}_{k}(x,t)$ are the scalar curvatures
 of the limits $(M_{\infty},\widetilde{g}^{\infty}_{k}(\cdot,t))$.
 Hence by the definition of convergence and the above curvature estimates,
  we can find a subsequence of
the rescaled solutions $\widetilde{g}^{\alpha_m\beta_m}_{ij}$
which converges in $C^{\infty}_{loc}$ topology to a smooth limit
$(M_{\infty},\widetilde{g}^{\infty}_{ij}(\cdot,t))$ on the maximal
time interval $(-t_{\max},0]$.

 We need to show $-t_{\max}=-\infty$. Suppose $-t_{\max}>-\infty$, there are
only the following two possibilities: either

(1) The curvature of the limiting solution
$(M_{\infty},\widetilde{g}^{\infty}_{ij}(\cdot,t))$ becomes
unbounded as $t\searrow -t_{\max}$; or

(2) For each small constant $\theta>0$ and each large integer
$m_0>0$, there is some $m\geq m_0$ such that the rescaled solution
$\widetilde{g}^{\alpha_m\beta_m}_{ij}$ has a surgery time
$T_m\in[-t_{\max}-\theta,0]$ and a surgery point $x_m$ lying in a
gluing cap at the times $T_m$ so that $d^2_{T_m}(x,\bar{x})$ is
uniformly bounded from above by a constant independent of $\theta$
and $m_0$.

     We next claim that the possibility (1) always occurs. Suppose not, then the curvature of the limiting solution
$(M_{\infty},\widetilde{g}^{\infty}_{ij}(\cdot,t))$ is uniformly
bounded by (some positive constant) $\hat{C}$ on $(-t_{\max},0]$. In
particular, for any $A>0$, there is a sufficiently large integer
$m_1>0$ such that  any rescaled solution
$\widetilde{g}^{\alpha_m\beta_m}_{ij}$ with $m\geq m_1$ on the
geodesic ball $\widetilde{B}_{0}(\bar{x},A)$ is defined on the time
interval $[-t_{\max}+\frac{1}{50}\eta^{-1}{\hat{C}}^{-1},0]$ and its
scalar curvature is bounded by $2\hat{C}$ there. (Here, without loss
of generality, we may assume that the upper bound $\hat{C}$ is so
large that $-t_{\max}+\frac{1}{50}\eta^{-1}{\hat{C}}^{-1} < 0$.) By
Assertion 1 in Step 2, as $m$ large enough, the rescaled solution
$\widetilde{g}^{\alpha_m\beta_m}_{ij}$ over
$\widetilde{B}_{0}(\bar{x},A)$ can be defined on the extended time
interval $[-t_{\max}-\frac{1}{50}\eta^{-1}{\hat{C}}^{-1},0]$ and
have the scalar curvature $\widetilde{R}\leq 10 \hat{C}$ on
$\widetilde{B}_{0}(\bar{x},A)\times
[-t_{\max}-\frac{1}{50}\eta^{-1}{\hat{C}}^{-1},0]$. So we can
extract a smooth limit from the sequence to get the limiting
solution which is defined on a larger time interval
$[-t_{\max}-\frac{1}{50}\eta^{-1}{\hat{C}}^{-1},0]$. This
contradicts with the definition of the maximal time $-t_{\max}$.

 We now remain to exclude the possibility (1).

 By using
Li-Yau-Hamilton inequality \cite{Ha4} again, we have
$$\widetilde{R}_{\infty}(x,t)\leq
\frac{t_{\max}}{t+t_{\max}}\widetilde{Q}.$$ So we only need to
control the curvature near $-t_{\max}$. Exactly as in the Step 4
of proof of Theorem 4.1, it follows from Li-Yau-Hamilton
inequality that
$$d_0(x,y)\leq d_t(x,y)\leq d_0(x,y)+30t_{\max}\sqrt{\widetilde{Q}}
\eqno (5.19)$$ for any $x,y\in M_{\infty}$ and $t\in
(-t_{\max},0]$.

Since the infimum of the scalar curvature is nondecreasing in
time, we have some point $y_{\infty}\in M_{\infty}$  and some time
$-t_{\max} < t_{\infty} < -t_{\max}+\frac{1}{50}\eta^{-1}$ such
that $\widetilde{R}_{\infty}(y_{\infty},t_{\infty})<{5}/{4}$. By
(5.19), there is a constant $\widetilde{A}>0$ such that
$d_t(\bar{x},y_{\infty})\leq \widetilde{A}/2$ for all $t \in
(-t_{\max},0]$.

Now we return back to the rescaled solution
$\widetilde{g}^{\alpha_m\beta_m}_{ij}$. Clearly, for arbitrarily
given small $\epsilon' > 0$, as $m$ large enough, there is a point
$y_m$ in the underlying manifold of
$\widetilde{g}^{\alpha_m\beta_m}_{ij}$ at time $0$ satisfying the
following properties
$$
\widetilde{R}(y_m,t_{\infty})<\frac{3}{2},\ \ \ \
\widetilde{d}_{t}(\bar{x},y_m)\leq \widetilde{A} \eqno (5.20)
$$
for $t\in [-t_{\max} + \epsilon',0]$. By the definition of
convergence, we know that for any fixed $A \geq 2\widetilde{A}$, as
$m$ large enough, the rescaled solution over
$\widetilde{B}_{0}(\bar{x},A)$ is defined on the time interval
$[t_{\infty},0]$ and satisfies
$$\widetilde{R}(x,t)\leq
\frac{2t_{\max}}{t+t_{\max}}\widetilde{Q}$$ on
$\widetilde{B}_{0}(\bar{x},A)\times [t_{\infty},0]$. Then by
Assertion 2 of Step 2, we have proved there is a sufficiently large
$\bar{m}_0$ such that as $m \geq \bar{m}_0$, the rescaled solutions
$\widetilde{g}^{\alpha_m\beta_m}_{ij}$ at $y_m$ can be defined on
$[-t_{\max}-\frac{1}{50}\eta^{-1},0]$, and satisfy
$$\widetilde{R}(y_m,t)\leq 15$$ for $t\in
[-t_{\max}-\frac{1}{50}\eta^{-1},t_{\infty}]$.

We now prove a statement analogous to Assertion 4 (i) of Step 3.

\vskip 0.2cm \noindent \textbf{Assertion 5.}\ \emph{ For the above
rescaled solutions $\widetilde{g}^{\alpha_m\beta_m}_{ij}$ and
$\bar{m}_0$, we have that for any $L>0$, there is a positive
constant $\omega(L)$ such that the rescaled solutions
$\widetilde{g}^{\alpha_m\beta_m}_{ij}$ satisfy
$$\widetilde{R}(x,t)\leq
\omega(L)$$ for all $(x,t)$ with $\tilde{d}_{t}(x,y_m)\leq L$ and
$t\in [-t_{\max}-\frac{1}{50}\eta^{-1},t_{\infty}]$ and for all $m
\geq \bar{m}_0$.} \vskip 0.2cm

\noindent \textbf{Proof of Assertion 5.}  We slightly modify the
argument in the proof of Assertion 4 (i). Let

\begin{eqnarray*}M(\rho)=\sup\{\widetilde{R}(x,t)&|& \tilde{d}_t(x,y_m)\leq\rho \mbox{
and } t \in
[-t_{\max}-\frac{1}{50}\eta^{-1},t_{\infty}]\\
& & \mbox{ in the rescaled solutions
 } \ \widetilde{g}^{\alpha_m\beta_m}_{ij}, m \geq \bar{m}_0\}\end{eqnarray*}
 and
$$\rho_0=\sup\{\rho > 0\ |\  M(\rho)<+\infty\}.$$  Note that the
estimate (5.17) implies that $\rho_0>0$. We only need to show
$\rho_0 = +\infty$.

We argue by contradiction. Suppose $\rho_0<+\infty$. Then, after
passing to subsequence, there are a sequence of
$(\tilde{y}_m,t_m)$ in the rescaled solutions
$\widetilde{g}^{\alpha_m\beta_m}_{ij}$ with $t_m \in
[-t_{\max}-\frac{1}{50}\eta^{-1},t_{\infty}]$ and
$\tilde{d}_{t_m}(y_m,\tilde{y}_m)\rightarrow\rho_0<+\infty$ such
that $\widetilde{R}(\tilde{y}_m,t_m)\rightarrow +\infty$. Denote
by $\gamma_m$ a minimizing geodesic segment from $y_m$ to
$\tilde{y}_m$ at the time $t_m$ and denote by
$\widetilde{B}_{t_m}(y_m,\rho_0)$ the geodesic open ball centered
at $y_m$ of radius $\rho_0$ on the rescaled solution
$\widetilde{g}^{\alpha_m\beta_m}_{ij}(\cdot,t_m)$.

 For any $0<\rho<\rho_0$ with $\rho$ near $\rho_0$, by applying Assertion 3 as before, we get that  the rescaled solutions on the
balls $\widetilde{B}_{t_m}(y_m,\rho)$ are defined on the time
interval $[t_m-\frac{1}{16}\eta^{-1}M(\rho)^{-1},t_m]$ for all large
$m$. And by Step 1 and Shi's derivative estimate, we further know
that the covariant derivatives of the curvatures of all order on
$\widetilde{B}_{t_m}(y_m,\rho-\frac{(\rho_0-\rho)}{2})\times[t_m-\frac{1}{32}\eta^{-1}M(\rho)^{-1},t_m]$
are also uniformly bounded. Then by the uniform
$\kappa$-noncollapsing and Hamilton's compactness theorem, after
passing to a subsequence, we can assume that the marked sequence
$(\tilde{B}_{t_m}(y_m,\rho_0),\widetilde{g}^{\alpha_m\beta_m}_{ij}(\cdot,t_m),y_m)$
converges in $C^{\infty}_{loc}$ topology to a marked (noncomplete)
manifold ($B_{\infty},\widetilde{g}^{\infty}_{ij},y_{\infty})$ and
the geodesic segments $\gamma_m$ converge to a geodesic segment
(missing an endpoint) $\gamma_{\infty}\subset B_{\infty}$ emanating
from $y_{\infty}$.

Clearly, the limit also has restrictive isotropic curvature pinching
(2.4). Then by repeating the same argument as in the proof of
Assertion 4 (i) in the rest, we derive a contradiction with
Hamilton's strong maximum principle. This proves Assertion 5. \vskip
0.2cm

 We then apply the second estimate of (5.20) and Assertion 5 to conclude that for any large constant
$0<A<+\infty$, there is a positive constant $C(A)$ such that for
any small $\epsilon'>0$, the rescaled solutions
$\widetilde{g}^{\alpha_m\beta_m}_{ij}$ satisfy
$$\widetilde{R}(x,t)\leq C(A), \eqno (5.21)$$
for all $x \in \widetilde{B}_{0}(\bar{x},A)$ and $t\in
  [-t_{\max} + \epsilon',0]$, and for all sufficiently large $m$. Then by applying Assertion 1
  in Step 2, we conclude that the rescaled solutions $\widetilde{g}^{\alpha_m\beta_m}_{ij}$ on the geodesic
  balls
  $\widetilde{B}_{0}(\bar{x},A)$ are also defined on the extended
time interval $ [-t_{\max} +
\epsilon'-\frac{1}{8}\eta^{-1}C(A)^{-1},0]$ for all sufficiently
large $m$. Furthermore, by the gradient estimates as in Step 1, we
have
$$\widetilde{R}(x,t)\leq 10C(A),$$
for $x \in \widetilde{B}_{0}(\bar{x},A)$ and $t\in
  [-t_{\max} + \epsilon'-\frac{1}{8}\eta^{-1}C(A)^{-1},0]$.
Since $\epsilon'>0$ is arbitrarily small, the rescaled solutions
$\widetilde{g}^{\alpha_m\beta_m}_{ij}$ on
  $\widetilde{B}_{0}(\bar{x},A)$ are defined on the extended
time interval $ [-t_{\max} -\frac{1}{16}\eta^{-1}C(A)^{-1},0]$ and
satisfy
$$\widetilde{R}(x,t)\leq 10C(A), \eqno (5.22)$$
for $x \in \widetilde{B}_{0}(\bar{x},A)$ and $t\in
  [-t_{\max} -\frac{1}{16}\eta^{-1}C(A)^{-1},0]$, and for all sufficiently large $m$.

   Now, by taking convergent subsequences from the rescaled solutions $\widetilde{g}^{\alpha_m\beta_m}_{ij}$, we see that
   the limit solution is defined smoothly on a space-time open subset of $M_{\infty}\times (-\infty,0]$
   containing $M_{\infty}\times [-t_{\max},0]$. By Step 4,
   we see that the limiting
metric
 $\widetilde{g}^{\infty}_{ij}(\cdot,-t_{\max})$ at time $-t_{\max}$
 has bounded curvature. Then by combining with the $2\varepsilon$-canonical neighborhood
 assumption we conclude that the curvature of the limit is uniformly bounded
 on the time interval
 $[-t_{\max},0]$. So we have excluded the possibility (1).

Hence we have proved a subsequence of the rescaled solutions
converges to an  ancient $\kappa$-solution.

Finally by combining with the canonical neighborhood theorem of
ancient $\kappa$-solutions with restricted isotropic curvature
pinching condition (Theorem 3.8) and the same argument in the
second paragraph of Section 4, we see that $(\bar{x},\bar{t})$ has
a canonical neighborhood with parameter $\varepsilon$, which is a
contradiction. Therefore we have completed the proof of the
proposition.$$\eqno \#$$ \vskip 0.3cm

Summing up, we have proved that for an arbitrarily given compact
four-manifold with positive isotropic curvature and with no
essential incompressible space form, there exist non-increasing
positive (continuous) functions $\widetilde{\delta}(t)$ and
$\widetilde{r}(t)$, defined on $[0,+\infty)$, such that for
arbitrarily given positive (continuous) function $\delta(t)$ with
$\delta(t)<\widetilde{\delta}(t)$ on $[0,+\infty)$, the Ricci flow
with surgery, with the given four-manifold as initial datum, has a
solution on a maximal time interval $[0,T)$, with $T\leq
2/R_{min}(0)<+\infty$, obtained by evolving the Ricci flow and by
performing ${\delta}$-cutoff surgeries at a sequence of times
$0<t_1<t_2<\cdots<t_i<\cdots<T$ with $\delta(t_i)\leq\delta\leq
\widetilde{\delta}(t_i)$ at each time $t_i$, so that the pinching
assumption and the canonical neighborhood assumption with
$r=\widetilde{r}(t)$ are satisfied. (At this moment we still do not
know whether the surgery times $t_i$ are discrete).

Clearly, the upper derivative of the volume in time satisfies
$$\frac{d}{dt}V(t)\leq 0$$ since the scalar curvature is
nonnegative. Thus $$V(t)\leq V(0)$$ for all $t\in [0,T)$. Also note
that at each time $t_i$, the volume which is cut down by
$\delta(t_i)$-cutoff surgery is at least an amount of $h^4(t_i)$
with $h(t_i)$ depending only on ${\delta}(t_i)$ and
$\widetilde{r}(t_i)$ (by Lemma 5.2). Thus the set of the surgery
times $\{t_i\}$ must be finite. So we have proved the following
long-time existence result.\vskip 0.3cm

$\underline{\mbox{\textbf{Theorem 5.6}}}$ \emph{ Given a compact
four-dimensional Riemannian manifold with positive isotropic
curvature and with no essential incompressible space form, and given
any fixed small constant $\varepsilon > 0$, there exist
non-increasing positive (continuous) functions
$\widetilde{\delta}(t)$ and $\widetilde{r}(t)$, defined on
$[0,+\infty)$, such that for arbitrarily given positive (continuous)
function $\delta(t)$  with $\delta(t) \leq \widetilde{\delta}(t)$ on
$[0,+\infty)$, the Ricci flow with surgery, with the given
four-manifold as initial datum, has a solution satisfying the the
pinching assumption and the canonical neighborhood assumption (with
accuracy $\varepsilon$) with $r=\widetilde{r}(t)$ on a maximal time
interval $[0,T)$ with $T<+\infty$ and becoming extinct at $T$, which
is obtained by evolving the Ricci flow and by performing a finite
number of cutoff surgeries with each $\delta$-cutoff at time $t\in
(0,T)$ having $\delta=\delta(t)$.} \emph{Consequently, the initial
manifold is diffeomorphic to a connected sum of a finite copies of
$\mathbb{S}^4$, $\mathbb{RP}^4$, $\mathbb{S}^3 \times \mathbb{S}^1$,
and $\mathbb{S}^3\widetilde{\times} \mathbb{S}^1$.}
$$\eqno \#$$ \vskip 0.3cm

Finally, the main theorem (Theorem 1.1) stated in Section 1 is a
direct consequence of the above theorem. \vskip 1cm

  \centerline{\large{\textbf{Appendix.  Standard Solutions}}}
 \vskip 0.3cm
 In this appendix, we will prove the curvature estimates for the standard
 solutions, and give a canonical neighborhood description for the standard
 solution in dimension four. We have used these estimates and the description in Section 5
 for the
 surgery arguments. The curvature estimate for the special case that
 the
 dimension is three and the initial metric is rotationally symmetric,
  was earlier claimed
 by Perelman in \cite{P2}.
\vskip 0.2cm

 $\underline{\mbox{\textbf{Theorem A.1}}}$\emph{
   Let $g_{ij}$ be a complete Riemannian metric on $\mathbb{R}^{n}$ $(n>2)$
  with
  nonnegative curvature operator and with positive scalar curvature which is
  asymptotic to a round
  cylinder of scalar curvature 1 at infinity.
  Then there is a complete
  solution $g_{ij}({\cdot,t})$ to the Ricci flow, with $g_{ij}$ as initial metric,
  which
  exists on the time interval $[0,\frac{n-1}{2})$, has bounded curvature in
  each closed time interval
   $[0,t] \subset [0,\frac{n-1}{2})$,
  and satisfies the estimate
  $$
  R(x,t)\geq \frac{C^{-1}}{\frac{n-1}{2}-t}
  $$
  for some $C$ depending only on the initial metric $g_{ij}$.}
\vskip 0.1cm

  $\underline{\mbox{\textbf{Proof}}}$. Since the initial metric has bounded
  curvature operator and has a positive lower bound on its scalar
  curvature, by \cite{Sh1} and the maximum principle,
  the Ricci flow has a
  solution $g(\cdot,t)$ on a maximal time interval $[0,T)$ with
  $T<\infty$. By Hamilton's maximum principle, the solution $g(x,t)$ has nonnegative curvature operator for $t>0$.
  Note that the injectivity
  radius of the initial metric has a positive lower bound, so by the same proof of
  Perelman's no local collapsing theorem I (in the section 7.3 of \cite{P1}, or see the proof of Theorem 3.5 of this paper),
   there is a $\kappa=\kappa(T,g_{ij})>0$ such
  that $g_{ij}({\cdot,t})$ is $\kappa$-noncollapsed on the scale
  $\sqrt{T}$.

  We will firstly prove the following assertion. \vskip 0.2cm

 $\underline{\mbox{\textbf{Claim 1}}}$  \emph{There is a positive function
  $\omega:[0,\infty)\longrightarrow[0,\infty)$ depending only on
  the initial metric and $\kappa$ such that
  $$
  R(x,t)\leq R(y,t)\omega(R(y,t)d_{t}^{2}(x,y))
  $$
  for all $x,y\in M^{n}=\mathbb{R}^{n}$, $t\in[0,T)$.
}\vskip 0.2cm

 The proof is similar to that of Proposition 3.3.
  Notice that the initial metric has nonnegative
  curvature operator and
 its scalar caurvature satisfies $$C^{-1}\leq R(x)
 \leq C \eqno\mbox{(A.1)}$$ for some positive constant
  $C>1$. By maximum
 principle, we know $T\geq\frac{1}{2nC}$ and $R(x,t)\leq 2C$
 for $t\in[0,\frac{1}{4nC}]$. The assertion is clearly true for
 $t\in[0,\frac{1}{4nC}]$.

 Now fix $(y,t_{0})\in M^{n}\times[0,T)$
 with $t_{0}\geq\frac{1}{4nC}$. Let $z$ be the closest point to
 $y$  with the property $R(z,t_{0})d^{2}_{t_{0}}(z,y)=1$ (at time
 $t_{0}$). Draw a shortest geodesic from $y$ to $z$ and choose a
 point $\tilde{z}$ on the geodesic satisfying $d_{t_{0}}(z,\tilde{z})
 =\frac{1}{4}R(z,t_{0})^{-\frac{1}{2}}$, then we have
 $$
 R(x,t_{0})\leq\frac{1}{(\frac{1}{2}R(z,t_{0})^{-\frac{1}{2}})^{2}},
    \ \ \ \ \ \ \ \  \mbox{on}\ \
    B_{t_{0}}(\tilde{z},\frac{1}{4}R(z,t_{0})^{-\frac{1}{2}})
     $$

Note that $R(x,t)\geq C^{-1}$ everywhere by the evolution equation
of the scalar curvature. Then by Li-Yau-Hamilton inequality
\cite{Ha4}, for all $(x,t)\in
B_{t_{0}}(\tilde{z},\frac{1}{8nC}R(z,t_{0})^{-\frac{1}{2}})
\times[t_{0}-(\frac{1}{8nC}R(z,t_{0})^{-\frac{1}{2}})^{2},t_{0}]
     $, we have
     \begin{equation*}
     \begin{split}
R(x,t)&\leq
(\frac{t_{0}}{t_{0}-(\frac{1}{8n\sqrt{C}})^{2}})\frac{1}{(\frac{1}{2}R(z,t_{0})^{-\frac{1}{2}})^{2}},\\
 &\leq [\frac{1}{8nC}R(z,t_{0})^{-\frac{1}{2}}]^{-2}
      \end{split}
     \end{equation*}
     Combining this with the $\kappa$-noncollapsing, we have
$$Vol (B_{t_{0}}(\tilde{z},\frac{1}{8nC}R(z,t_{0})^{-\frac{1}{2}}))
\geq \kappa (\frac{1}{8nC}R(z,t_{0})^{-\frac{1}{2}})^{n}$$ and
then
 $$
 Vol (B_{t_{0}}(z,8R(z,t_{0})^{-\frac{1}{2}}))
\geq \kappa (\frac{1}{64nC})^{n}(8R(z,t_{0})^{-\frac{1}{2}})^{n}
$$
So by Corollary 11.6 (b) of  \cite{P1}, there hold
$$
R(x,t_{0})\leq C(\kappa)R(z,t_{0}),\ \ \ \text{for all } x\in
B_{t_{0}}(z,2R(z,t_{0})^{-\frac{1}{2}}).
$$
Here in the following we denote by $C(\kappa)$ various positive
constants depending only on $\kappa,n$ and the initial metric.

Now by Li-Yau-Hamilton inequality \cite{Ha4} and local gradient
estimate of Shi \cite{Sh1}, we obtain
\begin{align*}
R(x,t)\leq C(\kappa)R(z,t_{0}),\ \  \text{and}\ \
|\frac{\partial}{\partial t}R|(x,t)\leq C(\kappa)(R(z,t_{0}))^{2}
\end{align*}
for all $ (x,t)\in
B_{t_{0}}(z,2R(z,t_{0})^{-\frac{1}{2}}))\times[t_{0}-
(\frac{1}{8nC}R(z,t_{0})^{-\frac{1}{2}})^{2},t_{0}]$. Therefore by
combining with the Harnack estimate \cite{Ha4}, we obtain
\begin{align*}
R(y,t_{0})&\geq C(\kappa)^{-1}R(z,t_{0}-C(\kappa)^{-1}R(z,t_{0})^{-1})\\
&\geq C(\kappa)^{-2} R(z,t_{0})
\end{align*}

Consequently, we have showed that there is a constant $C(\kappa)$
such that
$$
 Vol (B_{t_{0}}(y,R(y,t_{0})^{-\frac{1}{2}}))
\geq C(\kappa)^{-1}(R(y,t_{0})^{-\frac{1}{2}})^{n}
$$
and
$$
R(x,t_{0})\leq C(\kappa)R(y,t_{0}) \ \ \text{for all}\ \
 x\in
B_{t_{0}}(y,R(y,t_{0})^{-\frac{1}{2}}).
$$
In general, for any $r\geq R(y,t_{0})^{-\frac{1}{2}}$, we have
$$
 Vol (B_{t_{0}}(y,r))
\geq C(\kappa)^{-1}(r^{2}R(y,t_{0}))^{-\frac{n}{2}}r^{n}.
$$
By applying Corollary 11.6 of \cite{P1} again, there exists a
positive constant $\omega(r^{2}R(y,t_{0}))$ depending only on the
constant
 $r^{2}R(y,t_{0})$ and $\kappa$ such that
$$
R(x,t_{0})\leq R(y,t_{0})\omega(r^{2}R(y,t_{0})),\ \ \ \ \mbox{for
all}\  x\in B_{t_{0}}(y,\frac{1}{4}r).
$$
 This proves the desired Claim 1. \vskip 0.2cm

 Now we study the asymptotic behavior of the solution at
 infinity. For any $0<t_{0}<T$, we know that the metrics $g_{ij}(x,t)$ with
 $t\in [0,t_{0}]$ has uniformly bounded curvature by the definition of $T$.
 Let $x_{k}$ be
 a sequence of points with $d_{0}(x_{0},x_{k})\mapsto \infty$. By
 Hamilton's compactness theorem \cite{Ha5}, after taking a subsequence, $g_{ij}(x,t)$ around $x_{k}$ will
 converge to a solution to the Ricci flow on $\mathbb{R}\times \mathbb{S}^{n-1}$
 with round cylinder metric of scalar curvature 1 as initial
 data.  Denote the limit by $\tilde{g}_{ij}$. Then by the uniqueness
 theorem in \cite{CZ2}, we have
 $$
 \tilde{R}(x,t)=\frac{\frac{n-1}{2}}{\frac{n-1}{2}-t},\ \ \text{for all
 }\ \ t\in[0,t_{0}].
 $$
 It follows that $T\leq \frac{n-1}{2}$. In order to show
 $T=\frac{n-1}{2}$, it suffices to prove the following assertion \vskip 0.2cm

 $\underline{\mbox{\textbf{Claim 2.}}}$\emph{
 Suppose $T<\frac{n-1}{2}$. Fix a point $x_{0}\in M^{n} $, then there is a $\delta>0$, such that
 for any $x\in M$ with $d_{0}(x,x_{0})\geq \delta^{-1}$, we have
 $$
 R(x,t)\leq 2C+\frac{n-1}{\frac{n-1}{2}-t}\ \ \ \ \text{ for all
 }\ \ t\in[0,T)
 $$
 where $C$ is the constant in }(A.1). \vskip 0.2cm

  In view of Claim 1, if Claim 2 holds, then
  \begin{align*}
  \sup_{M^{n}\times[0,T)}R(y,t)&\leq
  \omega(\delta^{-2}(2C+\frac{n-1}{\frac{n-1}{2}-T}))(2C+\frac{n-1}{\frac{n-1}{2}-T})\\
  &< \infty
  \end{align*}
  which will contradict with the definition of $T$.

  To show Claim 2, we argue by contradiction. Suppose for each
  $\delta>0$, there is a $(x_{\delta},t_{\delta})$ with $0<t_{\delta}<T
  $ such that
  $$
  R(x_{\delta},t_{\delta})>2C+\frac{n-1}{\frac{n-1}{2}-t_{\delta}}\ \
   \mbox{and}\ \  d_{0}(x_{\delta},x_{0})\geq \delta^{-1}.$$
  Let
  $$
  \bar{t}_{\delta}=\sup\{t \ | \sup_{M^{n}\setminus B_{0}(x_{0},\delta^{-1})}
  R(y,t)<2C+\frac{n-1}{\frac{n-1}{2}-t}\}.
  $$
Since
$\lim\limits_{d_{0}(y,x_{0})\rightarrow\infty}R(y,t)={\frac{n-1}{2}}/
{(\frac{n-1}{2}-t)}$ and $\sup_{M\times[0,\frac{1}{4nC}]}R(y,t)\leq
2C$, we know $\frac{1}{4nC}\leq \bar{t}_{\delta}\leq t_{\delta}$ and
there is a $\bar{x}_{\delta}$ such that
$d_{0}(x_{0},\bar{x}_{\delta})\geq\delta^{-1}$ and
$R(\bar{x}_{\delta},\bar{t}_{\delta})=2C+{n-1}/{(\frac{n-1}{2}-\bar{t}_{\delta})}.$
By Claim 1 and Hamilton's compactness theorem \cite{Ha5}, as
$\delta\ \rightarrow 0$ and after taking subsequence, the metrics
$g_{ij}(x,t)$ on $B_{0}(\bar{x}_{\delta},\frac{\delta^{-1}}{2})$
over the time interval $[0,\bar{t}_{\delta}]$ will converge to a
solution $\tilde{g}$ on $\tilde{M}=\mathbb{R}\times
\mathbb{S}^{n-1}$ with standard metric of scalar curvature 1 as
initial datum over the time interval $[0,\bar{t}_{\infty}]$, and its
scalar curvature satisfies
\begin{eqnarray*}
 \tilde{R}(\bar{x}_{\infty},\bar{t}_{\infty})&=&2C+\frac{n-1}{\frac{n-1}{2}-\bar{t}_{\infty}},\\
 \tilde{R}(x,t)&\leq&2C+\frac{n-1}{\frac{n-1}{2}-\bar{t}_{\infty}},
 \ \ \ \  \mbox{for\  all}\  t\in [0,\bar{t}_{\infty}],
\end{eqnarray*}
where $(\bar{x}_{\infty},\bar{t}_{\infty})$ is the limit of
$(\bar{x}_{\delta},\bar{t}_{\delta})$. On the other hand, by the
uniqueness theorem in \cite{CZ2} again, we know
$$
\tilde{R}(\bar{x}_{\infty},\bar{t}_{\infty})=\frac{\frac{n-1}{2}}{\frac{n-1}{2}-\bar{t}_{\infty}}
$$
which is a contradiction.  Hence we have proved Claim 2 and then
have verified $T=\frac{n-1}{2}$.\vskip 0.2cm

Now we are ready to show
$$
R(x,t)\geq \frac{\tilde{C}^{-1}}{\frac{n-1}{2}-t}, \ \ \ \text{for
all } (x,t)\in M^{n}\times[0,\frac{n-1}{2}), \eqno\mbox{(A.2)}
$$
for some positive constant $\tilde{C}$ depending only on the
initial metric.

For any $(x,t)\in M^{n}\times[0,\frac{n-1}{2})$, by Claim 1 and
$\kappa$-noncollapsing, there is a constant $C(\kappa)>0$ such
that
$$
Vol_{t}(B_{t}(x,{R(x,t)}^{-\frac{1}{2}}))\geq
C(\kappa)^{-1}({R(x,t)}^{-\frac{1}{2}})^{n}.
$$
Then by the volume estimate of Calabi-Yau \cite{ScY} on manifolds
with nonnegative Ricci curvature, for any $a\geq 1$, we have
$$
Vol_{t}(B_{t}(x,a{R(x,t)}^{-\frac{1}{2}}))\geq
C(\kappa)^{-1}\frac{a}{8n}({R(x,t)}^{-\frac{1}{2}})^{n}.
$$
On the other hand, since $(M^{n},g_{ij}(\cdot,t))$ is asymptotic
to a cylinder of scalar curvature
${\frac{n-1}{2}}/{(\frac{n-1}{2}-t)}$ , for sufficiently large
$a>0$, we have
$$
Vol_{t}(B_{t}(x,a\sqrt{\frac{n-1}{2}-t}))\leq
C(n)a(\frac{n-1}{2}-t)^{\frac{n}{2}}.
$$
Combining these two inequalities, we have for  all sufficiently
large $a$:
\begin{eqnarray*}
C(n)a(\frac{n-1}{2}-t)^{\frac{n}{2}}&\geq&
Vol_{t}(B_{t}(x,a(\frac{\sqrt{\frac{n-1}{2}-t}}{R(x,t)^{-\frac{1}{2}}}){R(x,t)}^{-\frac{1}{2}}))\\
&\geq&
C(\kappa)^{-1}\frac{a}{8n}(\frac{\sqrt{\frac{n-1}{2}-t}}{{R(x,t)}^{-\frac{1}{2}}})
({R(x,t)}^{-\frac{1}{2}})^{n},
\end{eqnarray*}
which gives the desired estimate (A.2). Therefore we complete the
proof of the theorem.
$$\eqno \#$$ \vskip 0.3cm

 We now fix a standard capped infinite cylinder metric on $\mathbb{R}^{4}$ as follows.
 Consider the semi-infinite standard round cylinder $N_0 =
 \mathbb{S}^3 \times (-\infty,4)$ with the metric $g_0$ of scalar curvature 1.
 Denote by $z$ the coordinate of
 the second factor $(-\infty,4)$.
 Let $f$ be a smooth nondecreasing convex function on
$(-\infty,4)$ defined by
$$
   \left\{
   \begin{array}{lll}
    f(z) = 0, \ \ \ z\leq 0,
         \\[3mm]
    f(z) = ce^{-\frac{D}{z}}, \ \ \ z \in (0,3], \\[3mm]
    f(z) \mbox{ is strictly convex on } z \in [3,3.9], \\[3mm]
    f(z) = -\frac{1}{2}\log(16-z^2), \ \ \ z \in [3.9,4),
\end{array}
\right.
$$
where the small (positive) constant $c=c_0$ and big (positive)
constant $D=D_0$ are fixed as in Lemma 5.3. Let us replace the
standard metric $g_0$ on the portion $\mathbb{S}^3 \times [0,4)$ of
the semi-infinite cylinder by $\hat{g} = e^{-2f}g_0$. Then the
resulting metric $\hat{g}$ will be smooth on $\mathbb{R}^4$ obtained
by adding a point to $\mathbb{S}^3 \times (-\infty,4)$ at $z=4$. We
denote the manifold by $(\mathbb{R}^4,\hat{g})$. \vskip0.2cm

Next we will consider the ``canonical neighborhood" decomposition of
the fixed standard solution with $(\mathbb{R}^4,\hat{g})$ as initial
metric.

\vskip 0.3cm \noindent \textbf{Corollary A.2.} \emph{
   Let $g_{ij}(x,t)$ be
  the above fixed standard solution to the Ricci flow on $\mathbb{R}^4\times
  [0,\frac{3}{2})$.
  Then for any $\varepsilon>0$, there is a positive constant $C(\varepsilon)$ such that
  each point
$(x,t) \in \mathbb{R}^4\times
  [0,\frac{3}{2})$ has an open neighborhood $B$, with $B_t(x,r)\subset B \subset
  B_t(x,2r)$ for some
$0<r<C(\varepsilon)R(x,t)^{-\frac{1}{2}}$, which falls into one of
the following two categories: either}

 \emph{(a) $B$ is  an
$\varepsilon$-cap, or}

\emph{(b) $B$ is an $\varepsilon$-neck and it is the slice at the
time $t$ of the parabolic neighborhood
$P(x,t,\varepsilon^{-1}R(x,t)^{-\frac{1}{2}},-\min
\{R(x,t)^{-1},t\})$, on which the standard solution is, after
scaling with the factor $R(x,t)$ and shifting the time $t$ to zero,
$\varepsilon$-close (in $C^{[\varepsilon^{-1}]}$ topology) to the
corresponding subset of the evolving standard cylinder $\mathbb{S}^3
\times \mathbb{R}$ over the time interval $[-\min \{tR(x,t),1\},0]$
with scalar curvature $1$ at the time zero.}

\vskip 0.2cm \noindent{\textbf{Proof}}. First, we discuss the
curvature pinching of this fixed standard solution.
 Because the initial metric is asymptotic to a cylinder,
we have a uniform isotropic curvature pinching at initial, that is
to say, there is a universal constant $\Lambda'>0$ such that

$$\max\{a_3,b_3,c_3\}\leq \Lambda' a_1 \mbox{ and
}\max\{a_3,b_3,c_3\}\leq \Lambda' c_1.$$
 Moreover since the initial metric has nonnegative curvature operator,
 we have
 $b_3^{2}\leq a_1c_1.$
 By the pinching estimates of Hamilton \cite{Ha2} \cite{Ha7},
 $b_3^{2}\leq a_1c_1$ is
  preserved, and the following two estimates
  are also preserved
$$\max\{a_3,b_3,c_3\}\leq \max\{\Lambda',5\} a_1 \mbox{ and
}\max\{a_3,b_3,c_3\}\leq \max\{\Lambda',5\} c_1, $$ under the
Ricci flow.

The proof of the lemma is reduced to two assertions. We now state
and prove the first assertion which takes care of those points
with times close to $\frac{3}{2}$.

     \vskip 0.2cm\noindent\textbf{Assertion 1.}  \emph{ For any $\varepsilon>0$, there is a positive
 number $\theta=\theta(\varepsilon)$ with $0<\theta<\frac{3}{2}$
 such that for any $(x_0,t_0)\in \mathbb{R}^4\times
 [\theta,\frac{3}{2})$, the standard solution on the parabolic
 neighborhood $$P(x_0,t_0,\varepsilon^{-1}R(x_0,t_0)^{-\frac{1}{2}},
 -\varepsilon^{-2}R(x_0,t_0)^{-1})$$ is well-defined and is, after
 scaling with the factor $R(x_0,t_0)$, $\varepsilon$-close
 (in $C^{[\varepsilon^{-1}]}$ topology) to the
 corresponding subset of some oriented ancient-$\kappa$
 solution with restricted isotropic curvature pinching (2.4).}\vskip 0.2cm

We argue by contradiction. Suppose the Assertion 1 is not true, then
there exists $\varepsilon_{0}>0$ and a sequence of points
$(x_{k},t_k)$ with $t_{k}\rightarrow \frac{3}{2}$, such that the
standard solution on the parabolic neighborhoods
$$P(x_k,t_k,\varepsilon_{0}^{-1}R(x_k,t_k)^{-\frac{1}{2}},
 -\varepsilon_{0}^{-2}R(x_k,t_k)^{-1})$$ is not, after
 scaling by the factor $R(x_k,t_k)$, $\varepsilon_0$-close to the
 corresponding subset of any ancient $\kappa$-solution. Note that by Theorem A.1, there is a constant
 $C>0$ (depending only on the initial metric, hence it is universal) such that
 $R(x,t)\geq C^{-1}/(\frac{3}{2}-t)$. This implies
  $$\varepsilon_{0}^{-2}R(x_k,t_k)^{-1}\leq
  C\varepsilon_{0}^{-2}(\frac{3}{2}-t_k)<t_k ,$$
  and then the standard solution on the parabolic
neighborhoods
$$P(x_k,t_k,\varepsilon_{0}^{-1}R(x_k,t_k)^{-\frac{1}{2}},
 -\varepsilon_{0}^{-2}R(x_k,t_k)^{-1})$$ is well-defined as $k$
 large. By Claim 1 in Theorem A.1, there is a  positive function
  $\omega:[0,\infty)\rightarrow[0,\infty)$
   such that
  $$
  R(x,t_k)\leq R(x_k,t_k)\omega(R(x_k,t_k)d_{t_k}^{2}(x,x_k))
    $$
  for all $x\in \mathbb{R}^4$.
   Now by scaling the standard solution $g_{ij}(\cdot,t)$
  around $(x_k,t_k)$ with the factor $R(x_k,t_k)$ and shifting the time $t_k$ to zero,
  we get a sequence of the rescaled solutions to the
  Ricci flow  $
  \tilde{g}^{k}_{ij}(x,\tilde{t})=R(x_k,t_k)g_{ij}(x,t_k+\tilde{t}/R(x_k,t_k))$
defined on $\mathbb{R}^4$ with $\tilde{t}\in [-R(x_k,t_k)t_k,0]$. We
denote the scalar curvature and the distance of the rescaled metric
$\tilde{g}^{k}_{ij}$ by $\tilde{R}^{k}$ and $\tilde{d}$. By
combining with the Claim 1 in Theorem A.1 and the Li-Yau-Hamilton
inequality, we get
 \begin{eqnarray*}
  \tilde{R}^{k}(x,0)&\leq& \omega(\tilde{d}_{0}^{2}(x,x_k))\\
\tilde{R}^{k}(x,\tilde{t})&\leq&
\frac{R(x_k,t_k)t_k}{\tilde{t}+R(x_k,t_k)t_k}\omega(\tilde{d}_{0}^{2}(x,x_k))
    \end{eqnarray*}
    for any $x\in \mathbb{R}^4$ and $\tilde{t}\in (-R(x_k,t_k)t_k,0]$.
    Note that
    $R(x_k,t_k)t_k\rightarrow \infty$ by Theorem A.1. We have shown in the proof of Theorem A.1 that
    the standard solution
    is $\kappa$-noncollapsed on  all scales less than $1$ for some $\kappa>0$. Then from the
    $\kappa$-noncollapsing, the above curvature estimates and Hamilton's compactness theorem (Theorem 16.1 of \cite{Ha6}),
     we know $\tilde{g}^{k}_{ij}(x,\tilde{t})$ has a
    convergent subsequence (as $k\rightarrow\infty$) whose limit is an
    ancient,
    $\kappa$-noncollapsed, complete and oriented solution with nonnegative curvature operator. This
    limit must has bounded curvature by the same proof of Step
    3 in Theorem 4.1. It also satisfies the restricted isotropic pinching condition
    (2.4). This gives a contradiction. The Assertion
    1 is proved. \vskip 0.2cm

  We now fix the constant $\theta(\varepsilon)$ obtained in Assertion 1.
   Let $O$ be the tip of the
manifold $\mathbb{R}^4$ (it is rotationally symmetric about $O$ at
time $0$, it remains so as $t>0$ by the uniqueness Theorem
\cite{CZ2}). \vskip 0.2cm\noindent\textbf{Assertion 2} \emph{ There
are constants $B_1(\varepsilon)$, $B_2(\varepsilon)$ depending only
on $\varepsilon$, such that if $(x_0,t_0)\in M\times
 [0,\theta)$ with $d_{t_0}(x_0,O)\leq B_{1}(\varepsilon)$, then there is a $0<r<B_2(\varepsilon)$
 such that $B_{t_0}(x_0,r)$ is an $\varepsilon$-cap; if $(x_0,t_0)\in
M\times
 [0,\theta)$ with $d_{t_0}(x_0,O)\geq B_1(\varepsilon)$, then the parabolic neighborhood
 $$P(x_0,t_0,\varepsilon^{-1}R(x_0,t_0)^{-\frac{1}{2}},
 -\min \{R(x_0,t_0)^{-1},t_0 \})$$ is after
scaling with the factor $R(x_0,t_0)$ and shifting the time $t_0$ to
zero, $\varepsilon$-close (in $C^{[\varepsilon^{-1}]}$ topology) to
the corresponding subset of the evolving standard cylinder
$\mathbb{S}^3 \times \mathbb{R}$ over the time interval $[-\min
\{t_0R(x_0,t_0),1\},0]$ with scalar curvature $1$ at the time
zero.}\vskip 0.2cm

Since the standard solution exists on the time interval
$[0,\frac{3}{2})$, there is a constant $B_{0}(\varepsilon)$ such
that the curvatures on $[0,\theta(\varepsilon)]$ are uniformly
bounded by $B_{0}(\varepsilon)$. This implies that the metrics in
$[0,\theta(\varepsilon)]$ are equivalent.  Note that the initial
metric is asymptotic to a standard cylinder. For any
 sequence of points $x_{k}$ with $d_{0}(O,x_{k})\rightarrow
 \infty$,
 after taking a subsequence, $g_{ij}(x,t)$ around $x_{k}$ will
 converge to a solution to the Ricci flow on $\mathbb{R}\times \mathbb{S}^{3}$
 with round cylinder metric of scalar curvature 1 as initial
 data. By the uniqueness theorem \cite{CZ2}, the limit solution must be the standard evolving round cylinder.
  This implies that
  there is a
constant $B_{1}(\varepsilon)>0$ depending on $\varepsilon$ such that
for any $(x_0,t_0)$ with $t_0\leq \theta(\varepsilon)$ and
$d_{t_0}(x,O)\geq B_{1}(\varepsilon)$, the standard solution on the
parabolic neighborhood
$P(x_0,t_0,\varepsilon^{-1}R(x_0,t_0)^{-\frac{1}{2}},
 -\min \{R(x_0,t_0)^{-1},t_0\})$ is, after
 scaling with the factor $R(x_0,t_0)$,  $\varepsilon$-close to the
 corresponding subset of the evolving round cylinder. Since the solution is rotationally symmetric around $O$,
  the cap neighborhood structures of those points $x_0$ with $d_{t_0}(x_0,O)\leq B_{1}(\varepsilon)$ follows
  directly. The Assertion 2 is proved.

  Therefor we finish the proof of
  Corollary A.2.
$$\eqno \#$$ 


\vskip 1cm

\begin{thebibliography}{99}
\bibitem{BGP} Burago,Y.  Gromov, M. and Perelman, G., {\sl A. D. Alexandrov spaces with curvatures bounded
below}, Russian Math. Surveys {\bf 47} (1992), 1-58.
\bibitem{Cao85} Cao, H.-D., {\sl Deformation of K$\ddot{a}$hler metrics to K$\ddot{a}$hler-Einstein
metrics on compact K$\ddot{a}$hler manifolds}, Invent. Math.,
\textbf{81} (1985), no. 2, 359--372.
\bibitem{CaZ}  Cao, H. D. and Zhu, X. P.,  {\sl A complete proof of the Poincar\'e and geometrization conjectures
-- application of the Hamilton-Perelman theory of the Ricci flow},
Asian J. Math., \textbf{10} (2006), no. 2, 165-492.
\bibitem{C} Chow, B., {\sl The Ricci flow on 2-sphere}, J. Diff. Geom. {\bf
33} (1991), 325-334.
\bibitem{CC} Cheeger, J. and Colding, T. H., {\sl On the structure of the
spaces with Ricci curvature bounded below I.}
 J. Diff. Geom. {\bf 46} (1997), 406-480.
\bibitem{CE}Cheeger, J. and Ebin, D., {\sl Comparison theorems in
Riemannian geometry,} North-Holland (1975).
\bibitem{CG}Cheeger, J. and Gromoll, D., {\sl On the structure of complete
manifolds of nonnegative curvature,} Ann. of Math., \textbf{46}
(1972), 413-433.

\bibitem{CTZ}  Chen, B. L., Tang, S. H. and Zhu, X. P., {\sl A
uniformization theorem of complete noncompact K\"ahler surfaces with
positive bisectional curvature}, J. Diff. Geom. {\bf 67 }
(2004),519-570.
\bibitem{CZ}  Chen, B. L. and Zhu, X. P., {\sl On complete noncompact K\"ahler
manifolds with positive bisectional curvature}, Math. Ann. {\bf 327}
(2003), 1-23.
\bibitem{CZ2}  Chen, B. L. and Zhu, X. P., {\sl Uniqueness of the
Ricci flow on complete noncompact manifolds}, arXiv:math.
DG/0505447 v3 May 2005, preprint.
\bibitem{De} De Turck, D., {\sl Deforming metrics in the direction of their Ricci
tensors} J. Diff. Geom. {\bf 18} (1983), 157-162.
\bibitem{Di} Ding, Y., {\sl Notes on Perelman's second paper} http://www.math.
lsa.umich.edu/research/ricciflow/perelman.html
\bibitem{Ha1}  Hamilton, R. S., {\sl Three manifolds with positive
Ricci curvature }, J. Diff. Geom. {\bf 17} (1982), 255-306.
\bibitem{Ha2}  Hamilton, R. S., {\sl Four--manifolds with positive
curvature operator}, J. Diff. Geom. {\bf 24} (1986), 153-179.
\bibitem{Ha3} Hamilton, R. S., {\sl The Ricci flow on surfaces,}
 Contemporary Mathematics {\bf 71} (1988) 237-261.
\bibitem{Ha4}  Hamilton, R. S., {\sl The Harnack estimate for the Ricci flow},
J. Diff. Geom. {\bf 37} (1993), 225-243.
\bibitem{Ha5}  Hamilton, R. S., {\sl A compactness property for solution
of the Ricci flow}, Amer. J. Math. {\bf 117} (1995), 545-572.

\bibitem{Ha6}  Hamilton, R. S., {\sl The formation of singularities in
the Ricci flow}, Surveys in Diff. Geom. (Cambridge, MA, 1993),
{\bf 2}, 7-136, International Press, Combridge, MA,1995.
\bibitem{Ha7} Hamilton, R. S., {\sl Four manifolds with positive isotropic
curvature}, Comm. Anal. Geom.,{\bf 5}(1997),1-92. (or see, {\sl
Collected Papers on Ricci Flow}, Edited by H. D. Cao, B. Chow, S. C.
Chu and S. T. Yau, International Press 2002).
\bibitem{Ha8} Hamilton, R. S., {\sl Non-singular solutions to the Ricci flow on three manifolds},
Comm. Anal. Geom. {\bf 1} (1999), 695-729.
\bibitem{Hir} Hirsch, M. W., {\sl Differential Topology},
Springer-Verlag, 1976.
\bibitem{Hu} Huisken, G., {\sl Ricci deformation of the metric on a Riemanian mnifold}
J. Diff. Geom. {\bf 21 }(1985), 47-62.
\bibitem{KL} Kleiner, B. and Lott, J., {\sl Note on Perelman's
paper,}
http://www.math.lsa.umich.edu/research/ricciflow/perelman.html.
\bibitem{MiMo} Micallef, M. and Moore, J. D., {\sl Minimal two-spheres and the
topology of manifolds with positive curvature on totally isotropic
two-planes}, Ann. of Math. (2) {\bf 127} (1988), 199-227.
\bibitem{Mi} Milka, A. D., {\sl Metric structure of some class of spaces containing straight
lines}, Ukrain. Geometrical. Sbornik, \textbf{4}, 1967, 43-48.
\bibitem{Mor} Morgan, J. W., {\sl Recent progress on the
Poincar$\acute{e}$
 conjecture and the classification of 3-manifolds,} Bull. of the
 A. M. S.,  {\bf 42} (2004) no. 1, 57-78.
 \bibitem{P1} Perelman, G., {\sl The entropy formula for the Ricci flow and its geometric
applications}, arXiv:math.DG/0211159 v1 November 11, 2002, preprint.
\bibitem{P2} Perelman, G., {\sl Ricci flow with surgery on three
manifolds} arXiv:math.DG/0303109 v1 March 10, 2003, preprint.
\bibitem{P3} Perelman, G., {\sl Finite extinction time to the solutions to the Ricci flow on certain three
manifolds, } arXiv: math. DG/0307245 July 17, 2003, preprint.
\bibitem{ScY}  Schoen, R. and Yau, S. T., {\sl Lectures on differential
geometry}, in conference proceedings and Lecture Notes in Geometry
and Topology, Volume {\bf 1}, International Press Publications,
1994.
\bibitem{STW} Sesum, N., Tian, G. and Wang, X. D., {\sl Notes on Perelman's
paper on the entropy formula for the
Ricci flow and its geometric applications.}
\bibitem{Sh1}  Shi, W. X., {\sl Deforming the metric on complete
Riemannian manifold}, J. Diff. Geom., {\bf 30} (1989), 223-301.

\end{thebibliography}
\end{document}